\documentclass{article}
\usepackage{rfm}
\usepackage[margin=1.5in]{geometry}

\usepackage{amsmath}
\usepackage{amssymb}
\usepackage{comment}
\usepackage{subfigure}
\usepackage{bm,bbm}
\usepackage{url}	
\usepackage{graphicx}
\usepackage{epstopdf}
\usepackage{epsfig,psfrag}
\usepackage{verbatim} 
\usepackage{subfigure}
\usepackage{algorithmic}
\usepackage{enumerate}
\usepackage{xcolor}
\usepackage{pdfsync}
\usepackage{hyperref}
\hypersetup{colorlinks,citecolor=blue,urlcolor=blue,linkcolor=blue}

\usepackage{natbib}
 \bibpunct[, ]{(}{)}{,}{a}{}{,}%
 \def\bibfont{\small}%
 \def\bibsep{\smallskipamount}%
 \def\bibhang{24pt}%
\newtheorem{theorem}{Theorem}[section]
	
\newtheorem{lemma}[theorem]{Lemma}

\newtheorem{proposition}[theorem]{Proposition}

\newtheorem{definition}[theorem]{Definition}

\newtheorem{remark}[theorem]{Remark}

\begin{document}

\title{Reinforcement with Fading Memories}


\author{
Kuang Xu\\
Graduate School of Business\\
Stanford University\\
\texttt{kuangxu@stanford.edu} 
\and
Se-Young Yun\\
Industrial \& Systems Engineering\\
KAIST\\
\texttt{yunseyoung@kaist.ac.kr}
}
\date{}

\maketitle 

\begin{abstract}

We study the effect of imperfect memory on decision making in the context of a stochastic sequential action-reward problem. An agent chooses a sequence of actions which generate discrete rewards at different rates. She is allowed to make new choices at rate $\beta$, while past rewards disappear from her memory at rate $\mu$. We focus on a family of decision rules where the agent makes a new choice by randomly selecting an action with a probability approximately proportional to the amount of past rewards associated with each action in her memory. 

We provide closed-form formulae for the agent's steady-state choice distribution in the regime where the memory span is large ($\mu \to 0$), and show that the agent's success critically depends on how quickly she updates her choices relative to the speed of memory decay. If $\beta \gg \mu$, the agent almost always chooses the best action, i.e., the one with the highest reward rate. Conversely,  if $\beta \ll \mu$, the agent chooses an action with a probability roughly proportional to its reward rate.\footnote{September 2017; revised: September 2019.  An extended abstract of this paper appeared in ACM SIGMETRICS, 2018.}


\end{abstract}

\maketitle 

\section{Introduction}
\label{sec:intro}

Memories serve as a crucial link between the past and future in dynamic decision-making problems, allowing subsequent decisions to draw on earlier experiences. However, they are hardly perfect in reality: humans routinely rely on faulty memories when making choices, and a firm's time-varying states, such as an active user base that is prone to unpredictable attrition, could sway its strategic focus. This raises the question: to what extent can we make good decisions despite having an imperfect memory? In this paper, we investigate this question in the context of a stochastic action-reward model. While the subject can be approached from many different angles, we will focus on understanding the interplay between the rate of memory decay and the rate of decision updates; the former a measure of the quality of memory, and the latter a proxy of the decision maker's temporal adaptivity. Our main results demonstrate that the relative magnitude between the two rates can have profound performance implications. 

\begin{figure}[h]
\vspace{-10pt}
\centering
\includegraphics[scale=.9]{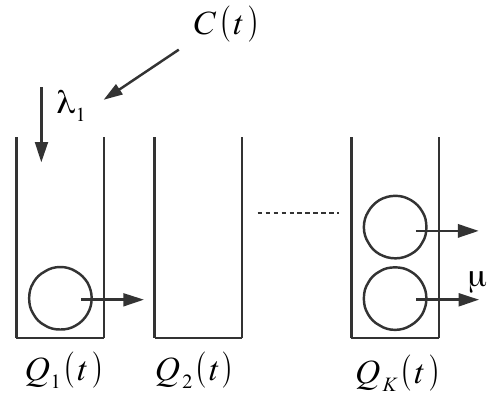}
\caption{A snapshot of the sequential decision-making problem, where the choice at time $t$ is action $1$, and the agent accrues discrete rewards according to Poisson process of rate $\lambda_1$. Each existing recallable reward departs from the system at rate $\mu$, as depicted by the horizontal arrows.}
\label{fig:scheme}
\end{figure}

Let us begin with an informal description of our model (Figure \ref{fig:scheme}). An agent operating in continuous time makes choices among a finite menu of actions. When action $k$ is chosen, the agent accrues discrete rewards according to a Poisson process with rate $\lambda_k$, while	 a unit of reward ``expires'' and disappears from her memory after a random amount of time that is exponentially distributed with rate $\mu$. Opportunities for the agent to make new choices arise according to a Poisson process with rate $\beta$. We will restrict our attention to a simple family of decision rules that the agent may use, referred to as the {\bf reward matching} rule, which {reinforces} actions with more positive past stimuli in proportional manner: she randomly samples an action such that the probability of choosing action $k$ is proportional to $\max\{Q_k, \alpha\}$, where $Q_k$ is the total units of rewards accrued through action $k$ that have yet to expire, and $\alpha>0$ an exploration parameter representing the agent's minimal willingness to experiment with any action. If $\alpha$ is set to 0, then the reward matching rule corresponds to the celebrated Luce's linear probabilistic choice model \citep{luce1959individual,erev1998predicting}, which has been extensively studied in behavioral economics \citep{erev1998predicting,beggs2005convergence} and evolutionary biology \citep{harley1981learning}. We are interested in the probability that the agent chooses a particular action in steady state.  We next examine two illustrative examples, which also serve as motivating applications for the more stylized model we study.

{\emph{Example 1 - Consumer Choice Modeling}. We may think of the agent as a consumer, the actions as products or service providers, and the opportunities to choose new actions as the times when the consumer is allowed to re-select her membership or subscription. The rewards act as positive experiences or impressions resulting from using a service, and the strength of each impression in the agent's memory diminishes as time progresses. When an opportunity to renew the subscription arises, the agent chooses a new service offering in a manner that is biased towards  the ones associated with more recallable impressions. }

{\emph{Example 2 - Dynamic Product Offering under Customer Attrition.} In this example, we slightly refine the model by letting the departure rates of rewards vary across actions. Consider a firm who operates a number of service contracts (actions), e.g., mobile phone plans, and needs to decide which contract to offer new customers (rewards) in a series of promotion periods, where one contract is offered for each period. Customers  arrive to the system and are subscribe to the service contract of the corresponding promotion period. Those who subscribe to contract type $k$ remain an active customer for an exponentially distributed amount of time with rate $\mu_k$, before departing from the system. The departure rate $\mu_k$ reflects some underlying qualities of the contract $k$. The reward matching rule corresponds to the service provider's choosing, at the end of each promotion period, a new contract with probabilities roughly proportional to the number of {active} customers  associated with each contract. 
}

\paragraph{Motivation for studying the reward matching rule.} An implied objective of our model is that the agent would like to maximize the steady-state total amount of rewards in system.  {If} all the parameters were known \emph{a priori}, it is not difficult to see that this objective could be trivially achieved by choosing at all times the {best} action, i.e., one with the highest reward rate $\lambda_k$ (Example 1) or the highest weighted reward rate, $\lambda_k/\mu_k$ (Example 2). The reward matching rule, therefore, can be viewed as an intuitive heuristic for approximately solving this optimization problem in the absence of the knowledge of the parameters. While the rule is by no means the only plausible heuristic for this purpose,  we believe that it has a number of advantages which merit a theoretical investigation such as the one carried out in this work:
\begin{enumerate}
\item Firstly, Luce's rule, after which the reward matching rule is modeled, has a long history as a fundamental building block in psychology and behavioral economics (see Section \ref{sec:litRev} for more details). The reward matching rule hence serves a natural starting point for  understanding other, potentially more complex, cognitive models. 
\item Secondly, the reward matching rule is useful not only for modeling human choice behavior. It could also be a policy deliberately chosen by a decision maker, such as in Example 2. While a conscious decision maker could in principle adopt more sophisticated decision rules or even use external information storage to mitigate memory decay, the reward matching rule offers a simple and intuitive rule of thumb that can be substantially easier for a practitioner to understand and implement, and is applicable in settings where external information storage is not possible. 
\item Last but not least, our results show that, under a large number of parameter values, the reward matching rule is in fact \emph{near optimal} and induces a steady-state choice probabilities that heavily concentrate on the best action. Along with  the  simplicity and conceptual appeal, the near-optimality therefore further justifies studying the reward-matching rule. 
\end{enumerate}

\paragraph{Preview of main results.} Our main results provide simple, closed-form formulae for the distribution of the agent's choice {in steady state}, in the regime where the agent's memory span is large ($\mu\to 0$). Using these formulae, we show that the agent's success critically depends on how fast her choices are updated ($\beta$) relative to the rate of memory decay ($\mu$), when $\alpha$ is relatively small compared to the overall memory size. In particular, if $\beta \gg \mu$,\footnote{The notation $\beta \gg \mu$ represents $\lim_{\mu \to 0} {\beta}/{\mu} = 0$.} the agent nearly always chooses the {best action}, i.e., the one with the highest reward rate $\lambda_k$. In contrast, if $\beta \ll \mu$,  the agent splits her attention more smoothly among all actions, with frequencies roughly proportional to the respective $\lambda_k$'s. 

\subsection{The Model} 

We now describe our model more formally, as is  depicted in Figure \ref{fig:scheme}. The system consists of an agent operating in continuous time, who  makes \emph{choices} by selecting from a set of $K$ \emph{actions}, $\calK=\{1, \ldots, K\}$. We denote by $C(t)\in \calK$ the action chosen by the agent at time $t$, and refer to $\{C(t)\}_{t\in \rp}$ as the \emph{choice process}. When the $C(t) = k$ for some $k\in \calK$, the agent accrues {discrete rewards}  according to a Poisson process with rate $\lambda_k>0$, where each reward can be treated as a token of unit size. We refer to $\lambda_k$ as the \emph{reward rate} of action $k$. We assume that there exists one reward rate that strictly dominates the rest, and, without loss of generality, we rank them in a decreasing order, so that $\lambda_1 > \lambda_2 \geq  \ldots \geq  \lambda_K > 0$. 
 
The agent's fading memory is modeled by the ``expiration'' of past rewards. Each reward is associated with a \emph{lifespan}, an exponential random variable with rate $\mu>0$ drawn independently from all other aspects of the system, and a reward {departs} permanently after staying in the system for a time duration equal to its lifespan. We refer to  $\mu$ as the \emph{rate of memory decay}.  {For simplicity of notation, the majority of our exposition and proofs will focus on the case where $\mu$ is uniform across all actions, but our main results extend to the case with action-dependent memory decay rates in a straightforward manner; such a refinement will be discussed in Section \ref{sec:heterorate}.}

For $k\in \calK$, we denote by $Q_k(t)$ the total units of rewards accrued from choosing action $k$ that have yet to depart by time $t$, and we refer to $\{Q(t)\}_{t\in \rp}$ as the \emph{recallable reward process}. Note that at time $t$ the recallable rewards associated with action $k$ depart at an aggregate rate of $\mu Q_k(t)$, and arrive at rate $\lambda_k$, if $C(t)=k$, or $0$, otherwise. We assume that the agent's initial recallable rewards at time $0$ is a bounded vector in $\zp^\calK$.  

The agent has the opportunity to make a new choice at a set of \emph{update points} scattered across time. The update points are distributed according to an exogenous Poisson process with rate $\beta$, which we refer to as the \emph{update rate}, so that the time between two adjacent update points is an exponentially distributed random variable with mean $\beta^{-1}$, independent from all other aspects of the system. Depending on the application context, $\beta$ could be interpreted as the level of ``activeness'' or ``savvyness,'' as it modulates the frequency at which the agent updates her decisions according to the present state of her recallable-memories (cf.~Examples 1 and 2 in the Introduction).

How does the agent make a new choice when an update point arises? In this work, we will focus on a family of choice heuristics referred to as the \emph{reward matching} rules. We assume that at time $t$ the agent makes a new choice by sampling from the distribution
\begin{equation}
\pb( \mbox{new choice } = k) = \frac{Q_k(t) \vee \alpha }{\sum_{i \in \calK} (Q_i(t) \vee \alpha) },  \quad k \in \calK, 
\label{eq:RMsampling}
\end{equation}
where $x \vee y \bydef \max\{x, y\}$.  Here, $\alpha$ is a positive constant referred to as the \emph{exploration parameter}, 
which captures the agent's willingness to experiment with an action even if there is currently  little or no recallable reward attached to it. The justification and background behind the choice of this decision rule will be further elaborated on in  Section \ref{sec:litRev}. The parameter $\alpha$ can also have an effect on the  trade-off between the system's steady-state versus transient performance; the reader is referred to Section \ref{sec:genChoi} for a discussion.


\begin{figure}
\centering
\includegraphics[scale=.55]{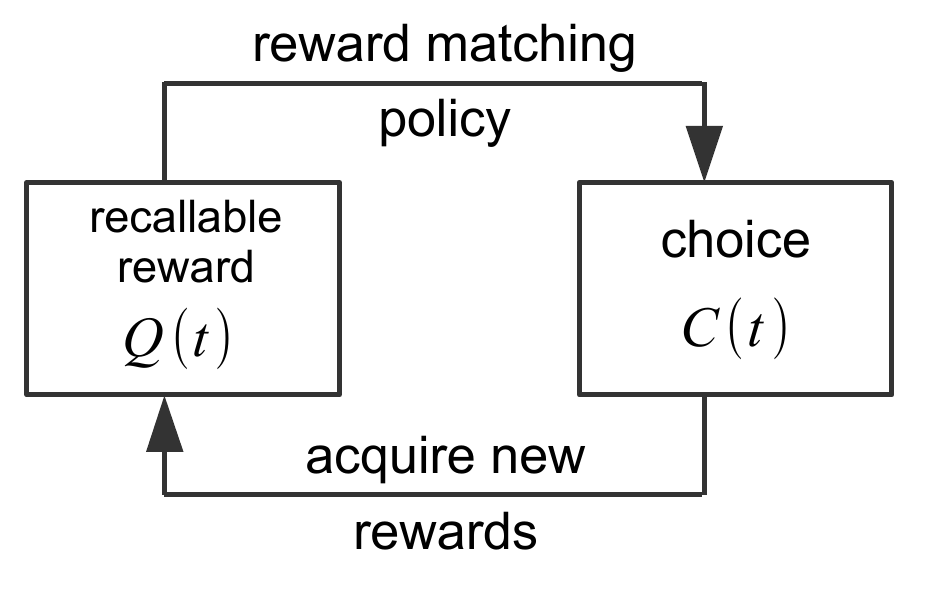}
\caption{A high-level illustration of the interaction between the choice and recallable reward processes.}
\label{fig:QandC}
\end{figure}

\subsection{Performance Metrics and Scaling Regime}
\label{sec:metrics}

The main goal of this paper is to study the agent's long-run behavior in the model described above, and a key quantity of interest is the steady-state distribution of the choice process, $C(\cdot)$, i.e., 
\begin{equation}
\lim_{t\to \infty}\pb(C(t) = k), \quad k \in \calK, 
\end{equation}
Unfortunately, the choice process is non-Markovian as its transitions are influenced by the recallable reward process $Q(\cdot)$ in a non-trivial way, and this makes it difficult to obtain explicit expressions for its steady-state distribution.  To circumvent this challenge, we will focus on the regime where the average lifespan of a  reward, $\mu^{-1}$, tends to infinity, as follows. Fix the action set, $\calK$, and the reward rates, $\{\lambda_k\}_{k \in \calK}$. We consider a sequence of systems, indexed by $m\in \N$, where the average lifespan of a reward in the $m$th system is equal to $m$, i.e.,
\begin{equation}
 \mu = 1/m.
 \end{equation} 
 We will denote by $\alpha_m$, $\beta_m$, $Q^m(\cdot)$, and $C^m(\cdot)$ the corresponding quantities in the $m$th system. Finally, as the rewards' average lifespan, $m$, tends to infinity, the total amount of recallable rewards in the system also scales as $\Theta(m)$.\footnote{To see this, note that if the choice process were to stay in action $1$ all the time, then $Q^m_1(\cdot)$ would evolve according to an $M/M/\infty$ queue with arrival rate $\lambda_1$ and departure rate $m^{-1}$, leading to a steady-state expected total recallable reward of $\lambda_1 m $.} For this reason, we will scale the exploration parameter $\alpha_m$ linearly with respect to $m$, by setting $\alpha_m = \alpha_0*m$, where $\alpha_0$ is a fixed parameter.

{\emph{Remark.} It is natural to ask  whether by taking a sequence of systems where the memory decay rate approaches zero, we will end up with the trivial case where there is no memory loss at all ($\mu=0$), which would have been at odds with our goal of studying imperfect memories. This not the case. Because we always look at the steady state of the system by taking the limit as $t\to \infty$ for every system in the sequence, the limiting regime is non-trivial and is not equivalent to simply setting $\mu=0$, as is evidenced by the existence of two {distinct} limiting steady-state distributions in our main theorems (This point is further elaborated on in Appendix \ref{app:perfect_memory}.) The  limit of $\mu\to 0$ is intended as an approximation for a pre-limit system where $\mu$ is positive but small. Numerical results in Figure \ref{fig:SS-Sim} and \ref{fig:tranResult} suggest that our theorems  provide a reasonably good approximation for $m$ as small as $200$.}

\subsection{Notation}

For $x, y \in \R$,  $x \vee y$ and $x \wedge y$ denote $\max\{x,y\}$ and $\min\{x,y\}$, respectively.  For a vector $x\in \rp^K$,  $\|x\|$ denotes the $L_1$ norm of $x$, $\|x\| = \sum_{i=1}^K|x_i|$.  

We use the notation $X(\cdot)$ and $X[\cdot]$ as a shorthand for processes $\{X(t)\}_{t\in \rp}$ and $\{X[n]\}_{n \in \N}$, respectively.  For two random variables $X$ and $Y$ taking values in $\R$, we use $X \stackrel{d}{=} Y$ to mean that $X$ and $Y$ have the same distribution. The expression $X \preceq Y$ means $X$ is stochastically dominated by $Y$, i.e., $\pb(X\geq x)\leq \pb(Y\geq x)$ for all $x\in \R$. For an event $\calA$, we use $X | \calA$ to denote the random variable with the cumulative distribution function $\pb(X \leq x \bbar \calA )$, $x\in \R$.  For a stochastic process, $\{X(t)\}_{t\in \rp}$, we denote by $X(\infty)$ the random variable distributed according to the stationary distribution of $\{X(t)\}_{t\in \rp}$, if it exists. 

For two functions $f$ and $g$, $f \gg g$ means that $\lim_{x \to \infty}\frac{f(x)}{g(x)} =\infty$; $f\ll g $ is defined analogously. When such a limiting notation is involved, the limit is assumed to be as the argument tends to $\infty$, unless explicitly specified otherwise.

\section{Main Results} 

We now present our main results, which provide exact expressions for the steady-state distribution of the choice process, as well as the expected value of the (scaled) recallable reward process, in the limit as $m\to \infty$. For the remainder of the paper, we will assume that the update rate $\beta_m$ satisfies either $\beta_m \gg 1/m$ or $\beta_m \ll 1/m$, as $m\to \infty$. 

\begin{theorem}[Steady-State Choice Probabilities]
\label{thm:Ck}
Fix $\alpha_0>0$.  The following limits exist
\begin{equation}
c_k = \lim_{m \to \infty}\pb(C^m(\infty) = k), \quad k \in \calK. 
\end{equation}
Furthermore, their expressions depend on the scaling of $\beta_m$, as follows. 
\begin{enumerate}
\item \emph{(Memory-abundant regime)} Suppose that $\beta_m \gg 1/m$, as $m\to \infty$. If $\alpha_0 \leq \lambda_1/K$, then
\begin{align}
c_1 =& 1- \frac{K-1}{\lambda_1}\alpha_0,  \nln
c_k =& \frac{\alpha_0}{\lambda_1} ,  \quad   k = 2, \ldots, K. 
\label{eq:ckFast1}
\end{align}
If $\alpha_0 > \lambda_1/K$, then
\begin{equation}
c_k = {1}/{K}, \quad k= 1,\ldots, K. 
\label{eq:ckFast2}
\end{equation}
\item \emph{(Memory-deficient regime)}  Suppose that $\beta_m \ll 1/m$, as $m\to \infty$. Then,
\begin{equation}
c_k =\frac{\lambda_k\vee \alpha_0 + (K-1)\alpha_0}{\sum_{i \in \calK}\lt[ \lambda_i \vee \alpha_0 + (K-1)\alpha_0 \right]}, \quad k = 1, \ldots, K. 
\label{eq:ckSlow}
\end{equation}
\end{enumerate}
\end{theorem}

The next  theorem characterizes the steady-state expected value of the recallable reward process, $Q^m(\cdot)$.  For $m\in \N$, define the \emph{scaled recallable reward process}: 
\begin{equation}
\qol (t) = \frac{1}{m}Q^m (m t), \quad t\in \rp. 
\end{equation}

\begin{theorem}[Steady-State Recallable Rewards]
\label{thm:main}
Fix $\alpha_0>0$. For all $m\in \N$, the scaled recallable reward process, $\{\qol(t)\}_{t\in \rp}$, admits a unique steady-state distribution, and the following limits exist: 
\begin{equation}
q_k = \lim_{m \to \infty}\E(\qol_k(\infty)), \quad k \in \calK. 
\end{equation}
Furthermore, their expressions depend on the scaling of $\beta_m$, as follows. 
\begin{enumerate}
\item \emph{(Memory-abundant regime)} Suppose that $\beta_m \gg 1/m$, as $m\to \infty$. If $\alpha_0 \leq \lambda_1/K$, then
\begin{align}
q_1 =& \lambda_1 - (K-1)\alpha_0,  \nln
q_k =& (\lambda_k/{\lambda_1})\alpha_0 ,  \quad  k = 2, \ldots, K,
\label{eq:MainFast1}
\end{align}
If $\alpha_0 > \lambda_1/K$, then
\begin{equation}
q_k = {\lambda_k}/{K}, \quad k= 1,\ldots, K. 
\label{eq:MainFast2}
\end{equation}
\item \emph{(Memory-deficient regime)}  Suppose that $\beta_m \ll 1/m$, as $m\to \infty$. Then,
\begin{equation}
q_k = \lambda_{k}  \frac{\lambda_k\vee \alpha_0 + (K-1)\alpha_0}{\sum_{i \in \calK}\lt[ \lambda_i \vee \alpha_0 + (K-1)\alpha_0 \right]}, \quad k = 1, \ldots, K. 
\end{equation}

\end{enumerate}
\end{theorem}


\subsection{Implications of the Theorems}
\label{sec:impThm}
A prominent feature of Theorems \ref{thm:Ck} and \ref{thm:main} is that the results critically depend on how quickly the update rate $\beta_m$ scales relative to the rate of memory decay $1/m$, but are otherwise {insensitive} to the exact expression of $\beta_m$. This leads to a natural dichotomy of the system dynamics into what we called the memory-abundant and memory-deficient regimes, and  we discuss below several notable properties of the two regimes. 

\paragraph{Near optimality and winner-takes-all.} Suppose that the constant $\alpha_0$ in the exploration parameter is substantially smaller than $\lambda_1/K$. Then, our results show that the best action attracts nearly all of the agent's attention, if $\beta_m \ll 1/m$ (memory-abundant regime). Specifically, in the limit as $m\to \infty$, the agent chooses action $1$ with a probability of $(1-\frac{K-1}{\lambda_1}\alpha_0)$, which is almost $1$ when $\alpha_0$ is very small compared to $\lambda_1/K$. Recall from Section \ref{sec:intro} that the optimal strategy for the agent, should she know all parameters of the problem, would be to choose action $1$ at all times. Therefore, in this regime, the reward matching rule is able to deliver {near optimal} performance.

Moreover, the system exhibits an interesting winner-takes-all phenomenon, whereby the degree to which the agent focuses on the best action is {independent} of the  reward rates of all other actions. For instance, if we were to increase the reward rate of the second-best action, $\lambda_2$, the probability that the agent chooses action $1$ would stay unchanged (and close to 1), and that of choosing action $2$ would still be close to zero no matter how close $\lambda_2$ is to $\lambda_1$. One consequence of this effect is that, if we imagine the actions as being ``competitive'' service providers vying for a consumer's attention, then the winner-takes-all phenomenon can be interpreted as a form of {extreme competition} among the providers, where the best attracts nearly all of the consumer's budget,  even if its superiority compared to the second best is not significant.

\begin{figure}[h]
 \centering
 \includegraphics[scale=.35]{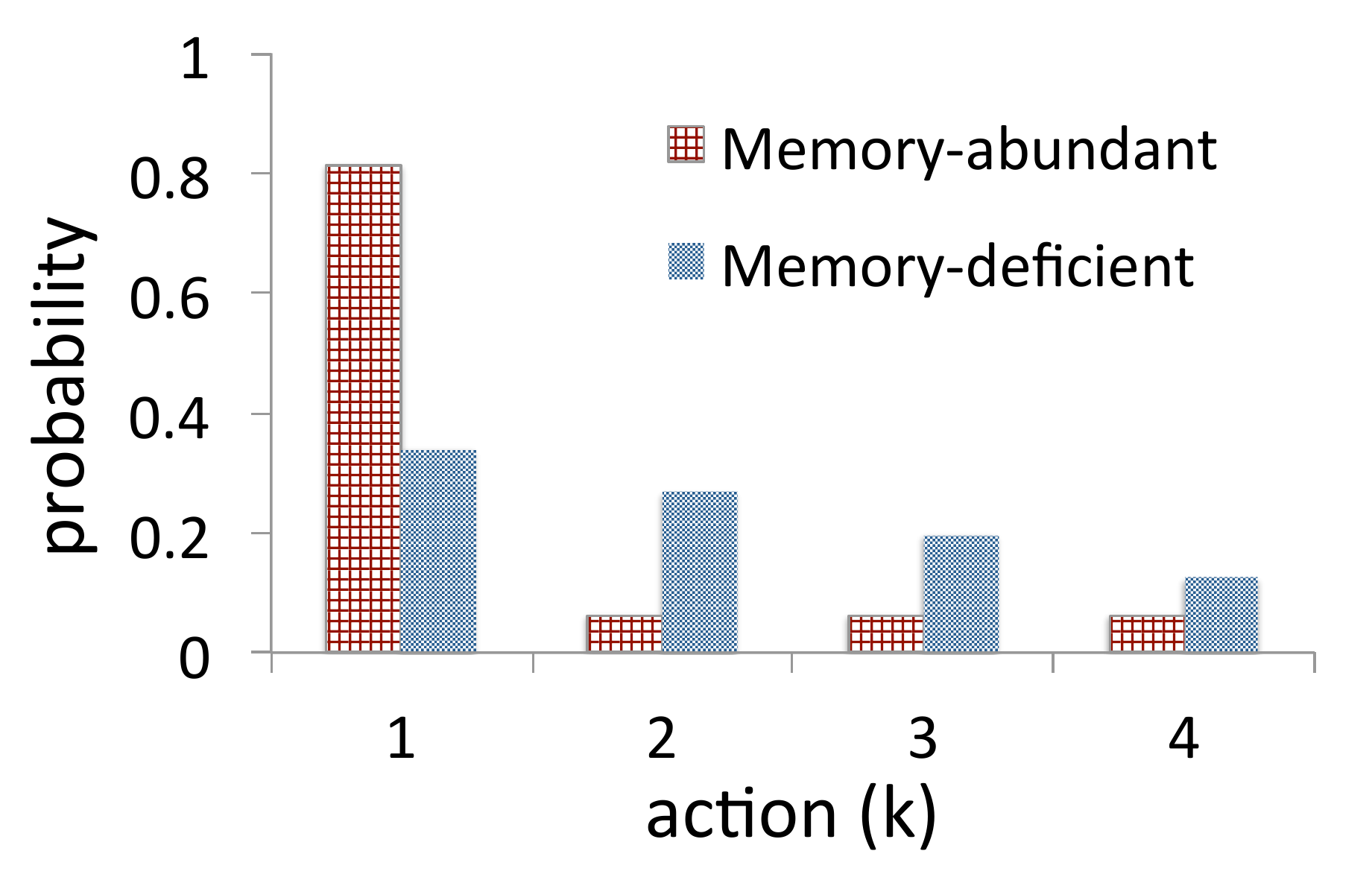}
 \caption{An example of the limiting steady-state distribution of the choice process $C^m(\cdot)$ in the two regimes given by Eqs.~\eqref{eq:ckFast1} through \eqref{eq:ckSlow}, with reward rates $(\lambda_1, \lambda_2, \lambda_3, \lambda_4)= (8, 6, 4, 2)$ and $\alpha_0=0.5$. The agent's  attention sharply concentrates on the best action in one regime, and is more spread out among all actions in the other.}
 \label{fig:ckSS}
\end{figure}

\paragraph{Reward-rate-proportional allocation.} In contrast to the memory-abundant regime, an update rate that scales substantially slower than the memory decay rate will lead to an allocation of choices that is much  smoother as a function of the reward rates (Figure \ref{fig:ckSS} provides a graphical comparison of the qualitative difference between the agent's behavior under the two regimes).  In this case, the probability that the agent chooses action $k$ is proportional to $[\lambda_k \vee \alpha_0 + (K-1)\alpha_0]$, which, when $\alpha_0$ is small, is essentially equal to the reward rate, $\lambda_k$. The agent's behavior thus exhibits a certain {weak reinforcement}, where the better actions do attract more attention but only proportional to their respective reward rates, and the above-mentioned winner-takes-all effect disappears. 


The intuition behind this reward-rate-proportional allocation can be roughly explained as follows. When the update rate is slow, by the time a new choice is to be made the agent would have forgotten the rewards associated with all actions other than her most recent choice.  Nevertheless, she still exhibits a  mild preference over actions with higher reward rates in steady state, because she {does} retain a good memory of the rewards of her most recent choice, $C(t)$, which can be shown to be proportional to $\lambda_{C(t)}$.  The reward rate $\lambda_{C(t)}$ dictates how likely the agent will make the same choices,  which in turn leads to the reward rates' approximately linear appearance in the steady-state distribution of the choice process. 

\paragraph{Rapid updates lead to complete oblivion.} The theorems also show that, perhaps surprisingly, a slower update rate is not always bad.  While it may be tempting to conclude that faster update rates always lead to better outcomes for the agent,  this turns out to be true only when the exploration parameter is relatively small. In fact, as soon as $\alpha_0$ grows beyond $\lambda_1/K$, in the memory-abundant regime, the agent becomes {completely oblivious} and chooses actions uniformly at random (Eq.~\eqref{eq:ckFast2}). In contrast, the agent is {always} biased towards the best action in the memory-deficient regime, albeit mildly, as long as $\alpha_0<\lambda_1$ (beyond which point the exploration parameter is so large that it trivially overwhelms even the best action). The true culprit behind the complete oblivion is a lack of patience: when the update rate is fast and the exploration parameter high, the agent tends to switch  her choices rapidly, before she is able to collect enough rewards to ``learn'' the reward rate of any individual action. Consequently,  in the long run, not a single component of the recallable reward, not even the best action, can distinguish itself by consistently rising above the exploration parameter, and hence complete oblivion becomes the only possible outcome.

\begin{remark}
Note that Theorems \ref{thm:Ck} and \ref{thm:main}   have left out the regime where $\beta$ and $\mu$ scale at the same rate, i.e., when  $\beta = \Theta(\mu)$ as $\mu\to 0$. This has turned out to be a difficult regime to tackle, and an intuitive explanation for this difficulty is as follows. Our current proof technique exploits the fact that one of the two processes in $Q(\cdot)$ and $C(\cdot)$ becomes asymptotic Markovian property as $\mu\to 0$, which significantly simplifies the dynamics (see Section \ref{sec:proof_over_view}). Unfortunately, this is no longer the case when $\beta=\Theta(\mu)$, as neither of the two processes $Q(\cdot)$ and $C(\cdot)$ becomes asymptotically Markovian by itself, and the stochastic dynamics of the system remain complex despite taking the limit of $\mu\to 0$. We may illustrate this phenomenon by the following back-of-the-envelope calculation: the probability that any individual unit of reward will depart from the system between two consecutive update points is on the order of $\beta\cdot \mu^{-1}$. When $\beta=\Theta(\mu)$, the number of reward departures between two update points forms a {constant} fraction of existing rewards and the same is true for the number of new rewards arriving during this period. Therefore, the changes in the recallable-reward process $Q(\cdot)$ remain highly variable from one update point to the next, and they depend both on the process $Q(\cdot)$ itself and the choice process, $C(\cdot)$.  In contrast, in the memory-abundant or memory-deficient regimes, the profile of the recallable rewards either changes very little between two update points (memory-abundant regime), or substantially but in a predictable manner (memory-deficient regime). 
\end{remark}

\subsection{Numerical Results}

\begin{figure}[h]
\begin{minipage}{0.46\textwidth}
\centering
\includegraphics[scale=.31]{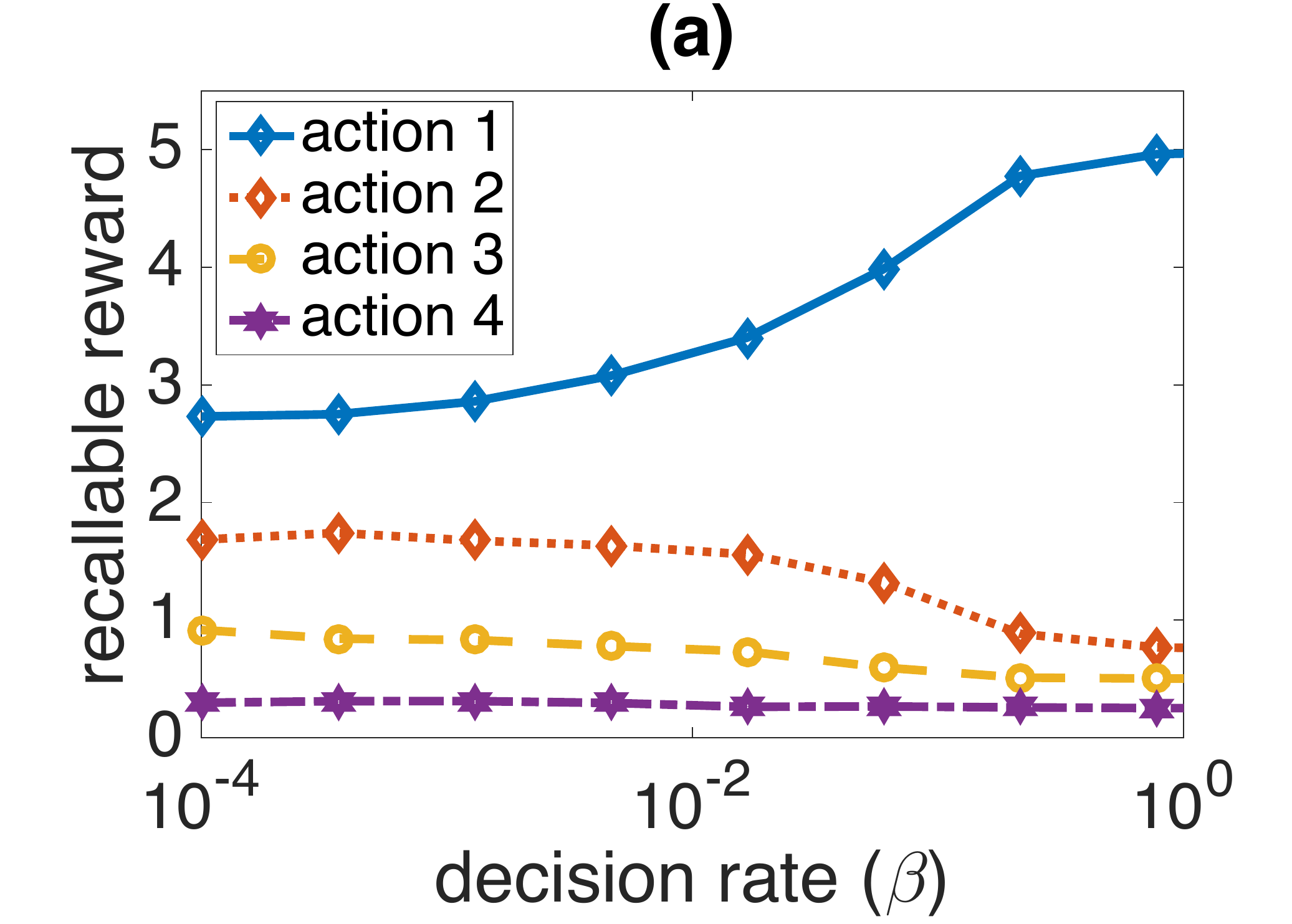}
\end{minipage}\hfill
\begin{minipage}{0.46\textwidth}
\centering
\includegraphics[scale=.29]{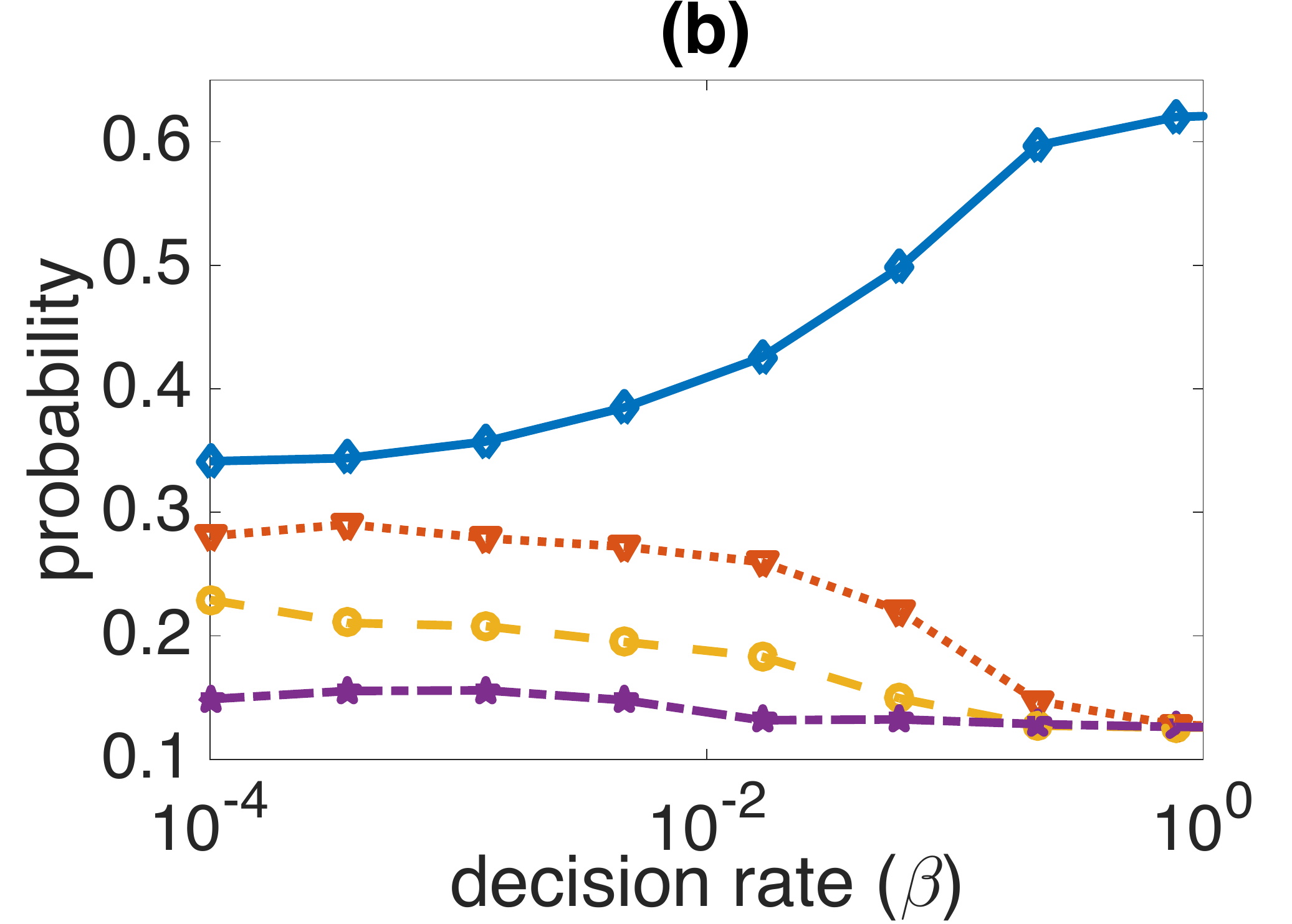}
\end{minipage}
\caption{Time-average scaled recallable rewards (plot $(a)$) and distribution of choices (plot $(b)$) as a function of the update rate, $\beta$, with $K=4$, $(\lambda_1,\lambda_2, \lambda_3, \lambda_4) =(8, 6, 4, 2) $, $\alpha_0=1$, and $m=200$. Each point is averaged over $5 \times 10^7$ to $5\times 10^5$ time units for $\beta$ ranging from $10^{-4}$ to $1$, respectively (longer duration for a slower update rate).}
\label{fig:SS-Sim}
\end{figure}

While Theorems \ref{thm:Ck} and \ref{thm:main}  only apply in the limit where the average reward lifespan tends to infinity, the simulation results in Figures  \ref{fig:SS-Sim} and \ref{fig:tranResult} show that our theoretical results also provide fairly accurate predictions, quantitatively and qualitatively, in systems with a moderate, finite lifespan. In Figure \ref{fig:SS-Sim}, we fix $m$, while varying
 the update rate, $\beta$, and we observe  that concentration on the best action intensifies as $\beta$ increases. Figure \ref{fig:tranResult} shows example sample paths of the scaled recallable reward processes, $\qol(\cdot)$. One can discern without difficulty that, even for a modest $m$, the process in the memory-abundant regime (plots $(a)$ through $(c)$) stays close to a set of smooth ``fluid'' trajectories, whose invariant state (the flat portion of the trajectories) coincide with the limiting value of $\E(\qol(\infty))$ predicted in Eq.~\eqref{eq:MainFast1}, while in the memory-deficient regime the process  changes more abruptly over time. As will become clear in the sequel, this qualitative discrepancy is a direct consequence of the overall system dynamics being dominated by the evolution of {either} $\qol(\cdot)$ or $C^m(\cdot)$, depending on the scaling of the update rate $\beta_m$, a key fact that we will exploit in our proof of the main theorems. 

\begin{figure}[h]
\centering
\begin{minipage}{0.47\textwidth}
\centering
\includegraphics[scale=.34]{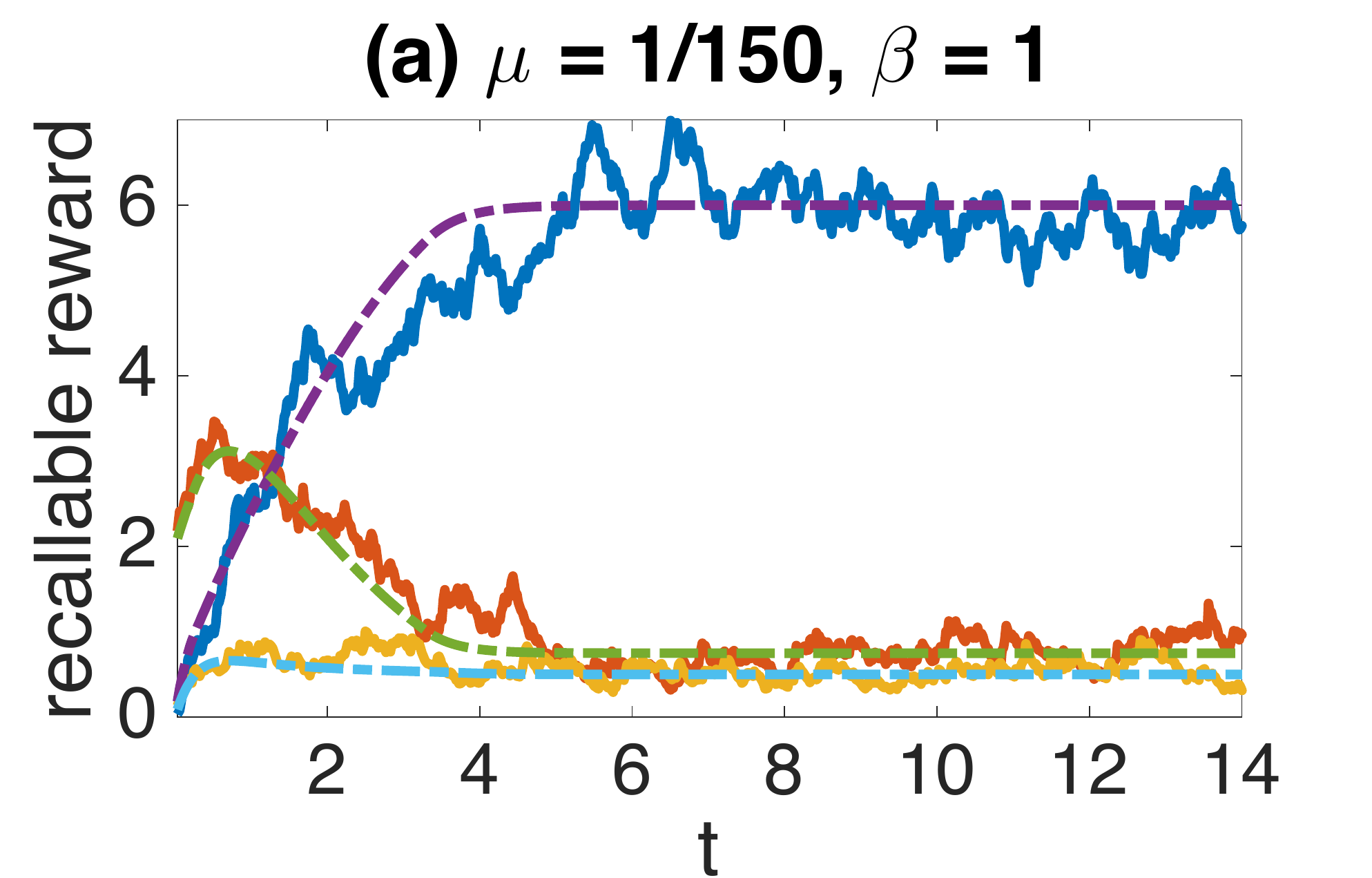}
\end{minipage}\hfill
\begin{minipage}{0.47\textwidth}
\centering
\includegraphics[scale=.32]{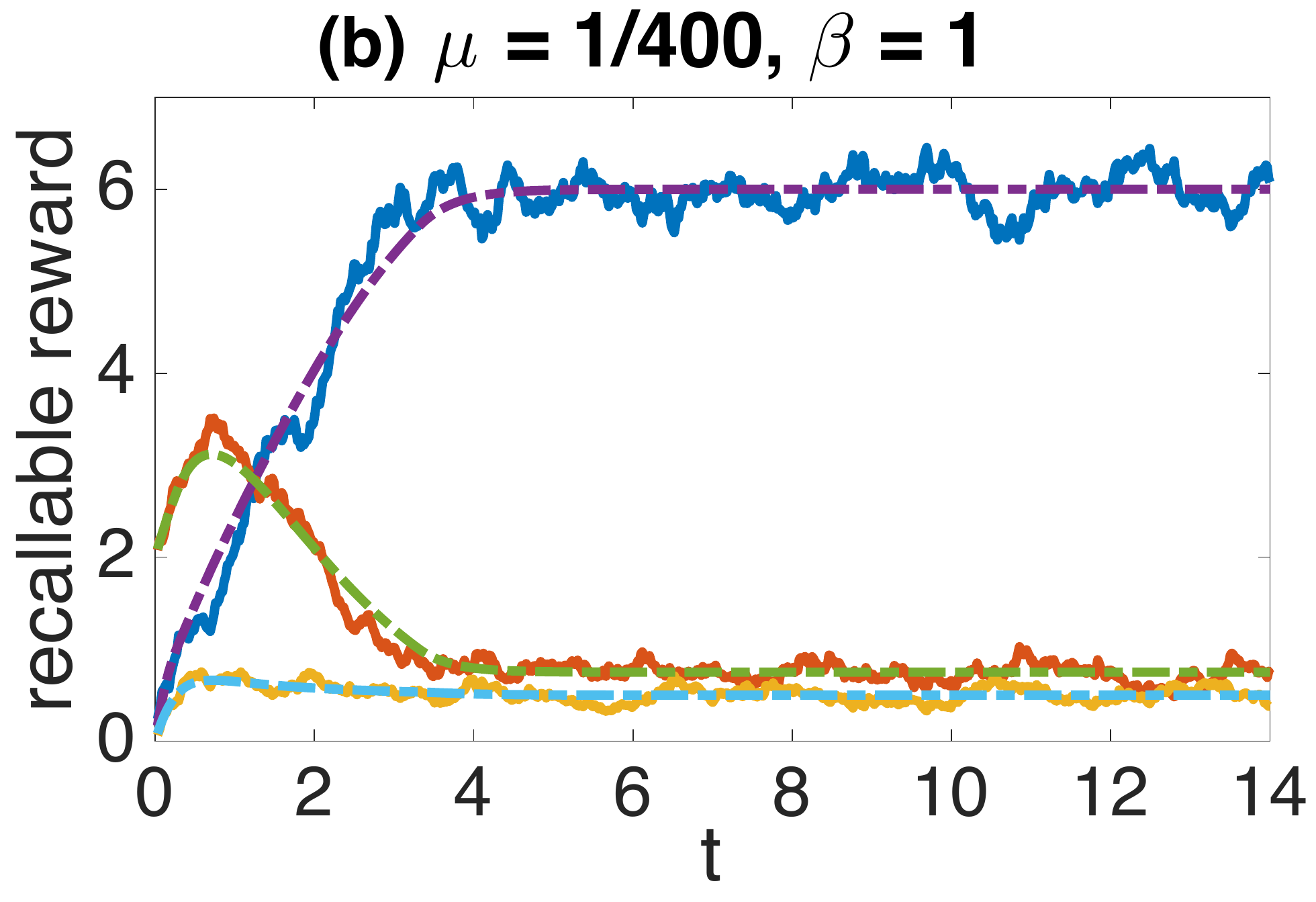}
\end{minipage}
\label{fig:tranResult}
\\
\vspace{10pt}

\begin{minipage}{0.47\textwidth}
\centering
\includegraphics[scale=.34]{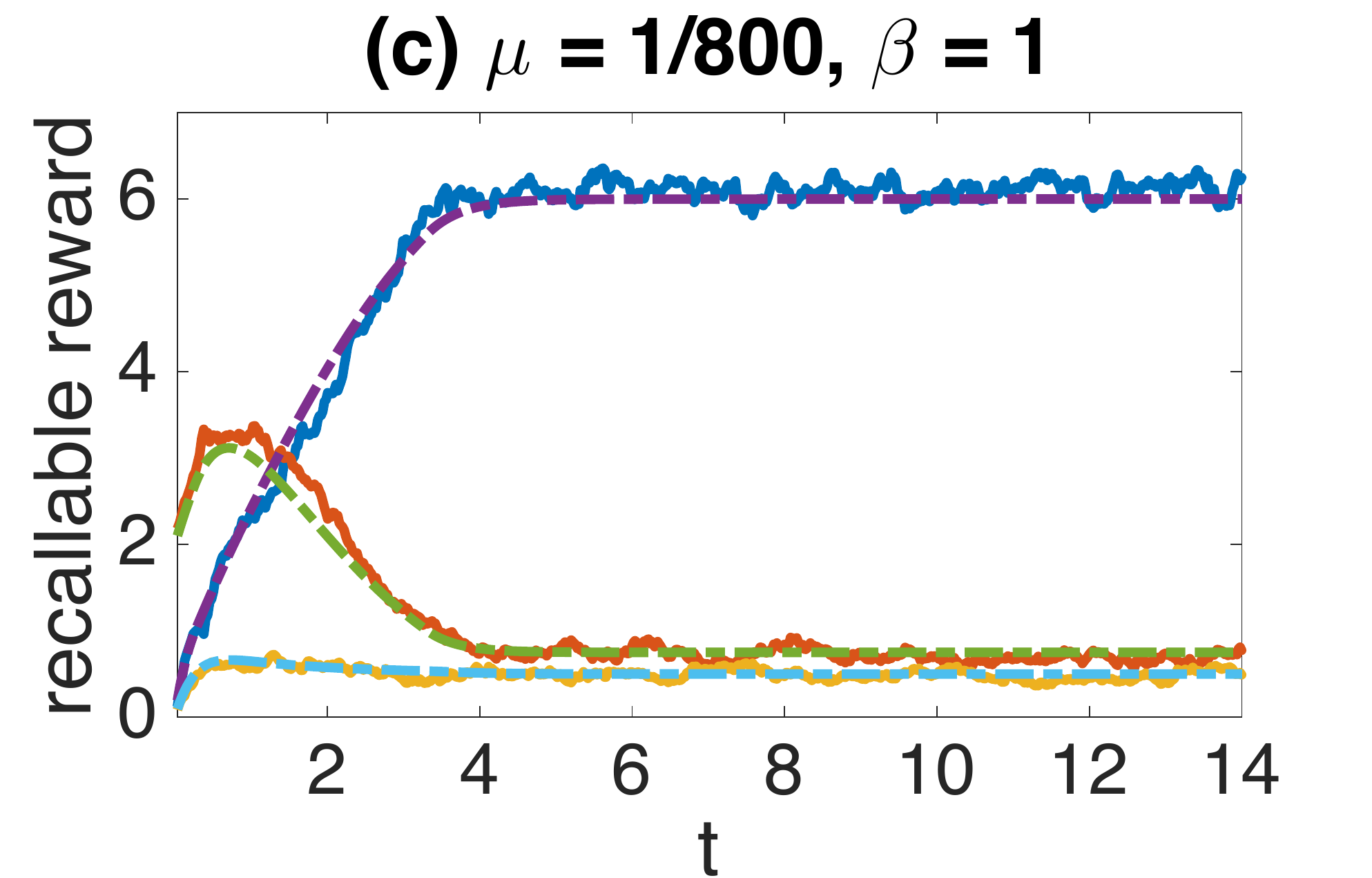}
\end{minipage}\hfill
\begin{minipage}{0.47\textwidth}
\includegraphics[scale=.295]{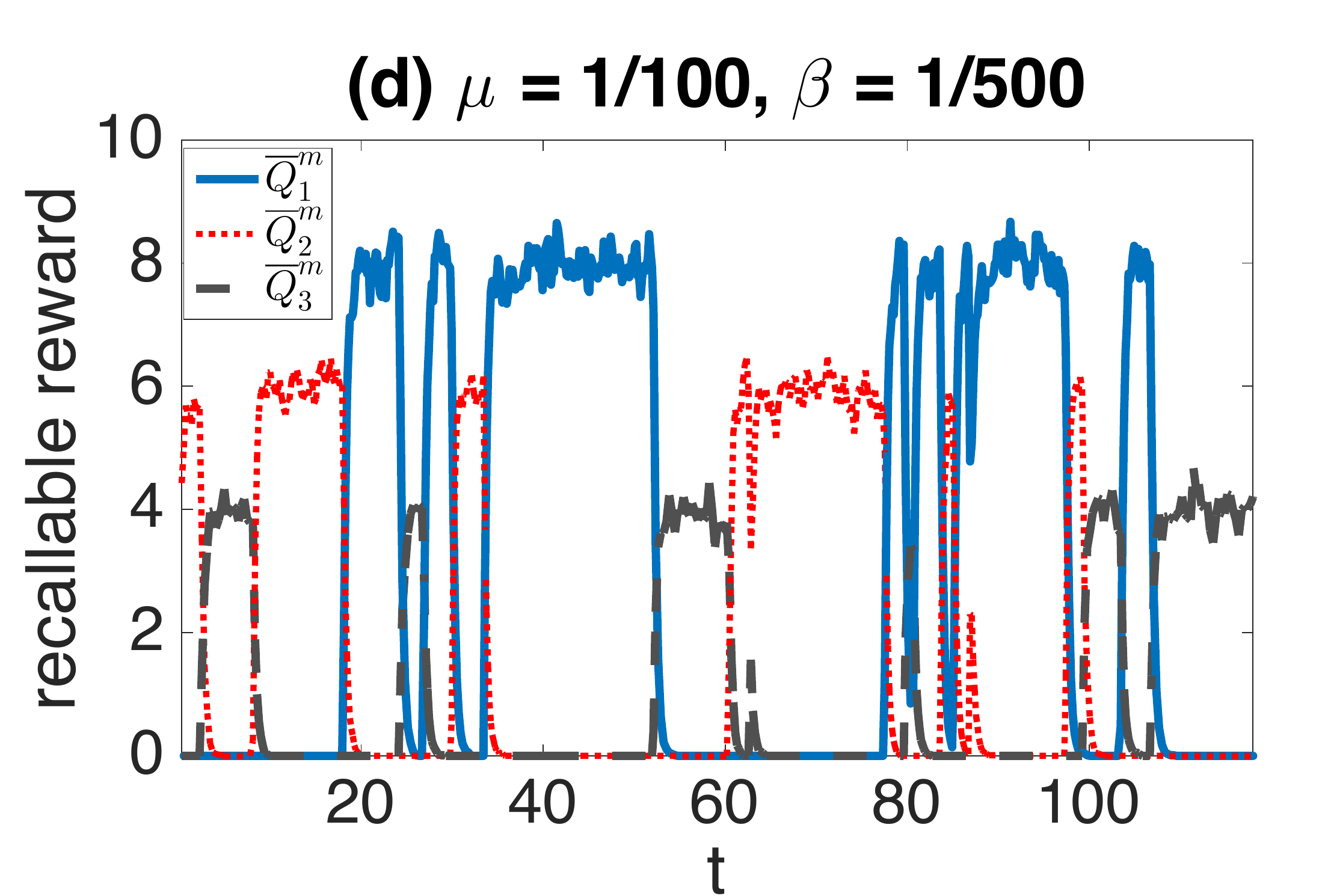}
\end{minipage}
\caption{Example sample paths under different memory decay and update rates, with $K=3$, $(\lambda_1,\lambda_2, \lambda_3) =( 8, 6, 4) $, $\alpha_0=1$. Plots $(a)$ through $(c)$ correspond to the memory-abundant regime with $m=150, 400$ and $800$, respectively, where the solid lines are the stochastic scaled recallable reward processes, $\qol(\cdot)$, and the dashed lines the fluid solution, $q(\cdot)$ (Eq.~\eqref{eq:drift}); the lines, from top to bottom, correspond to actions $1$ through $3$, respectively. Plot $(d)$ corresponds to the memory-deficient regime, where the update rate $\beta$ is much smaller than $\mu$.}
\label{fig:tranResult}
\end{figure}

\section{Related Literature}
\label{sec:litRev}


{The assumption of a reward's lifespan being exponentially distributed is modeled after the exponential memory decay theory of Ebbinghaus (\cite{ebbinghaus1913memory}), which posits that the recall probability of a past event decays exponentially as time passes. The reward matching decision rule (Eq.~\eqref{eq:RMsampling}) can be viewed as a generalization of Luce's linear probabilistic choice rule (\cite{luce1959individual}), originally created to model empirical observations where humans or animal subjects make choices with probabilities proportional to associated amount of past rewards or stimuli. Luce's rule has been studied extensively in a variety of disciplines, ranging across cognitive science (\cite{herrnstein1970law}), behavioral economics (\cite{erev1998predicting}), and evolutionary biology (\cite{harley1981learning}); more recently, in a different context, a similar proportional-sampling inference algorithm highly related to the Luce's rule has been used in establishing the sample complexity of Bayesian active learning with privacy constraints (\cite{xu2018query}). Noteably, the celebrated work of \cite{erev1998predicting} demonstrates that a discrete-time version of Luce's rule can be a powerful  reinforcement learning model for explaining a player's behavior in sequential games, and subsequent theoretical analysis showed that the empirical frequency of the best action under Luce's rule converges to one in single player games (\cite{rustichini1999optimal, beggs2005convergence}, and more recently,  \cite{mertikopoulos2016learning}). These results, however, assume the agent has perfect memory, while models involving memory decay have received relatively little attention. \cite{beggs2005convergence} analyzes a variation of the model involving ``forgetting'' as a means to speed up convergence to mixed equilibria, with a crucial feature that the forgetting be performed in a noiseless fashion by weighing distant past experience with a deterministic discount factor. As such, among other differences, the stochastic memory decay model that we adopt exhibits fundamentally different dynamics than the deterministic discount model; for instance, under our model there is a non-negligible probability that all past rewards could be entirely erased, while the deterministic version never truly forgets the past, no matter how distant. }

The impact of limited memory to sequential decision-making has been examined in the statistics literature dating back to the seminal work of \cite{cover1970two,hellman1970learning}, which proposed algorithms for solving multi-arm bandit and sequential hypothesis testing problems when the decision maker has access to only a finite number of bits of memory. In the queueing theory literature, \cite{mitzenmacher2002load} and \cite{gamarnik2016delay} study load balancing problems with space-bounded routing algorithms. In contrast to the present paper, the memories in these models are assumed to be perfect and immune to random erasures. 

The dynamic induced in our model is  related to the celebrated P{\'o}lya's urn scheme (cf.~\cite{pemantle2007survey}). 
In its basic form, one ball is drawn uniformly at random from an urn containing balls of different colors. The chosen ball is then returned to the urn along with a new ball of the same color, and the procedure repeats. Since the probability that a certain color being chosen is proportional to the number of balls of that color present in the urn, by viewing the balls as rewards and colors as actions, our model can be thought of as the urn scheme with the modifications where (a) each ball disappears from the urn after some random amount of time, corresponding to memory loss, and (b) the number of new balls added is random,  whose distribution depends on the color, corresponding to the variations in the reward rates. Note that in P{\'o}lya's urn scheme the total number of balls in the urn eventually tends to infinity as the number of draws increases, while in our model the recallable rewards always stay a finite random variable as a result of memory loss. As such, our model can be viewed as  a stationary variation of the urn scheme, which admits a very different dynamic.

On the technical end, our analysis of the memory-abundant regime relies on a certain fluid model to approximate the evolution of the original, stochastic recallable reward process. While the use of fluid models is a popular technique in the literature of queueing networks \cite[cf.][]{kurtz1970solutions, bramson1998state, tsitsiklis2012power, massoulie2018capacity}, our work departs from conventional fluid models in a notable way: the recallable reward process is not Markovian, and the analysis hence must take into account the dynamics of the auxiliary choice process. It is worth noting that there has been a large-deviation theory developed for obtaining fluid-like approximation results where the transition rates of a continuous-time jump process are modulated by a finite-state auxiliary chain  (cf.~Theorem 1 of \cite{mitzenmacher2002load}, which is based on Theorem 8.15 of \cite{shwartz1995large}). Unfortunately, these results  do not apply in our model for two reasons. First, they require the transition rate of the jump process of interest to be bounded, a condition not met here due to the aggregate departure rate of the rewards being unbounded over the state space. Second, the existing results typically have one scaling parameter, $m$, whereby time is sped up by a factor of $m$, and space scaled down by a factor of $1/m$. We use a similar scaling for the recallable reward process, $Q^m(\cdot)$, while the update rate $\beta_m$ serves as a second scaling parameter that is not captured by the model in \cite{mitzenmacher2002load,shwartz1995large}. The addition of $\beta_m$ is non-trivial: in fact, if $\beta_m$ is too small, the fluid approximation does not hold, and a different (memory-deficient) regime arises. The dynamics of our model in the memory-abundant regime are also related to the so-called averaging principle in establishing fluid approximation for Markov processes, where the evolution of a Markov process at a slow time scale is influenced by average behavior of a modulating process at a faster time scale (cf.~\cite{perry2011ode,mandelbaum1998strong,evdokimova2018coupled,coutin2010markov}; see also \cite{darling2008differential} for an overview.) However, in these models the underlying slow process is Markovian (e.g., queue lengths)  and the modulating process a function of its value (e.g., the differences in queue lengths in \cite{perry2011ode}). In contrast, the recallable reward process in our model is not Markovian due to the existence of the choice process. The non-Markovian nature of our model means that certain martingale properties typically used to establish convergence do not directly follow from existing results. As a result of these aforementioned differences, we will develop our fluid approximation results starting from first principles using arguments developed by \cite{kurtz1978strong}, involving elementary properties of continuous-time martingales and certain {quasi-fluid processes}.

\section{Proof Overview}
\label{sec:proof_over_view}
The remainder of the paper is devoted to the proof of Theorem \ref{thm:main}; the proof of Theorem \ref{thm:Ck} is a straightforward extension of Theorem \ref{thm:main} and will be given in Section \ref{sec:thmCk}. We discuss in this section some high-level ideas before delving into the details. Our proof consists of two main components corresponding to the memory-abundant and memory-deficient regimes, respectively. The system dynamics turn out to be quite different depending on the scaling of $\beta_m$, and as a result, so are our proof techniques. 
Fix $m \in \N$, and define the joint reward-choice process 
\begin{equation}
W^m(t) = (Q^m(t), C^m(t)), \quad t\in \rp. 
\end{equation}
It is not difficult to check that the joint process $\{W^m(t)\}_{t\in \rp}$ is Markovian, while, for any fixed $m$, the component processes $\{Q^m(t))\}_{t\in \rp}$ and $\{C^m(t)\}_{t\in \rp}$ depend on one another and are {not} individually Markovian (see Figure \ref{fig:QandC}). However, a key observation is that one of the two component processes becomes {asymptotically Markovian} as $m\to \infty$: 

\begin{enumerate}
\item (Memory-abundant) If $\beta_m \gg 1/m$ as $m\to \infty$, a large number of choice updates takes place before substantial change occurs in the recallable reward process. As a result, the recallable reward process,  $\{Q^m(t))\}_{t\in \rp}$ , becomes a sufficient Markov representation of the system dynamics, as $m\to \infty$.

\item (Memory-deficient) If $\beta_m \ll 1/m$  as $m\to \infty$, the interval between two consecutive update points is large. By the time a new choice is to be made, the agent will have ``forgotten'' nearly all of the rewards associated with the actions other than the most current choice, and the rewards associated with the current choice will have become quite predictable.  Here, the choice process  $\{C^m(t))\}_{t\in \rp}$ becomes asymptotically Markovian as $m\to \infty$.
\end{enumerate}

Our proofs for the two regimes will be tailored to the dichotomy outlined above. \begin{enumerate}
\item For the memory-abundant regime, our analysis focuses on characterizing the scaled recallable reward process, $\qol(\cdot)$, and uses a certain fluid model to show that, as $m\to \infty$, the evolution of $\qol(\cdot)$  is well approximated by the solution to a system of ordinary differential equations (ODE), and its stationary distribution converges to the unique invariant state of the ODEs. This portion of the proof will be presented in Section \ref{sec:memAbundant}. 

\item For the memory-deficient regime, we turn instead to the choice process, $C^m(\cdot)$, and use its asymptotic Markovian property to obtain an explicit formula for its steady-state distribution, in the limit as $m\to \infty$. We then leverage this formula to calculate the expected steady-state recallable reward using the stationarity properties of the continuous-time processes. The proof for this portion of the theorem is presented in Section \ref{sec:memDeficient}. 
\end{enumerate}

\section{The Memory-abundant Regime}
\label{sec:memAbundant}
We establish the first claim of Theorem \ref{thm:main} in this section, and we will do so by actually proving a stronger statement, that the steady-state distribution of $\qol_k(\cdot)$ converges in $L^1$ to the corresponding expressions in Eqs.~\eqref{eq:MainFast1} and \eqref{eq:MainFast2}, as $m \to \infty$. As was alluded to earlier, our main tool relies on certain fluid solutions, defined as the solutions to a set of ordinary differential equations, whose dynamics will be shown to closely approximate that of $\qol(\cdot)$ when $m$ is large. We begin by introducing the fluid solutions and some of their basic properties. 

\subsection{Fluid Solutions and Their Basic Properties}

Define the function $p: \rp^\calK \to \rp^\calK$, where
\begin{equation}
p_k(q) = \frac{q_k \vee \alpha_0}{\sum_{i\in \calK}q_i\vee \alpha_0}, \quad q\in \rp^\calK, k \in \calK.
\label{eq:pdef} 
\end{equation}

\begin{definition}[Fluid Solutions]
\label{def:fluidSol}
Fix $q^0 \in \rp^{\calK}$. A continuous function, $q: \rp \to \rp^\calK$, is called a \emph{fluid solution} with initial condition $q^0$, if it satisfies the following:  
\begin{enumerate}
\item $q(0)  = q^0$; 
\item for almost all $t\in \rp$, the function is differentiable along all coordinates, with 
\begin{equation}
\dot{q}_k(t) = \lambda_k p_k(q(t)) - q_k(t), \quad \forall k \in \calK, 
\label{eq:drift}
\end{equation}
where $p(\cdot)$ is defined in Eq.~\eqref{eq:pdef}. 
\end{enumerate}
We say that $q^I \in \rp^K$ is an \emph{invariant state} of the fluid solutions, if by setting $q^0 = q^I$, we have that $q(t) = q^I$ for all $t\in \rp$. 
\end{definition}

The first term on the right-hand side of Eq.~\eqref{eq:drift} corresponds to  the instantaneous arrival rate of the rewards, where $p_k(q(t))$ is the probability that action $k$ is chosen given the current states of recallable rewards, and $\lambda_k$ the reward rate for that action. The second term corresponds to the instantaneous departure rate of rewards, which is proportional to the amount of recallable rewards associated with action $k$.  The following lemma states some basic properties of the fluid solutions. The proof is given in Appendix \ref{app:lem:fluidUniqu}. 

\begin{lemma}
\label{lem:fluidUniqu}
For all $q^0 \in \rp^K$,  there exists a unique fluid solution $\{q(q^0, t)\}_{t\in \rp}$ with initial condition $q^0$. Furthermore, for all $t\in \rp$, $q(q^0, t)$ is a continuous function with respect to $q^0$. 
\end{lemma}

The next result states that the fluid solutions admit a unique invariant state. Note that the expressions of the invariant state coincide with those in Eqs.~\eqref{eq:MainFast1} and \eqref{eq:MainFast2} in Theorem \ref{thm:main}. 

\begin{theorem}
\label{thm:invStat}
Fix $\alpha_0 >0$. The fluid solutions admit a unique invariant state, $q^I$, whose expressions depend on the value of $\alpha_0$, as follows. 

\begin{enumerate}
\item Suppose $\alpha_0\in (0, \lambda_1/K]$. Then,
 \begin{align}
q^I_k = \left\{ \begin{array}{ll}
          \lambda_1 - (K-1)\alpha_0, & \quad k=1,\\
          (\lambda_k/{\lambda_1})\alpha_0 , & \quad  k = 2, \ldots, K. \\
         \end{array} \right. 
         \label{eq:invrStates}
\end{align}
\item Suppose $\alpha_0 >  \lambda_1/K$. Then,
\begin{equation}
q^I_k = \frac{\lambda_k}{K}, \quad k = 1, \ldots, K. 
\label{eq:invStates2}
\end{equation}
\end{enumerate}
\end{theorem}

\emph{Proof.~} It can be easily verified  that the expressions in Eqs.~\eqref{eq:invrStates} and \eqref{eq:invStates2} are valid invariant states in their respective regimes of $\alpha_0$ values. We next show that they are indeed the unique invariant states in both cases. Let $q^I$ be an invariant state of the fluid solutions. Then, $q^I$ must satisfy: 
\begin{equation}
\lambda_k \frac{q^I_k \vee \alpha_0 }{\sum_{i \in \calK}  q^I_i \vee\alpha_0  } - q^I_k = 0, \quad \forall k \in \calK. 
\label{eq:condInvStat}
\end{equation}
Define a partition of $\calK$ into $ \calK^+ $ and $\calK^-$, where: 
\begin{align}
\calK^+ = & \lt\{ k \in \calK: q^I_k > \alpha_0 \rt\}, \quad \calK^- =  \lt\{ k \in \calK: q^I_k \leq \alpha_0 \rt\}. 
\end{align}
Note that by definition,
\begin{align}
\label{eq:qiinK+}
q^I_k \vee \alpha_0 = q^I_k, \quad k\in \calK^+, \\
q^I_k \vee \alpha_0 = \alpha_0, \quad k\in \calK^-. 
\label{eq:qiinK-}
\end{align}

Suppose that $\alpha_0 \in (0, \lambda_1/K]$. We next show that in this case $\calK^+$ is the singleton set, $\{1\}$. First, suppose that $\calK^+$ is empty, that is, $q^I_k \leq \alpha_0$ for all $k \in \calK$. We have that
\begin{equation}
q^I_1= \lambda_1 \frac{\alpha_0}{(\sum_{i \in \calK^+}q^I_i) + \alpha_0|\calK^-|} = \lambda_1\frac{\alpha_0}{K\alpha_0} = \frac{\lambda_1}{K} > \alpha_0,
\label{eq:qI1contra2}
\end{equation}
where the last inequality follows from the assumption that $\alpha_0<\lambda_1/K$. This leads to a contradiction, implying that $\calK^+ \neq \emptyset$. Fix $k\in \calK^+$. Combining Eqs.~\eqref{eq:condInvStat}, \eqref{eq:qiinK+} and \eqref{eq:qiinK-}, we have that
\begin{equation}
\lambda_k \frac{q^I_k}{(\sum_{i \in \calK^+}q^I_i) + \alpha_0|\calK^-|}= q^I_k, 
\end{equation}
which, after rearrangement, can be written as
\begin{equation}
\sum_{i \in \calK^+}q^I_i = \lambda_k - \alpha_0|\calK^-|. 
\label{eq:K+set}
\end{equation}
Note that in Eq.~\eqref{eq:K+set} only the term $\lambda_k$ on the right-hand depends on $k$. This implies that, if $\calK^+ \neq \emptyset$, then all coordinates in $\calK^+$ must have the same arrival rate of rewards. Thus, we can assume that there exists  $\tilde{\lambda} \in \rp$, such that $\lambda_k = \tilde{\lambda}$, for all $k\in \calK^+$,  and it suffices to show that $\tilde{\lambda}= \lambda_1$. 

For the sake of contradiction, suppose that $\tilde{\lambda}\neq \lambda_1$, and hence $\lambda_1>\tilde \lambda$. This implies that $1 \notin \calK^+$, and hence $q^I_1 \leq \alpha_0$. Invoking Eq.~\eqref{eq:condInvStat} again, we have that
\begin{equation}
q^I_1= \lambda_1 \frac{\alpha_0}{(\sum_{i \in \calK^+}q^I_i) + \alpha_0|\calK^-|} \sk{a}{=} \alpha_0 \frac{\lambda_1}{\tilde{\lambda}} \sk{b}{>} \alpha_0, 
\label{eq:qI1contra}
\end{equation}
where step $(a)$ follows from Eq.~\eqref{eq:K+set}, and $(b)$ from the fact that $\lambda_1>\tilde \lambda$. Clearly, Eq.~\eqref{eq:qI1contra} contradicts with the fact that $q^I_1 \leq \alpha_0$. We thus conclude that $\tilde{\lambda}=\lambda_1$, $\calK^+ = \{1\}$, and $\calK^- = \{2, \ldots, K\}$. By setting $k=1$ in  Eq.~\eqref{eq:K+set}, we further conclude that
\begin{equation}
q^I_1 = \sum_{i \in \calK^+}q^I_i = \lambda_1 - \alpha_0|\calK^-| = \lambda_1 - (K-1)\alpha_0. 
\label{eq:qI1proof}
\end{equation}
Now, fix $k\in \{2, \ldots, K\}$. We have, from Eq.~\eqref{eq:condInvStat}, that 
\begin{equation}
q^I_k = \lambda_k \frac{\alpha_0}{(\sum_{i \in \calK^+}q^I_i) + \alpha_0|\calK^-|}  = \alpha_0\frac{\lambda_k}{\lambda_1}, 
\label{eq:qIkproof}
\end{equation}
where the last equality follows from the fact that 
\begin{equation}
\lt(\sum_{i \in \calK^+}q^I_i \rt) + \alpha_0|\calK^-| = q^I_1 + (K-1)\alpha_0 = [\lambda_1 - (K-1)\alpha_0] +(K-1)\alpha_0 = \lambda_1. 
\end{equation}
Combining Eqs.~\eqref{eq:qI1proof} and \eqref{eq:qIkproof}, we have shown that the expressions in Eq.~\eqref{eq:condInvStat} are the unique invariant state of the fluid model, when $\alpha_0 \in (0, \lambda_1/K]$. 

We  turn next to the second case of Theorem \ref{thm:invStat}, and assume that $\alpha_0 > \lambda_1/K$. First, suppose that $\calK^+\neq \emptyset$. By the arguments leading to Eq.~\eqref{eq:qI1contra}, which do not depend on the value of $\alpha_0$, we know that $\calK^+ = \{1\}$, and 
\begin{equation}
q^I_1 = \lambda_1 - (K-1)\alpha_0 \leq \lambda_1 - (1-1/K)\lambda_1 = \frac{\lambda_1}{K} \leq \alpha_0, 
\end{equation}
where both inequalities follow from the assumption that $\alpha_0 > \lambda_1/K$. This contradicts with the assumption that $1 \in \calK^+$ and hence $q^I_1 > \alpha_0$, and we conclude that $\calK^+$ must be empty. This implies that $q^I_k \leq \alpha_0$ for all $k\in \calK$, and by Eq.~\eqref{eq:condInvStat}, we have that
\begin{equation}
q^I_k = \lambda_k \frac{\alpha_0}{ \alpha_0 K }  = {\lambda_k}{K}, \quad k = 1, \ldots, K.  
\label{eq:qIkproof2}
\end{equation}
This completes the proof of Theorem \ref{thm:invStat}. \qed

The following theorem is the main result of this section, which would imply the first claim of Theorem \ref{thm:main}.  

\begin{theorem}
\label{thm:steady}
 Fix $m\in \N$. The process $\qol(\cdot)$ admits a unique steady-state distribution, denoted by $\pi^m$. Suppose that $\beta_m \gg 1/m$ as $m\to \infty$. We have that
\begin{equation}
\pi^m \mbox{ converges weakly to } \delta_{q^I}, \quad \mbox{as $m\to \infty$}, 
\label{eq:piQconverg}
\end{equation}
where $q^I$ is the unique invariant state of the fluid model, and $\delta_{q^I}$ is the probability measure with unit mass on $q^I$. Furthermore, we have that
\begin{equation}
\lim_{m \to \infty}  \E(|\qol_k(\infty)-q^I_k|) =  0, \quad \forall k \in \calK. 
\label{eq:convergInL1final}
\end{equation}
\end{theorem}

\begin{figure}[h]
\centering
\includegraphics[scale=.95]{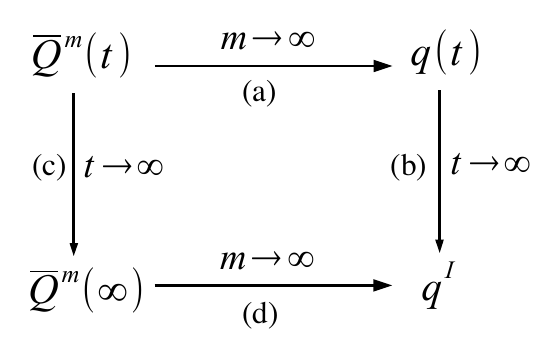}
\caption{An illustrative summary of the convergence results for the memory-abundant regime. The arrows $(a)$ and $(b)$ correspond to Theorems \ref{thm:MFE} and  \ref{thm:fluidToInvar}, respectively, and $(c)$ and $(d)$ to Theorem \ref{thm:steady}.  }
\label{fig:convergePlot}
\vspace{-10pt}
\end{figure}

The remainder of Section \ref{sec:memAbundant} is devoted to showing Theorem \ref{thm:steady} and consists of three main parts, as illustrated in Figure \ref{fig:convergePlot}.
\begin{enumerate}
\item Theorem \ref{thm:MFE} of Section \ref{sec:MFE} shows that the scaled recallable-reward process converges uniformly to a fluid solution over any finite time horizon, as $m\to \infty$. This is the more technical portion of the proof, and the key idea is to use a {quasi-fluid} process, one where the evolution of the reward process is deterministic while the choice process remains random, as an intermediary to bridge the gap between the fluid solution and the original stochastic process. 
\item Theorem \ref{thm:fluidToInvar} of Section \ref{sec:globalStab} shows that starting from any bounded initial condition, $q^0$, the fluid solution $q(q^0,t)$ converges exponentially fast to the unique invariant state, as $t\to \infty$. The proof is centered around the evolution of a simple piece-wise linear potential function, which ensures that $q(q^0,t)$ rapidly enters a suitable subset of the state space, from which point exponential convergence takes place. 
\item Finally, in Section \ref{sec:SSconverge}, by combining Theorems \ref{thm:MFE} and \ref{thm:fluidToInvar}, we show that the sequence of steady-state scaled recallable rewards converges in $L_1$ to a unique limiting point concentrated on the invariant state of the fluid solution, completing the proof of Theorem \ref{thm:steady}. 
\end{enumerate}

\subsection{Convergence to the Fluid Solution over Finite Horizon}
\label{sec:MFE}

In this subsection, we show that, with high probability, the scaled recallable reward process converges to the fluid solution uniformly over any finite horizon, as $m\to \infty$, as summarized in the following theorem (See plots $(a)$ through $(b)$ of Figure \ref{fig:tranResult} for a sample-path example of how this convergence takes place.) 

\begin{theorem}
\label{thm:MFE}
Suppose that $\beta_m \gg 1/m$, as $m\to \infty$. Fix $q^0 \in \rp^\calK$, and suppose that 
\begin{equation}
  \lim_{m \to \infty} \pb(\| {\qol(0)} - q^0\|>\epsilon ) = 0, \quad \forall \epsilon>0. 
\end{equation}
Let $q(q^0,\cdot)$ be the unique fluid solution with initial condition $q^0$, as defined in Lemma \ref{lem:fluidUniqu}. We have that, for all $T>0$, 
\begin{equation}
\lim_{m\to\infty} \pb\left(  \sup_{0\le t\le T}  \left\| \qol(t)-q (q^0,t) \right\| > \epsilon\right) =0, \quad \forall \epsilon>0. 
\end{equation}
\end{theorem}

A main challenge in proving Theorem \ref{thm:MFE} is that the processes  $Q^m(\cdot)$ and  $C^m(\cdot)$ interact in an intricate way: the recallable reward vector influences how the subsequent choices are chosen, while the current choice in turn dictates the arrival rate of new rewards associated with each action. The main idea of the proof is to ```disentangle'' their interactions by means of an intermediate, \emph{quasi-fluid} process, where the dynamics of the arrivals and departures of rewards becomes deterministic and ``fluid-like'', while the choice process remains random.

We now construct the quasi-fluid process. Fix $q\in \rp^\calK$ and $c\in \calK$, and define
\begin{equation}
 G_k(q,c) = \lambda_k \mathbb{I}(c = k) - q_k, \quad k \in \calK. 
 \label{eq:gkdef}
 \end{equation} 
In words, $G_k(q,c)$ represents the rate of change in recallable rewards associated with action $k$ when the current recallable reward vector and action are $q$ and $c$, respectively. Define the scaled action process: 
\begin{equation}
\col(t) = C^m(mt), \quad t\in \rp.  
\end{equation}
Note that, unlike $\qol(\cdot)$, the scaling of $\col(\cdot)$  only occurs in time and not in space. We then construct a \emph{quasi-fluid solution}, $\{V^m(t)\}_{t \in \rp}$, by letting 
\begin{equation}
V^m_k(t) = q^0_k + \int_0^t G_k(V^m(s),\col(s))ds, \quad k \in \calK, t\in \rp. 
\label{eq:dif-1}
\end{equation}
In words, $V^m(\cdot)$ corresponds to a system in which  recallable rewards are governed by deterministic arrival and departure processes, as specified by Eq.~\eqref{eq:gkdef}, and where the choice process is scaled in time by a factor of $m$, but remains stochastic. 

We note some basic properties of $V^m_k(\cdot)$ that will become useful. First, since $C^m(\cdot)$ is a Markov jump process and the function $G_k(\cdot,\cdot)$ uniformly Lipschitz continuous, analogously to Lemma \ref{lem:fluidUniqu}, it is not difficult to show that $V^m(\cdot)$, expressed as a solution to the above integral equation, is uniquely defined almost surely. Also, it is not difficult to verify that 
\begin{equation}
0\leq  V^m_k(t)  \leq q^0_k + \lambda_kt, \quad \forall t\in \rp, 
\label{eq:vmkbound}
\end{equation}
where the right-hand side of the second inequality corresponds to the scenario where $\col(s)=k$ for all $s \in [0,t]$ and the negative drift term $(-q_k)$ is absent from the expression of $G_k(q,c)$.

The proof of Theorem \ref{thm:MFE} proceeds in two main steps. We first show that $\{\qol(t)\}_{t\ge 0}$ converges to the quasi-fluid process, $\{V (t)\}_{t\ge 0}$, as $m \to \infty$, as is summarized in the following proposition. 

\begin{proposition} \label{prop:bd1}
 Fix $q^0 \in \rp^\calK$, and suppose that 
\begin{equation}
  \lim_{m \to \infty} \pb(\| {\qol(0)} - q^0\|>\epsilon ) = 0, \quad \forall \epsilon>0.
  \label{eq:prop1initial} 
\end{equation}
Then, 
\begin{equation}
\lim_{m\to\infty} \pb\lt( \sup_{0\le t\le T}  \left\|\qol (t) - V^m(t) \rt\| > \epsilon \right) =0, \quad \forall \epsilon>0. 
\end{equation}
\end{proposition}

We then show that, when $\beta_m \gg  1/m$, the process $\{V^m (t)\}_{t\ge 0}$ converges to the fluid solution, $\{q (q^0, t)\}_{t\ge 0}$: 
\begin{proposition} 
\label{prop:bd2}
Suppose, in addition to the condition of Eq.~\eqref{eq:prop1initial} in Proposition \ref{prop:bd1}, that $m \beta_m \to \infty$, as $m\to \infty$. For all $\epsilon>0$, 
\begin{equation*}
 \lim_{m \to \infty} \pb\left( \sup_{0\le t\le T} \left\| V^m(t)-q (q^0, t) \right\| > \epsilon \right) =0.
\end{equation*}
\end{proposition}
Note, in particular, that Proposition \ref{prop:bd1} holds {independently} of the scaling of $\beta_m$, while Proposition \ref{prop:bd2} requires that the system be in the memory-abundant regime.  Together, Propositions~\ref{prop:bd1} and \ref{prop:bd2} imply Theorem \ref{thm:MFE}. 

The proof for  Propositions~\ref{prop:bd1} and \ref{prop:bd2} are given in Appendix~\ref{sec:proof-lemma-refl} and \ref{sec:proof-lemma-refl-1}, respectively. The main rationale in executing the proof of Theorem \ref{thm:MFE} via two steps is as follows. In Proposition \ref{prop:bd1}, by {fixing} an arbitrary sample path of the choice process, $C(\cdot)$, the arrivals and departures from $Q^m(\cdot)$ become locally Poisson with a rate modulated by the current value of $C(t)$, and thus allowing us to establish a martingale property that will lead to the concentration of $Q^m(\cdot) $ around the quasi-fluid solution (Eq.~\eqref{eq:dif-1}). With Proposition \ref{prop:bd1}, the only remaining randomness in the quasi-fluid solution is in the choice process $C(\cdot)$ alone, and in Proposition \ref{prop:bd2} we show that its limiting behavior can be captured by a purely deterministic fluid solution in the memory-abundant regime. This two-step procedure is to be contrasted with analyzing the discrepancy between $Q^m(t)$ and the fluid solution directly, a more often adopted approach in the literature: in our model, the non-Markovian nature of $Q^m(t)$ implies that such differences do not immediately form a martingale. 

\emph{Remark.} We believe that there could be an alternative approach to establishing Theorem \ref{thm:MFE} by analyzing the limiting behavior of a discrete-time embedded Markov chain of the reward process $Q(\cdot)$ using stochastic approximation theory, and yet there are non-trivial difficulties along this direction. In particular, let $S_i$ be the $i$th update point (when the choice process $C(\cdot)$ can be updated), then the embedded reward process $\{Q[i] \bydef Q(S_i)\}_{i \in \N}$ can be shown to be a Markov chain, even though the original continuous-time process $Q(t)$ is not Markovian. One may then hope to invoke stochastic approximation theory to demonstrate that the  scaled version of $Q[\cdot]$ converges to a fluid solution, such as the approach taken in \cite{beggs2005convergence}. There are however several difficulties. First, the jumps in $Q[\cdot]$ can be unbounded due to the Poisson reward structure, while most standard results in the stochastic approximation literature require bounded updates between adjacent steps (e.g.,  Theorem 1 of \cite{benveniste2012adaptive}); it is likely possible to remove this requirement but it will require a more delicate argument. Second, even assuming such convergence has been established, it would only apply to the update points themselves, and to obtain the type of uniform convergence over an entire continuous time interval in the form of Theorem \ref{thm:MFE} would require additional procedures to bound the maximal fluctuation of the reward process between event points. In light of these observations, we tend to believe that the approach adopted in this paper is more direct and easier to carry out.

\subsection{Exponential Convergence of Fluid Solutions to Invariant State} 
\label{sec:globalStab}

In this section, we show that the fluid solution converges exponentially fast to the unique invariant state, starting from any bounded set of initial conditions. 

\begin{theorem} 
\label{thm:fluidToInvar}
Fix $h_0>0$. There exist $a,b >0$, such that for all $u \in \rp^K$, $\lnorm{u}\leq h_0$, 
\begin{equation}
 \lnorm{q(u,t)-q^I} \leq a\exp(-b t), \quad \forall t \in \rp. 
 \label{eq:qconvergeInv}
 \end{equation} 
\end{theorem}

\emph{Main idea.} The main difficulty in the proof of Theorem \ref{thm:fluidToInvar} is that, depending on the initial condition, the individual coordinates of $q(\cdot)$ may not converge to their equilibrium state monotonically, as can be seen in  the trajectory of $q_2(\cdot)$ in the first three plots of Figure \ref{fig:tranResult}. The main idea behind our proof is to obtain a more accurate description of $q(\cdot)$ with the aid of the following piece-wise linear {potential function}:
\begin{equation}
g(q) = \sum_{k \in \calK} q_k \vee \alpha_0,  \quad q\in \rp^\calK.
\label{eq:potFunct}
 \end{equation} 
By analyzing the evolution of $g(q(\cdot))$, our proof shows that: (1) for every bounded initial condition of $q(\cdot)$, the fluid solution rapidly enters a desirable subset of the state space, and (2) starting from this subset exponential convergence to the invariant state takes place. 

While potential functions are often used to establish convergence results, to the best of our knowledge, the use of the potential function $g(\cdot)$ in Eq.~\eqref{eq:potFunct} is new. This is largely due to the fact that our potential function is specifically designed to mirror the total ``weights'' across all actions under the reward matching rule. In addition, our use of the potential function also diverges from the their conventional role in proving convergence results: we do not use $g(\cdot)$ directly for characterizing the rate of convergence; the primary role of $g(\cdot)$ is for controlling the position of the trajectory in the state-space. Once we are assured that the trajectory has entered a suitable subset of the state space, the rate of convergence can be more readily obtained via elementary arguments.

\emph{Proof of Theorem \ref{thm:fluidToInvar}}.  Fix $u \in \rp^\calK$, such that $\lnorm{u}\leq h_0$. For simplicity of notation, we will omit the dependence of the initial condition and write $q(\cdot)$ in place of $q(u, \cdot)$,  whenever doing so does not cause any confusion.  The proof of the theorem is divided into two parts, depending on the value of $\alpha_0$. For the first part of the proof, we will assume that $\alpha_0\leq  \lambda_1/K$. 

Let $g(\cdot)$ be the potential function defined in Eq.~\eqref{eq:potFunct}. We will first partition $\rp^\calK$ into three disjoint subsets, depending on the value of $g(\cdot)$. Let $\epsilon$
\begin{equation}
\epsilon = \min\lt\{\frac{\lambda_1-\lambda_2}{2}, \lambda_1 - (K-1)\alpha_0 \rt\},
\end{equation}
and define
\begin{align}
\calQ_1 = & \{q \in \rp^\calK: g(q) \leq \lambda_1 - \epsilon \}, \nln
\calQ_2 = & \{q \in \rp^\calK:  \lambda_1- \epsilon <  g(q)  \leq \lambda_1 \}, \nln
\calQ_3 = & \{q \in \rp^\calK:  g(q) > \lambda_1 \}. 
\end{align}

We observe some useful properties of the fluid solutions.
\begin{lemma} 
\label{lem:partitionProp}
Fix $t>0$ so that all coordinates of $q(\cdot)$ are differentiable. The following holds: 
\begin{enumerate}
\item If $q(t) \in \calQ_1 \cup \calQ_2$, then
\begin{equation}
 \dot{q}_1(t) \geq 0.
 \end{equation} Furthermore, there exists a constant $d_1>0$, such that if $q(t) \in \calQ_1$, then
 \begin{equation}
  \dot{q}_1(t) \geq d_1.
  \end{equation} 
\item Denote by $\calK+(t)$ the coordinates over which $q(t)$ is at least $\alpha_0$: 
\begin{equation}
\calK^+(t) = \{k \in \calK: q_k(t) \geq \alpha_0\}. 
\end{equation}
There exists a constant $d_2>0$, such that if $q(t) \in \calQ_2 \cup \calQ_3$, then for all $k \in \calK^+(t)$, $k\geq 2$, 
\begin{equation}
\dot{q}_k(t) \leq -d_2. 
\label{eq:qkdecQ3}
\end{equation}
\item If $q(t) \in \calQ_3$, then
\begin{equation}
 \frac{d}{dt}g(q(t)) < 0.
 \end{equation} 
\end{enumerate}
\end{lemma}
\bpf
Suppose that $q(t) \in \calQ_1\cup \calQ_2$. We have that
 \begin{align}
\dot q _1(t)  =& \lambda_1\frac{q_1(t) \vee \alpha_0}{g(q(t))} - q_1(t) \nln
\geq &  \lambda_1\frac{q_1(t) \vee \alpha_0}{g(q(t))} - q_1(t)\vee \alpha_0 \nln
 =& (q_1(t)\vee \alpha_0 ) \lt( \frac{\lambda_1}{g(q(t))} -1  \rt) \nln
\sk{a}{\geq} & \alpha_0\lt(\frac{\lambda_1}{\lambda_1}-1\rt)  = 0, 
\end{align}
where step $(a)$ follows from the assumption that $q(t) \in \calQ_1\cup \calQ_2$. In a similar fashion, now suppose that $q(t) \in \calQ_1$ and hence $g(q(t))\leq \lambda_1-\epsilon$. We obtain that
\begin{equation}
\dot q _1(t){\geq} \alpha_0\lt(\frac{\lambda_1}{\lambda_1-\epsilon}-1\rt)  \bydef d_1>0. 
\end{equation}
This proves the first claim of the lemma. 

For the second claim, fix $q(t)\in \calQ_2 \cup \calQ_3$ and  $k \in \calK^+(t)$, $k\geq 2$.  We have that
\begin{align}
\dot q_k(t) =  &\lambda_k\frac{q_k(t) \vee \alpha_0}{g(q(t))}- q_k(t) \nln
=&  q_k(t) \lt( \frac{\lambda_k}{g(q(t))}-1 \rt) \nln
\sk{a}{\leq}   & \alpha_0 \lt( \frac{\lambda_k}{\lambda_1-\epsilon}-1 \rt) \nln
\sk{b}{\leq}   & \alpha_0 \lt( \frac{\lambda_2}{\lambda_2 + (\lambda_1-\lambda_2)/2}-1 \rt) <0, 
\end{align}
where step $(a)$ follows from the assumption that $q(t) \in \calQ_2 \cup \calQ_3$, and $(b)$ from the fact that $\lambda_2 \geq \lambda_k$ for all $k\geq 2$, and $\epsilon \leq \frac{\lambda_1-\lambda_2}{2}$. This proves the second claim. 

Finally, for the last claim, note that since $\alpha_0 \leq\lambda_1/K$, the fact that $q(t) \in \calQ_3$ implies that at least one component of $q(t)$ is no less than $\alpha_0$, and hence $\calK^+(t) \neq \emptyset$. We have that 
\begin{align}
\frac{d}{dt}g(q(t)) =&  \frac{d}{dt}\lt( \sum_{k\in \calK} q_k(t) \vee \alpha_0\rt)\nln
 = &   \sum_{k\in \calK^+(t)} \dot q_k(t) \nln
= &   \sum_{k\in \calK^+(t)} q_k(t) \lt(  \frac{\lambda_k}{g(q(t))} -1 \rt)  \nln
\sk{a}{<} &   \sum_{k\in \calK^+(t)} \alpha_0 \lt( \frac{\lambda_k}{\lambda_1} - 1 \rt) \nln
{\leq} & 0. 
\end{align}
where the strict inequality in step $(a)$ follows from the definition of $\calQ_3$, the fact that $\calK^+(t)$ is non-empty, and that $\lambda_k\leq \lambda_1$ for all $k$. This completes the proof of Lemma \ref{lem:partitionProp}. 
\qed

We proceed by considering two cases for the initial condition of $q(\cdot)$. 

{\bf Case 1: $q(0) \in  \calQ_1 \cup \calQ_2$.} Claim  3 in Lemma \ref{lem:partitionProp} implies that 
\begin{equation}
 q(t) \notin \calQ_3, \quad \forall t \in \rp. 
 \label{eq:qnotinQ3}
 \end{equation} 
By Claim 1 of the same lemma, we hence conclude that $q_1(t)$ is non-decreasing in $t$ for all $t \in \rp$. By the same claim, we observe that $\dot{q}_1(t)$  is bounded from below by the constant $d_1$ whenever $q(t)$ lies in $\calQ_1$.  Since $g(q(t)) \geq q_1(t)$,  we further conclude that in order for Eq.~\eqref{eq:qnotinQ3} to be valid, the amount of time $q(t)$ spends in $\calQ_1$ must satisfy: 
\begin{equation}
\int_{t\geq 0} \mathbb{I}(q(t) \in \calQ_1) d t \leq  \frac{\lambda_1}{d_1}, 
\label{eq:qtinQ1}
\end{equation}
which, combined with Eq.~\eqref{eq:qnotinQ3}, implies that for all $t\geq 0$: 
\begin{equation}
\int_{0}^t \mathbb{I}(q(s) \in \calQ_2) d s \geq t-  \frac{\lambda_1}{d_1}. 
\label{eq:timeinQ2}
\end{equation}
Claim 2 of Lemma \ref{lem:partitionProp} states that if $q(t)\in \calQ_2$, then any component of $q(t)$ in $\calK^+(t)\backslash \{1\}$ will have a negative drift that is bounded away from zero. Fix $t \in \rp$, and suppose there exists $k' \in \calK^+(t) \backslash \{1\}$. Since $\dot{q}_k(t)\leq \lambda_1$ for all $k$ and $t$, by Eq.~\eqref{eq:timeinQ2}, we have that
\begin{equation}
q_{k'}(t) \leq q_{k'}(0)+ \lambda_1\lt(\frac{\lambda_1}{d_1} \rt) - \lt(t-\frac{\lambda_1}{d_1} \rt)d_2 \leq  h_0+ \frac{\lambda_1^2}{d_1} - \lt(t-\frac{\lambda_1}{d_1} \rt)d_2.
\end{equation}
Since the components of $q(t)$ must be non-negative, the above equation implies that for all $k\geq 2$, 
\begin{equation}
q_k(t) \leq \alpha_0, \quad \forall t \geq T_0, 
\label{eq:qkbelowalpha}
\end{equation}
where $T_0 = \frac{h_0+ \lambda_1^2/d_1}{d_2}+ \frac{\lambda_1}{d_1}. $ 

Eq.~\eqref{eq:qkbelowalpha} will let us strengthen our previous observation. For all $t\geq T_0$, we have that
\begin{equation}
g(q(t)) = \sum_{k \in \calK} q_k(t) \vee \alpha_0 = q_1(t)\vee \alpha_0 + (K-1)\alpha_0. 
\end{equation}
Since $q_1(t)$ is non-decreasing in $t$ by the first claim of Lemma \ref{lem:partitionProp}, we conclude that $g(q(t))$ is non-decreasing in $t$ after $T_0$. Combined with Eq.~\eqref{eq:qtinQ1}, this implies that 
\begin{equation}
q(t) \in \calQ_2, \quad \forall t \geq T_1, 
\label{eq:qtinQ2}
\end{equation}
where  $T_1 = T_0+\frac{\lambda_1}{d_1}$.  Fix $t >  T_1$. Observe that by the definition of $\epsilon$, we have that $q_1(t) \geq \alpha_0$ whenever $q(t) \in \calQ_2$.  We have that
\begin{equation}
g(q(t)) = q_1(t) + (K-1)\alpha_0, 
\end{equation}
and hence
\begin{align}
\frac{d}{dt}(q^I-q_1(t)) = - \dot{q}_1(t) = - & \lambda_1\frac{q_1(t) }{q_1(t) +(K-1)\alpha_0} + q_1(t) \nln
 =& - q_1(t)\lt( \frac{\lambda_1}{q_1(t) + (K-1)\alpha_0}-1\rt)\nln
 \sk{a}{\leq} & - \alpha_0 \lt( \frac{\lambda_1}{q_1(t) + (K-1)\alpha_0}-1\rt)\nln
= & - \alpha_0 \lt( \frac{1}{1-(\lambda_1-(K-1)\alpha_0  - q_1(t) ) /\lambda_1}-1\rt)\nln
= & - \alpha_0 \lt( \frac{1}{1-(q^I_1  - q_1(t) ) /\lambda_1}-1\rt)\nln
\sk{b}{\leq} & - \frac{\alpha_0}{\lambda_1} ( q^I_1 -q_1(t)),
\label{eq:diffLowerbound}
 \end{align} 
where step $(a)$ follows from the fact that $q_1(t)\geq \alpha_0$. For step $(b)$, note that since $q(t) \in \calQ_2$ and $g(q(t)) = q_1(t) + (K-1)\alpha_0$, we have that
\begin{equation}
q_1^I - q_1(t) =  q_1^I -  \lt[ g(q(t)) - (K-1)\alpha_0 \rt] \in [0, \epsilon) \subset [0,\lambda_1/2). 
\end{equation}
The inequality in step $(b)$ thus follows from the fact that $1/(1-x) \geq 1+x$, for all $x\in [0,1)$, where we let $x = (q^I_1  - q_1(t))/\lambda_1$. 

Recall that, for $a>0$, the solution to the ordinary differential equation $\dot x(t) = -a x(t)$ is given by 
\begin{equation}
 x(t) = x(0)\exp(-at). 
 \label{eq:diffEsol1}
 \end{equation} 
From Eq.~\eqref{eq:diffLowerbound}, we conclude that for all $s >0$, 
\begin{align}
q^I_1 - q_1(T_1+s) \leq (q^I_1-q_1(T_1)) \exp\lt(- \frac{\alpha_0}{\lambda_1}s\rt) \leq \epsilon\exp\lt(- \frac{\alpha_0}{\lambda_1}s\rt), 
\label{eq:q1ODEbound}
\end{align}
where the last inequality follows from the fact that $q(T_1)\in \calQ_2$ and  $q_k(T_1)\leq \alpha_0$ for all  $k\geq 2$, and hence $q_1(T_1) \geq \lambda_1-(K-1)\alpha_0 - \epsilon = q^I_1-\epsilon$.  Let $a_1 = q^I_1\vee \epsilon$, and $b_1 = \alpha_0/\lambda_1$. Since $q_1(t)$ is non-decreasing by Lemma \ref{lem:partitionProp}, from Eqs.~\eqref{eq:qtinQ2} and \eqref{eq:q1ODEbound}, we have that if $q(0)\in \calQ_1\cup\calQ_2$, then
\begin{equation}
|q_1(t) - q^I_1| \leq a_1\exp(-b_1t), \quad \forall t \in \rp, 
\label{eq:q1convCase1}
\end{equation}

We now show the convergence of $q_k(\cdot)$ for $k\geq 2$. Fix $k\in \{2, \ldots, K\}$.  Define
\begin{equation}
\delta_k(t) =  q_k(t)-q^I_k, \quad k \in \calK, 
\end{equation}
and recall from Theorem \ref{thm:invStat} that $q^I_k = \frac{\lambda_k}{\lambda_1}\alpha_0$, $k = 2, \ldots, K$. 
By Eq.~\eqref{eq:qkbelowalpha}, we have that $q_k(t)\leq \alpha_0$ for all $t\geq T_0$. Therefore, for all $t\geq T_1$, 
\begin{align}
\dot{q}_k(t) = &\lambda_k \frac{\alpha_0}{g(q(t))} - q_k(t) \nln
= & \lambda_k \frac{\alpha_0}{q_1(t) + (K-1)\alpha_0} - q_k(t) \nln
= & \lambda_k \frac{\alpha_0}{\lambda_1 + \delta_1(t)} - q_k(t). 
\label{eq:qk2dot}
\end{align}
Using the Taylor expansion of the function $f(x) = \frac{\alpha_0}{\lambda_1+x}$ around $x=0$ and noting that $\delta_1(t) \leq a_1\exp(-b_1 t)$ (Eq.~\eqref{eq:q1convCase1}), we have that there exists $T_2 \geq T_1$, such that for all $t\geq T_2$, 
\begin{equation}
\dot{\delta}_k(t)  \in - \delta_k(t) \pm \frac{2\lambda_k\alpha_0}{\lambda_1}|\delta_1(t)| =  - \delta_k(t) \pm \gamma_k\exp(-b_1t),
\label{eq:dotdelta2Sandwi} 
\end{equation}
where $\gamma_k \bydef \frac{2\lambda_k\alpha_0 a_1}{\lambda_1}$, and we use the notation $x \in y\pm z$ to mean $y-z \leq x \leq y+z$. Given a fixed value of $\delta_k(T_2)$, Eq.~\eqref{eq:dotdelta2Sandwi} implies that the value of $\delta_k(T_2+t)$ must lie between the solutions to the ODEs
$\dot{x}(t)  = - x(t) +\gamma_k\exp(-b_1t)$ and $\dot{x}(t)  = - x(t) - \gamma_k\exp(-b_1t)$, with initial condition $x(0) = \delta_k(T_2)$, respectively. It can be verified that the solution to the ODE $\dot x (t) = - x(t)+c_1 \exp(c_2 t)$, with initial condition $x(0)=x_0$, is given by $x(t) = \frac{c_1}{1-c_2}\exp(-c_2 x)+\lt(x_0-\frac{c_1}{1-c_2}\rt)\exp(-x)$. Setting $c_1=\gamma_k$ and $c_2 = b_1$, we thus conclude that there exist $\overline{a}_k,{b}_k>0$, such that 
\begin{equation}
 |\delta_k(T_2+t)| \leq  (\delta_k(T_2)+\overline{a}_k)\exp(-b_kt). 
\end{equation}
Note that $\dot \delta_k(t)$ is always bounded. We thus conclude hat there exists $a_k>0$ such that,
\begin{equation}
|q_k(t) - q^I_k| = |\delta_k(t)| \leq a_k\exp(-b_k t), \quad \forall t \in \rp. 
\label{eq:qkconvexpokgeq2}
\end{equation}
Since the above equation holds for all $k \geq 2$, together with Eq.~\eqref{eq:q1convCase1}, we have proven  Eq.~\eqref{eq:qconvergeInv} for the first case, by assuming that $q(0) \in \calQ_1\cup \calQ_2$. 

{\bf Case 2: $q(0)\in \calQ_3$}. We now consider the second case where the initial condition, $q(0)$, belongs to the set $\calQ_3$. Recall from Eq.~\eqref{eq:qkdecQ3} of Lemma \ref{lem:partitionProp} that when $q(t) \in \calQ_2 \cup \calQ_3$, all components of $q(t)$ in $\calK^+(t)$ exhibit a negative drift of magnitude at least $d_2$. Letting $T_3 = h_0 / d_2$, we thus conclude that one of the following scenarios must occur: 
\begin{enumerate}
\item  $q(s') \in \calQ_1 \cup \calQ_2$ for some $s' \in [0, T_3]$. 
\item  $q(t) \in \calQ_3$ for all $t\in [0,T_3]$,  and there exists $s' \in [0,T_3]$ such that 
\begin{equation}
q_k(s') \leq \alpha_0, \quad \forall k\geq 2. 
\end{equation}
\end{enumerate}
For the first scenario, it is equivalently to ``re-initializing'' the system at time $t=s'$ at a point in $\calQ_1$.  The proof for Case 1 presented earlier thus applies and Eq.~\eqref{eq:qconvergeInv} follows. In what follows, we will focus on the second scenario. Without loss of generality, it suffices to show the validity of Eq.~\eqref{eq:qconvergeInv} by fixing an initial condition: 
\begin{equation}
q(0) \in \calQ_3, \mbox{ and }  \max_{k \geq 2} q_k(0) \leq \alpha_0. 
\end{equation}
Not that, by the definition of $\calQ_3$, this implies that 
\begin{equation}
q_1(0) >  q^I_1>\alpha_0. 
\end{equation}

We next show that $q_1(t) > q^I_1$ for all $t \in \rp$.  Fix $t>0$ such that 
\begin{equation}
q(s) \in \calQ_3, \quad \forall s \in [0,t]. 
\end{equation}
Such $t$ exists because the set $\calQ_3$ is open and  $q(\cdot)$ is continuous. From Lemma \ref{lem:partitionProp}, we know that if $q(t)\in \calQ_3$, then any component in $\calK^+(t) \backslash \{1\}$ must have a negative drift. It thus follows that 
\begin{equation}
q_k(t) < \alpha_0, \quad \forall k\geq 2. 
\label{eq:qk+tillT4}
\end{equation}
We have that, whenever $q^I_1< q_1(t)<q^I_1+1$, 
\begin{align}
\frac{d}{dt}\lt(q_1(t) - q^I_1\rt) = \dot q_1(t) \sk{a}{=} & - q_1(t)\lt(1- \frac{\lambda_1}{q_1(t) + (K-1)\alpha_0}\rt)  \nln
= & - q_1(t)\lt(1- \frac{1}{1+(q_1(t) - q^I_1)/\lambda_1} \rt) \nln
\geq & - (q^I_1+1)\lt(1- \frac{1}{1+(q_1(t) - q^I_1)/\lambda_1} \rt) \nln
\sk{b}{\geq} & - \frac{q^I_1+1}{2\lambda_1}(q_1(t) - q^I_1)
\label{eq:diffq1mq1I}
\end{align}
where step $(a)$ follows from the fact that $q_k(t)\leq \alpha_0$ for all $k\geq 2$, and $(b)$ from $1/(1+x) \leq 1-x/2$ for all $x\in [0, 1]$. In light of Eq.~\eqref{eq:diffEsol1}, the above inequality implies if the value of $q_1(t) - q^I_1$ is ever to be below $1$, then from that point on it will be bounded from below by an exponential function that is strictly positive. This shows that $q_1(t) > q^I_1$ for all $t>0$. 

We now prove the exponential rate of convergence of $q_1(t)$ to $q^I_1$. Using again Eq.~\eqref{eq:diffq1mq1I}, we obtain
\begin{align}
\frac{d}{dt}\lt(q_1(t) - q^I_1\rt) = & - q_1(t)\lt(1- \frac{\lambda_1}{q_1(t) + (K-1)\alpha_0}\rt) \nln
\sk{a}{\leq} & - q^I_1\lt(1- \frac{1}{1+(q_1(t) - q^I_1)/\lambda_1} \rt) \nln
\sk{b}{\leq} & - \frac{q^I_1}{\lambda_1}(q_1(t) - q^I_1),
\end{align}
where step $(a)$ follows from the fact that $q_1(t) > q^I_1$, and $(b)$ from $1/(1+x) \geq 1-x$ for all $x\geq 0$. Using the reasoning identical to that in Eqs.~\eqref{eq:diffLowerbound} through \eqref{eq:q1convCase1}, we conclude that there exists $a_1, b_1>0$, such that $|q_1(t) - q^I_1| \leq a_1\exp(-b_1t)$, for all $ t \in \rp. $

Finally, the fact that $q_1(t)> q^I_1$  implies that $g(q(t))> \lambda_1$, and hence $q(t)\in \calQ_3$ for all $t>0$. By Eq.~\eqref{eq:qk+tillT4}, this further implies that for any $k\geq 2$,  $q_k(t) \leq \alpha_0$ for all $t \in \rp$.  The exponential rate of convergence of $q_k(t)$ for $k\in \{2, \ldots, K\}$ thus follows from the same argument as in Eqs.~\eqref{eq:qk2dot} through \ref{eq:qkconvexpokgeq2}. This completes the proof of Theorem \ref{thm:fluidToInvar}, assuming that $\alpha_0 < \lambda_1/K$. 

We now turn to the case where $\alpha_0 > \lambda_1/K$. We will use a similar proof strategy by partitioning the state space based on the value of $g(q(t))$, albeit with a different partitioning. In particular, define the sets
\begin{align}
\tilde \calQ_1  = & \{q \in \rp^\calK: g(q) =  K\alpha_0\}, \nln
\tilde \calQ_2  = & \{q \in \rp^\calK: g(q) >  K\alpha_0\}. 
\end{align}
Note that by definition $g(q)\geq K\alpha_0$ for all $q\in \rp^\calK$, so the above sets constitute a partition of $\rp^\calK$. We have the following lemma.

\begin{lemma} 
\label{lem:qdriftildeQ2}
Fix $t>0$ so that all coordinates of $q(\cdot)$ are differentialble. If $q(t)\in \tilde \calQ_2$, then there exists a constant $d_3>0$, such that $\frac{d}{dt}g(q(t)) < -d_3$. \end{lemma}

\bpf Suppose that $q(t)\in \tilde \calQ_2$.  Since $g(q(t))>K\alpha_0$, we have that $\calK^+(t)\neq \emptyset$.  We have that 
\begin{align}
\frac{d}{dt}g(q(t)) =&  \frac{d}{dt}\lt( \sum_{k\in \calK} q_k(t) \vee \alpha_0\rt) =    \sum_{k\in \calK^+(t)} \dot q_k(t) \nln
= &   \sum_{k\in \calK^+(t)} q_k(t) \lt(  \frac{\lambda_k}{g(q(t))} -1 \rt)  \nln
\sk{a}{\leq} &  -  \sum_{k\in \calK^+(t)} \alpha_0 \lt( 1- \frac{\lambda_k}{K\alpha_0}  \rt) \nln
\sk{b}{\leq}&  - \alpha_0\lt( 1- \frac{\lambda_1}{K\alpha_0}  \rt) \bydef -d_3, 
\end{align}
where step $(a)$ follows from the fact that $g(q(t))>K\alpha_0$,  and $(b)$ from $\calK^+(t)\neq \emptyset$. Finally, note that since $\alpha_0 > \lambda_1/K$, $d_3$ is strictly positive. This completes the proof of the lemma. 
\qed

Lemma \ref{lem:qdriftildeQ2} implies that if $q(0) \in \tilde \calQ_2$, then we must have that 
\begin{equation}
q(t) \in \tilde \calQ_1, \quad \forall t \geq h_0/d_3. 
\end{equation}
Since the rate of change in every coordinate of $q(t)$ is bounded, in light of the above equation, it suffices to establish Eq.~\eqref{eq:qconvergeInv} by considering only the scenario where $q(t) \in \tilde \calQ_1$ for all $t \in \rp$. Fixing $k\in \calK$, we have that
\begin{align}
& \frac{d}{dt} \lt( q_k(t) - q^I_k \rt)   \nln
= &  \lambda_k \frac{q_k(t) \vee \alpha_0}{g(q(t))}- q_k(t) \nln
 \sk{a}{=} & \lambda_k \frac{\alpha_0}{g(q(t))}- q_k(t) \nln
 \sk{b}{=}  & \lambda_k/K - q_k(t) \sk{c}{=}  - (q_k(t)-q^I_k),
\label{eq:qkconvregime2}
\end{align}
where step $(a)$ follows from the fact that $q_k(t)\leq \alpha_0$ whenever $q(t) \in \tilde \calQ_1$, $(b)$ from the fact that $g(q(t)) = K\alpha_0$, and $(c)$ from the definition $q^I_k = \lambda_k/K$. In light of Eq.~\eqref{eq:diffEsol1}, Eq.~\eqref{eq:qkconvregime2} implies that 
\begin{equation}
|q_k(t) - q^I_k| = |q_k(0) - q^I_k| \exp(-t) \leq (h_0 \vee q^I_k) \exp(-t), \quad t \in \rp. 
\end{equation}
Since the above equation holds for all $k$, we have thus established an exponential rate of convergence of $q(t)$ to $q^I$, as $t\to \infty$. This completes the proof of Theorem \ref{thm:fluidToInvar}. 
\qed

\subsection{Convergence of Steady-State Distributions}
\label{sec:SSconverge}

We complete the proof of Theorem \ref{thm:steady} in this subsection. We will begin by establishing a useful distributional upper bound on $Q^m_k(\infty)$. Fix $t\in \rp$, and denote by $U^m_\lambda(t)$ the number of jobs in system at time $t$ in an $M/M/\infty$ queue with arrival rate $\lambda$ and departure rate $1/m$, and by $U^m_\lambda$ its steady-state distribution. It is well known that $U^m_\lambda$ is a Poisson random variable with mean $m\lambda$, and it follows that, for all $\lambda >0$, 
\begin{equation}
m^{-1}U^m_{\lambda} {\rightarrow} \lambda, \quad \mbox{almost surely, as $m\to \infty$}. 
\label{eq:Umlamconverg}
\end{equation}
Define the set 
\begin{equation}
\calW = \rp^\calK \times \calK,
\label{eq:calWdef}
\end{equation}
and the process
\begin{equation}
\wol(t) = (\qol(t), C^m(t)), \quad  t\in \rp. 
\end{equation} We have the following lemma, whose proof is given in Appendix \ref{app:lem:stochDom1}. 
\begin{lemma} 
\label{lem:stochDom1}
The process $\{\wol(t)\}_{t\in \rp}$  is positive recurrent, and  $\qol_k(\infty) \preceq   m^{-1}U^m_{\lambda_1}$, for all $k\in \calK$. \end{lemma}

Denote by $\pi_W^m$ the probability distribution over $\calW$ (Eq.~\eqref{eq:calWdef}) that corresponds to the steady-state distribution of $\wol(\cdot)$. Lemma \ref{lem:stochDom1} can be used to to show that the sequence $\{\pi^m_W \}_{m \in \N}$ is tight, as formalized in the following result, whose proof is given in Appendix \ref{app:lem:disTight}. 
\begin{lemma}
\label{lem:disTight}
For every $x \in \rp$,  define $ \cqol_x = \{ q \in \rp^\calK: \max_{k\in \calK} q_k \leq x \} \times \calK$, $x \in \rp. $ Then, for every $\epsilon > 0$, there exists $M\in \N$ such that 
\begin{equation}
\pi_W^m(\cqol_{M}) \ge 1- \epsilon, \quad \forall m\in \N, 
\label{eq:tightQ1}
\end{equation}
\end{lemma}

We say that $y$ is a limit point of  $\{x_i\}_{i \in \N}$ if there exists a sub-sequence of $\{x_i\}$ that converges to $y$.  It is not difficult to verify that the space $\calW$ is separable as a result of both $\rp^\calK$ and $\calK$ being separable. By Prohorov's theorem, the tightness of $\{\pi_W^m\}_{m \in \N}$ thus implies that any sub-sequence of $\{\pi_W^m\}_{m \in \N}$ admits a  limit point with respect to the topology of weak convergence. Let $\{\pi_W^{m_i}\}_{i \in \N}$ be a sub-sequence of $\{\pi_W^m\}_{m \in \N}$, and $ \pi_W$ a limit point of the sub-sequence. Denote by $ \pi$ and $\pi^{m_i}$ the marginals of ${\pi}_W$ and $\pi^{m_i}_W$ over the first $K$ coordinates  (corresponding to $\qol(\cdot)$), respectively. In the remainder of the proof, we will show that $\pi$ necessarily concentrates on $q^I$.

We first show that ${\pi}$ is a stationary measure with respect to the deterministic fluid solution. To this end, we will use the method of  continuous test functions (cf.~Section 4 of \cite{ethier2005markov}).  Let $\overline{\calC }$ be the space of all bounded continuous functions from $\rp^\calK$ to $\rp$. We will demonstrate that, for all $f \in \overline{\calC}$, 
\begin{equation}
\left|  \E_{\pi} \lt( f(q(q^0,t)) \rt)  - \E_{\pi} \lt(  f(q^0)\rt)  \right| =0, \quad \forall t \in \rp, 
\label{eq:hatpiinvariant}  
\end{equation}
where, from here onward, we will use the subscript, $\E_\pi(\cdot)$, to indicate the distribution of $q^0$. 

Define $\calQ^m = \{{x}/{m}, \,  x\in \zp^\calK\}$, $m \in \N.$ Fix $f\in \overline{\calC}$. Let $\wol(0)$ be distributed according to $\pi^m_W$, and define 
\begin{align}
F^m(q^0,t)  = & \E\left( \left. f(\qol (t)) \right| \qol(0)= q^0 \right),  \quad t\in \rp, q^0 \in \calQ^m. 
\end{align}
Fix $t >0 $.  We have that
\begin{align}
&\left|  \E_{\pi} \lt( f(q(q^0, t)) \rt)  - \E_{\pi} \lt(  f(q^0)\rt)  \right| \cr
 \stackrel{(a)}{\le} &  \limsup_{i \to \infty} \left|  \E_{\pi} \lt( f(q(q^0,t)) \rt) -  \E_{\pi^{m_i}} \lt( f(q(q^0,t)) \rt)  \right| \nln
 &  + \limsup_{i \to \infty} \left| \E_{\pi^{m_i}} \lt( f(q(q^0,t)) \rt)  - \E_{\pi^{m_i}}\lt( F^{m_i}(q^0,t)\rt)  \right|\nln
 & + \limsup_{i \to \infty} \left|  \E_{\pi^{m_i}}\lt( F^{m_i}(q^0,t)\rt)  - \E_{ \pi}\lt( f(q^0) \rt)  \right| \nln
 \stackrel{(b)}{\le} &  \limsup_{i \to \infty} \left|  \E_{\pi} \lt( f(q(q^0,t)) \rt) -  \E_{\pi^{m_i}} \lt( f(q(q^0,t)) \rt)  \right|  \nln
  &  + \limsup_{i \to \infty} \left| \E_{\pi^{m_i}} \lt( f(q(q^0,t)) \rt)  - \E_{\pi^{m_i}}\lt( F^{m_i}(q^0,t)\rt)  \right|\nln
 & + \limsup_{i \to \infty} \left|  \E_{\pi^{m_i}}\lt( f(q^0)\rt)  - \E_{ \pi}\lt( f(q^0) \rt)  \right| \nln
\stackrel{(c)}{=} &  \limsup_{i \to \infty} \left| \E_{\pi^{m_i}} \lt( f(q(q^0,t)) \rt)  - \E_{\pi^{m_i}}\lt( F^{m_i}(q^0,t)\rt)  \right|.
\label{eq:Tf}
\end{align}
Step $(a)$ follows from the triangle inequality, and $(b)$ from the fact that $Q^m(\cdot)$ is stationary when initialized according to $\pi^m$, and therefore, 
\begin{align}
\E_{\pi^m}\lt( F^m (q^0,t)\rt) =&\int_{q^0 \in \rp^\calK} \E\left( \left. f(\qol (t)) \right| \qol (0) = q^0 \right)  d\pi^m = \E_{\pi^m}\lt( f(q^0)\rt). \label{eq:steady-m}
\end{align}
 For step $(c)$, we note that since $q(q^0,t)$ depends continuously on $q^0$ (Lemma \ref{lem:fluidUniqu}), the function $g(x) = f(q(x,t)), \, x\in \rp^\calK$, belongs to $\overline{\calC}$. Therefore, the step follows from the fact that $\pi^{m_i}$ converges weakly to $ \pi$  as $i\to \infty$.

Fix $\epsilon>0$. There exists $M>0$ such that the the right-hand side of Eq.~\eqref{eq:Tf} satisfies
\begin{align}
&  \limsup_{i \to \infty} \left| \E_{\pi^{m_i}} \lt( f(q(q^0,t)) \rt)  - \E_{\pi^{m_i}}\lt( F^{m_i}(q^0,t)\rt)  \right| \cr
\le&\limsup_{i \to \infty} \left|  \int_{q^0 \in \cqol_{M}} f(q(q^0,t)) d\pi^{m_i} - \int_{q^0 \in \cqol_{M}}  F^m (q^0,t) d\pi^{m_i}  \right|+\cr
&\limsup_{i \to \infty} \left|  \int_{q^0 \in \rp^\calK \setminus \cqol_{M}} f(q(q^0,t)) d\pi^{m_i} - \int_{q^0 \in \rp^\calK \setminus \cqol_{M}}   F^m (q^0,t) d\pi^{m_i}  \right| \cr
\stackrel{(a)}{\le} & \limsup_{i \to \infty}  \int_{q^0 \in \cqol_{M}} \left|  f(q(q^0,t)) -   F^{m_i} (q^0,t)  \right|  d\pi^{m_i}  + 2 \epsilon  \sup_{x\in \rp^\calK} f(x) \nln
\sk{b}{=} & 2 \epsilon  \sup_{x\in \rp^\calK} |f(x)|. \label{eq:Tf-3}
\end{align}
where $(a)$ follows from the tightness property of Eq.~\eqref{eq:tightQ1}, and step $(b)$ involves an argument that allows for interchanging the order of taking the limit and integration. We isolate step $(b)$ in the following lemma; the proof leverages Theorem \ref{thm:MFE} and is given in Appendix \ref{app:lem:interchange1}.
\begin{lemma} 
\label{lem:interchange1}
Fix a compact set $S \subset \rp^\calK$ and $f\in \overline{\calC}$, we have that 
\begin{equation}
\limsup_{m \to \infty}  \int_{q^0 \in S} \lt| f(q(q^0,t)) -   F^{m} (q^0,t)  \right| d\pi^{m}   =0. 
\end{equation}
\end{lemma}

Fix $\delta>0$. We have that
\begin{align}
& \lim_{m\to \infty}   \int_{q^0 \in S} \left| f(q(q^0,t)) -   F^m (q^0,t)  \rt| d\pi^{m}    \nln
\sk{a}{\leq} &  \lim_{i \to \infty}  \int_{q^0 \in S}    \E \lt( \lt|  f(q(q^0,t)) - f(\overline{Q}^{m_i}(t)) \rt| \Bbar \overline{Q}^{m_i}(0)=q^0 \rt) d\pi^{m}  \nln
{\leq} & w_f(S, \delta) + \lt(\sup_{x \in \rp^\calK}|f(x)| \rt)  \lim_{i \to \infty} \sup_{x\in S\cap \calQ^{m}} \pb\lt( \|q(x,t)- \qol(x,t)\|>\delta\rt) \nln
\sk{b}{=} & w_f(S, \delta), 
\label{eq:uniformApprox3}
\end{align}
where $w_f$ is the modulus of continuity of $f$ in $S$: $w_{f}(S,\delta) = \sup_{x,y \in S, \|x-y\|\leq \delta} |f(x)-f(y)|$. Step $(a)$ follows from the fact that $\pi^m(\calQ^m)=1$, and $(b)$ from Eq.~\eqref{eq:uniformApprox}. Because $S$ is compact, $f$ is uniformly continuous in $S$ and hence $\lim_{\delta \to 0} w_{f}(S,\delta)=0$. The lemma follows by taking the limit as $\delta\to 0$ in Eq.~\eqref{eq:uniformApprox3}. \qed

Since Eq.~\eqref{eq:Tf-3} holds for all $\epsilon>0$, combining it with Eq.~\eqref{eq:Tf}, we conclude that 
\begin{equation}
\left|  \E_{\pi} \lt( f(q(q^0,t)) \rt)  - \E_{\pi} \lt(  f(q^0)\rt)  \right| \leq  \limsup_{i \to \infty} \left| \E_{\pi^{m_i}} \lt( f(q(q^0,t)) \rt)  - \E_{\pi^{m_i}}\lt( F^{m_i}(q^0,t)\rt)  \right| =0, 
\end{equation}
and this proves Eq.~\eqref{eq:hatpiinvariant}.  We now show that Eq.~\eqref{eq:hatpiinvariant} implies that $\pi = \delta_{q^I}$, i.e., $\delta_{q^I}$ is the unique invariant measure with respect to the dynamics induced by the fluid solution. Define the truncated norm: 
\begin{equation}
\|x\|_M \bydef  \min\{\|x\|,  M\}, \quad M \in  \rp, 
\end{equation}
and the set
\begin{equation}
\calR_L \bydef \{q \in \rp^\calK: \|q\|\leq L\} , \quad L \in  \rp, 
\end{equation}

Let $q^0$ be a random vector in $\rp^\calK$ distributed according to $\pi$. Suppose, for the sake of contradiction, that $\E(\|q^0-q^I\|_M)\bydef \mu_0>0$ for some $M>0$, where the expectation is taken with respect to the randomness in $q^0$. Fix $L>0$ such that $\pb(q^0 \notin \calR_L)\leq \frac{\mu_0}{4}\cdot \frac{1}{M}$. By Theorem \ref{thm:fluidToInvar}, there exists $T>0$ such that $\sup_{x \in \calR_L}\|q(x,T)-q^I\|_M < \mu_0/4$. Therefore, 
\begin{equation}
\E(\|q(q^0,T)-q^I\|_M) < \mu_0/4 + M\cdot \pb(q^0 \notin \calR_L) \leq \mu_0/2 = \E(\|q^0-q^I\|_M)/2, 
\end{equation}
where the first inequality uses the fact that $\|q\|_M\leq M$ for all $q\in \rp^\calK$. Because $\E(\|q^0-q^I\|_M)$ was assumed to be strictly positive, this  is in contradiction with the stationarity of $\pi$ with respect to the fluid solution, which would imply that $\E(\|q(q^0,T)-q^I\|_M)=\E(\|q^0-q^I\|_M)$ for all $T\in \rp$. We thus conclude that 
\begin{equation}
\E(\|q^0-q^I\|_M) = 0, \quad \forall M\in \rp, 
\end{equation}
which in turn implies that $\pi = \delta_{q^I}$. This proves Eq.~\eqref{eq:piQconverg}. 

Finally, to show Eq.~\eqref{eq:convergInL1final}, note that by Lemma \ref{lem:stochDom1}, for all $m\in \N$ and $k\in \calK$, the random variable $\qol_k(\infty)$ is non-negative and stochastically dominated by $m^{-1}U^m_{\lambda_1}$. Using similar steps as those in Eq.~\eqref{eq:tightShow}, it is not difficult to show that $\sup_{m \in \N}\E(|m^{-1}U^m_{\lambda_1}|^2)<\infty$, which, implies the uniform integrability of $\{\qol_k\}_{m\in \N}$. In combination with Eq.~\eqref{eq:piQconverg}, this proves the convergence in $L^1$ stated in Eq.~\eqref{eq:convergInL1final}. This completes the proof of Theorem \ref{thm:steady} as well as the first claim of Theorem \ref{thm:main}.

\section{The Memory-Deficient Regime}
\label{sec:memDeficient}
We now prove the second claim of Theorem \ref{thm:main} concerning the memory-deficient regime, where $\beta_m \ll 1/m$, as $m\to \infty$. The statement is repeated below for easy reference. 
\begin{theorem}
\label{thm:qBetaSlow}
Suppose that $\beta_m \ll 1/m$, as $m\to \infty$. Then, 
\begin{equation}
q_k = \lambda_{k}  \frac{\lambda_k\vee \alpha_0 + (K-1)\alpha_0}{\sum_{i \in \calK}\lt[ \lambda_i \vee \alpha_0 + (K-1)\alpha_0 \right]}, \quad k = 1, \ldots, K, 
\end{equation}
where $q_k = \lim_{m \to \infty} \E\lt(m^{-1} Q^m_k(\infty)\rt)$.
\end{theorem}

Fix $m\in \N$. Let $S^m_n$ be the $n$th update point,  and set $S^m_0 \bydef 0$. Denote by
\begin{equation}
C^m[n] = C^m(S^m_n), \quad n \in \zp, 
\label{eq:disctProc}
\end{equation}
 the embedded discrete-time process associated with the choice process $\{C^m(t)\}_{t\in \rp}$. The discrete-time processes $\{Q^m[n]\}_{n \in \zp}$ and $\{W^m[n]\}_{n \in \zp}$ are defined analogously. 

The key to proving Theorem \ref{thm:qBetaSlow} hinges upon the following property: the discrete-time choice process, $\{C^m[n]\}_{n \in \N}$, becomes {asymptotically Markovian} in the limit as $m$ tends to infinity. Note that $\{C^m[n]\}_{n \in \N}$ is not Markovian for a finite $m$, because the choices are sampled from a distribution that is a function of the current recallable reward vector, which in turn depends on past choices. However, the dependence of the reward vector on past choices is greatly weakened in the memory-deficient regime: because the life span of a unit reward is so small compared to the time between two adjacent choices, by the time a new choice is to be made, the agent will have essentially forgotten all of the past rewards associated with any action other than the one she currently uses. Therefore, when $m$ is large, knowing the choice $C^m[n]$ for some $n\in \N$, we can predict with high accuracy the recallable reward vector at the $n+1$ update point, rendering the choice process approximately Markovian. The approximate Markovian property will allow us to explicitly characterize the steady-state distribution of $C^m[\cdot]$ in the limit as $m\to \infty$, which we then leverage to analyze the steady-state distribution of the recallable reward process, $\{Q^m(t)\}_{t\in \N}$. Following this line of thinking, our proof will consist of three main parts: 
\begin{enumerate}
\item (Proposition \ref{prop:EQconditional}) Fix $n\in \N$. We first show that, in the limit as $m\to \infty$, the value of the scaled recallable reward process at the $(n+1)$-th update point converges in probability to a {deterministic} vector that only depends on $C^m[n]$. 
\item (Proposition \ref{prop:sitesSSProb}) Using the convergence of the recallable reward vector, we derive an explicit expression for the steady-state distribution of the discrete choice process, $\{C^m[n]\}_{n \in \N}$, in the limit as $m\to \infty$. 
\item Using the stationarity of the original recallable reward process, we combine the above two results to arrive at a characterization of the steady-state distribution of $\{Q^m(t)\}_{t\in \rp}$, thus completing the proof of Theorem \ref{thm:qBetaSlow}. 
\end{enumerate}

We begin by stating some basic properties of $Q^m[\cdot]$ and $C^m[\cdot]$, which will allow us to relate the steady-state behavior of these processes to their continuous-time counterparts; the proof of the following lemma is given in Appendix \ref{app:lem:pasta}. 

\begin{lemma} 
\label{lem:pasta}
Fix $m\in \N$. The discrete-time process $\{W^m[n]\}_{n \in \N}$ is positive recurrent, whose steady-state distribution satisfies
\begin{equation}
Q^m_k[\infty]  \preceq U^m_{\lambda_1}, \quad \forall k \in \calK, 
\label{eq:Qmdom-discrete}
\end{equation}
and
\begin{equation}
 C^m[\infty] \stackrel{{d}}{=} C^m(\infty), \quad \mbox{ and }  \quad Q^m[\infty] \stackrel{{d}}{=} Q^m(\infty). 
 \label{eq:Cm-discrete}
\end{equation}
\end{lemma}
 
For the remainder of the proof, we will fix $h: \zp\to \zp$ to be an increasing function, such that $1/\beta_m \gg h(m) \gg m$, and, in particular, 
\begin{equation}
\lim_{m \to \infty} h(m)/ m = \infty, \mbox{ and   } \lim_{m \to \infty} h(m)\beta_m = 0. 
\end{equation}
Define the scaled discrete-time recallable reward process:  $\qol[n] =  m^{-1}Q^m(S^m_n) = m^{-1}Q^m[n]$, $n \in \N$.  The following proposition states that when $m$ is large, the value of $\qol[1]$ converges to a deterministic value that depends solely on $C^m[0]$:  $\qol_k[1]$ is equal to $\lambda_k$ for $k = C^m[0]$, and zero, otherwise. In other words, the agent will have forgotten all past rewards other than those associated with her current choice. The proof is given in Appendix \ref{app:prop:EQconditional}. 
\begin{proposition} 
\label{prop:EQconditional}
Define the $K \times K$ matrix $q^*$, where  
\begin{equation}
q^*_{k,i} = \mathbb{I}(k=i)\lambda_{k}, \quad k,i\in \calK,
\end{equation} 
Suppose that the continuous-time process $\{W^m(t)\}_{t\in \N}$ is initialized  at $t=0$ according to its steady-state distribution. Then, for all $i, k \in \calK$,
\begin{equation}
\lim_{m \to \infty}\pb\lt( \big| \qol[1]- q^*_{k,i} \big| \leq \epsilon \bbar C^m[0]=k\rt)  = 1, \quad \forall \epsilon >0, 
\label{eq:convgProbQmh}
\end{equation}
and
\begin{equation}
\lim_{m \to \infty}\E\lt( \big| m^{-1} Q^m(h(m))- q^*_{k,i} \big| \bbar C^m[0] = k\rt) = 0. 
\label{eq:convgExpQmh}
\end{equation}
\end{proposition} 

Using Proposition \ref{prop:EQconditional}, we will be able to obtain an explicit expression for the limit $\lim_{m\to \infty} \pb(C^m[\infty] = k)$, as is stated in the following proposition.
\begin{proposition} 
\label{prop:sitesSSProb}
Define $p^m_k = \pb(C^m[\infty] = k)$, $k \in \calK.$ We have that
\begin{equation}
\lim_{m \to \infty} p^m_k =  \frac{\lambda_k\vee \alpha_0 + (K-1)\alpha_0}{z_0}, \quad k \in \calK, 
\label{eq:pmklim}
\end{equation}
where $z_0 = \sum_{i \in \calK}\lt[ \lambda_i \vee \alpha_0 + (K-1)\alpha_0 \right]$ is a normalizing constant. 
\end{proposition}
\bpf Suppose that $\{W^m[n]\}_{n \in \N}$ is initialized at $n=0$ in its steady-state distribution. We have that, for all $i\in \calK$, 
\begin{align}
\pb(C^m[1] = i) = \sum_{k\in \calK} \pb(C^m[0] = k)\sum_{u \in \calQ^m} \lt(\frac{u_i\vee \alpha_0}{\sum_{ \in \calK} (u_j\vee \alpha_0)} \pb\lt(\qol[1] = u \bbar C^m[0]=k\rt)\rt), 
\label{eq:Cm1decomp}
\end{align}
where  $\calQ^m \bydef \frac{1}{m}\zp^\calK$. 
By Eq.~\eqref{eq:convgProbQmh} in Proposition \ref{prop:EQconditional}, Eq.~\eqref{eq:Cm1decomp} implies that, there exists a sequence of $K\times 1$ vectors $\{\delta^m\}_{m \in \N}$ with $
\lim_{m \to \infty} \max_{i \in \calK} |\delta^m_{i}| = 0$,  such that for all $ i \in \calK$ and $m \in \N$, 
\begin{align}
\pb(C^m[1] = i) = & \delta^m_i+ \sum_{k\in \calK} \pb(C^m[0] = k) \frac{q^*_{k,i} \vee \alpha_0}{\sum_{j \in \calK}( q^*_{k,j} \vee \alpha_0)},
\label{eq:Cmtransition}
\end{align}
where $q^*_{k,i} = \mathbb{I}(k=i)\lambda_{k}$. Note that by the stationarity of $C^m[\cdot]$, we have that $C^m[0]\stackrel{d}{=}C^m[1]\stackrel{d}{=}C^m[\infty]$. Therefore, by treating $p^m = (p^m_1, p^m_2, \ldots, p^m_K)^T$ as a column vector, and defining the $K\times K$  matrix, $R$, where
\begin{equation}
R_{k,i} = \frac{q^*_{k,i} \vee \alpha_0}{\sum_{j \in \calK}( q^*_{k,j} \vee\alpha_0)}, \quad \forall i,k\in \calK, 
\end{equation}
Eq.~\eqref{eq:Cmtransition} can be written more compactly as
\begin{equation}
p_i^m = \delta^m_i + (R p^m)_i, \quad \forall i \in \calK.  
\label{eq:pmR1}
\end{equation}
It is not difficult to verify that  $R$  is row-stochastic and corresponds to the transition kernel of a discrete-time, irreducible  Markov chain over $\calK$. Therefore,  the equation  $ x = Rx $ admits a unique solution, $x^*$, in the $K$-dimensional simplex. Furthermore, it is not difficult to check that the following choice of $x^*$ satisfies $x^* = Rx^*$, and is hence the unique solution: 
\begin{equation}
x^*_k = \frac{\lambda_k\vee \alpha_0 + (K-1)\alpha_0}{\sum_{i \in \calK}\lt[ \lambda_i \vee \alpha_0 + (K-1)\alpha_0 \right]}, \quad k = 1,2, \ldots, K. 
\label{eq:pmR2}
\end{equation}
Since $\lim_{m \to \infty} \max_{i \in \calK} |\delta^m_{i}| = 0$, combining Eqs.~\eqref{eq:pmR1} and \eqref{eq:pmR2}, we conclude that $p^m$ converges to $x^*$ coordinate-wise, as $m\to \infty$.  This completes the proof of Proposition \ref{prop:sitesSSProb}. \qed

We are now ready to complete the proof of Theorem \ref{thm:qBetaSlow}. Fix $m\in \N$. Suppose that we initiate the continuous-time process $\{W^m(t)\}_{t\in \N}$ in its steady-state distribution so that it is stationary. The  stationarity implies that for all $k \in \calK$, 
\begin{align}
& \E(m^{-1}Q^m_k(0)) \nln
= & \E(m^{-1}Q^m_k(h(m))) = \sum_{i \in \calK}\pb(C^m(0)= i)\E\lt(m^{-1}Q^m_k(h(m)) \bbar C^m(0)= i\rt) \nln
\sk{a}{=} & \sum_{i\in \calK}\pb(C^n[\infty]=i) \E\lt(m^{-1}Q^m_k(h(m)) \bbar  C^m[0] = i\rt) \nln
=&   \sum_{i\in \calK}p^m_i\E\lt(m^{-1}Q^m_k(h(m)) \bbar  C^m[0] = i\rt),\nonumber
\end{align}
where $p^m_i$ is as defined in Proposition \ref{prop:sitesSSProb}, and step $(a)$ follows from the fact that $C^m[\infty] = C^m(\infty)$ (Lemma \ref{lem:pasta}).    Taking the limit as $m \to \infty$, we have that 
\begin{align}
q_k = & \lim_{m \to \infty}\E(m^{-1}Q^m_k(0))   {=}  \lim_{m \to \infty} \sum_{i\in \calK}p^m_i\E\lt(m^{-1}Q^m_k(h(m)) \bbar  C^m[0] = i\rt)\nln
\sk{a}{=}& \sum_{i\in \calK}\frac{\lambda_i \vee \alpha_0 + (K-1)\alpha_0}{z_0} q^*_{k.i}= \lambda_{k}  \frac{\lambda_k\vee \alpha_0 + (K-1)\alpha_0}{\sum_{i \in \calK}\lt[ \lambda_i \vee \alpha_0 + (K-1)\alpha_0 \right]}. 
\end{align}
where step $(a)$ follows from Eq.~\eqref{eq:convgExpQmh} in Proposition \ref{prop:EQconditional} and Eq.~\eqref{eq:pmklim} in Proposition \ref{prop:sitesSSProb}. This completes the proof of Theorem \ref{thm:qBetaSlow}. 

\section{Proof of Theorem \ref{thm:Ck}}
\label{sec:thmCk}
We prove Theorem \ref{thm:Ck} in this section. Let $W^m[\cdot] = (C^m[\cdot], Q^m[\cdot] )$ be the discrete-time embedded Markov chain defined in Eq.~\eqref{eq:disctProc}. By Lemma \ref{lem:pasta}, $C^m[\infty]$ and $C^m(\infty)$ admit the same distribution, and it hence suffices to characterize the former. Because the limiting expressions for $C^m[\infty]$  in the memory deficient regime have already been established in Proposition \ref{prop:sitesSSProb}, we will focus on the memory abundant regime, where $\beta_m \gg 1/m$. By Theorem \ref{thm:steady}, for all $\epsilon>0$
\begin{equation}
\lim_{m \to \infty}\pb(\|m^{-1}Q^m[\infty] - q^I \| > \epsilon) = \lim_{m \to \infty} \pb(\|m^{-1} Q^m(\infty) - q^I \| > \epsilon) = 0. 
\label{eq:scaledQmconverge}
\end{equation}
where the first step follows from Eq.~\eqref{eq:Cm-discrete} in Lemma \ref{lem:pasta}.  Fix $k \in \calK$. We have that
\begin{align}
\lim_{m\to \infty} \pb(C^m[\infty] = k )= \lim_{m\to \infty} \sum_{q \in \rp^\calK}\frac{q_k \vee (m \alpha_0)  }{\sum_{i}q_i\vee (m \alpha_0) } \pb(Q^m[\infty ] = q)
= & \frac{q^I_k \vee  \alpha_0  }{\sum_{i}q^I_i\vee \alpha_0 }, 
\label{eq:ckmconverge}
\end{align}
where the last equality follows from Eq.~\eqref{eq:scaledQmconverge} and the fact that $\frac{q_k \vee (m \alpha_0)  }{\sum_{i}q_i\vee (m \alpha_0) }$ is uniformly bounded over $q\in \rp^{\calK}$. Eq.~\eqref{eq:ckmconverge},  along with the expression for $q^I$ in Theorem \ref{thm:invStat}, leads to Eqs.~\eqref{eq:ckFast1} and \eqref{eq:ckFast2}, depending on the value of $k$ and $\alpha_0$. This completes the proof of Theorem \ref{thm:Ck}.

\section{Extensions and Generalizations}

We examine in this section a number of extensions to the original model. We will focus on key ideas with the intention of illustrating possible future directions, and the discussion will be more exploratory in nature, and less formal than that of our main results. 

\subsection{Fluid Solutions with Polynomial Choice Models} We consider in this sub-section a generalization of the reward-matching choice rule. Recall that under reward matching, the probability of choosing action $k$ when the recallable reward vector is $Q$ is {linearly} proportional to $Q_k \vee \alpha$. A natural generalization would be to consider a {polynomially weighted} reward-matching rule, where the probability of choosing action $k$ is proportional to $\lt( Q_k \vee \alpha\rt)^\eta$, where $\eta >0$ is a fixed parameter. 

An interesting question is how the agent's behavior will depend on the value of $\eta$, and  we will provide some preliminary analysis for this question in this sub-section. We will focus our attention on the invariant state(s) of the corresponding fluid solutions for the recallable rewards, as it serves as a good proxy for the stationary distribution of the pre-limit process in the memory-abundant regime, and from it the stationary choice distributions may also be inferred. The memory-deficient regime is in fact much easier to analyze, and we will postpone the discussion of that regime till the end of this sub-section (Section \ref{sec:genChoi}). 

In the fluid limit, the polynomial reward matching rule would lead to a modified function, $p$ (originally defined in Eq.~\eqref{eq:pdef}), given by 
\begin{equation}
p_k(q) = \frac{(q_k \vee \alpha_0)^\eta}{\sum_{i\in \calK}(q_i\vee \alpha_0)^\eta}, \quad q\in \rp^\calK, k \in \calK.
\label{eq:pdepoly} 
\end{equation}
It turns out that the dynamics of the fluid model will highly depend on whether $\eta \leq 1$ or not, and we will divide our analysis into two cases as such.  

{\bf Case 1: $\eta\leq 1$}. The next proposition shows that when $\eta< 1$ the fluid solutions admit a \emph{unique} invariant state, and gives a characterization of its form.\footnote{Our original model already covers the case of $\eta=1$, and we shall therefore focus on the case where $\eta<1$.} The proof of Proposition \ref{prop:etaless1} is given in Appendix \ref{app:prop:etaless1}. 
\begin{proposition}
\label{prop:etaless1}
Fix $\alpha_0 >0$ and $\eta \in (0,1)$.  The fluid solutions admit a unique invariant state, $q^I$, whose expression depends on the value of $\alpha_0$, as follows. 
\begin{enumerate}
\item  Suppose that $\alpha_0 > \lambda_1/K$. Then $q^I_k = \frac{\lambda_k}{K}$, for all $ k \in \calK.$ 
\item Suppose that $\alpha_0\leq \lambda_K \frac{\lambda_K^\frac{\eta}{1-\eta}}{\sum_{i \in \calK}
\lambda_i^\frac{\eta}{1-\eta}}$. Then
\begin{equation}
q^I_k = \lambda_k \frac{\lambda_k^{\frac{\eta}{1-\eta}}}{\sum_{i \in \calK}\lambda_i^{\frac{\eta}{1-\eta}}}, \quad \forall k \in \calK. 
\label{eq:qik-smeta}
\end{equation}
\item Suppose that $ \lambda_K \frac{\lambda_K^\frac{\eta}{1-\eta}}{\sum_{i \in \calK}
\lambda_i^\frac{\eta}{1-\eta}} <\alpha_0 \leq \lambda_1/K$. There exists $i^* \in \{1,\ldots, K-1\}$, such that $q^I_i \geq \alpha_0$ for $i = 1, \ldots, i^*$, and  $q^I_i < \alpha_0$ for $i=i^*+1, \ldots, K$. In particular, define the function 
\begin{equation}
g(i) =  \alpha_0\lt[ (K-i) + \sum_{j=1}^{i} \left( \frac{\lambda_j}{\lambda_i}\right)^{\frac{\eta}{1-\eta}} \rt],  \quad i =1, \ldots, K-1. 
\end{equation}
Then, $i^*$ is the unique index that satisfies  
\begin{equation}
\lambda_{i^*} \geq g(i^*), \quad \mbox{and}\quad \lambda_{i^*+1} < g(i^*).
\label{eq:gitoistar}
\end{equation}
Furthermore,  $q^I_{i^*}$ is the unique solution to the following equation in $\rp$
\begin{equation}
\lt( q^I_{i^*} \rt)^{1-\eta}= \frac{\lambda_{i^*}}{(K-i^*)\alpha_0^\eta + \lt(q^I_{i^*} \rt)^\eta \sum_{j=1}^{i^*} \left( \frac{\lambda_j}{\lambda_{i^*}}\right)^{\frac{\eta}{1-\eta}}},
\label{eq:qstarval}
\end{equation}
and the remaining coordinates of $q^I$ are given by: 
 \begin{align}
q^I_k = \left\{ \begin{array}{ll}
          q^I_{i^*} \left( \frac{\lambda_k}{\lambda_{i^*}}\right)^{\frac{1}{1-\eta}}, & \quad k=1, \ldots, i^*-1,\\
          \lt(q^I_{i^*} \rt)^{1-\eta}\frac{\lambda_k}{\lambda_{i^*}}\alpha_0^\eta , & \quad  k = i^*+1, \ldots, K. \\
         \end{array},  \right. \label{eq:thm-smalleta}
\end{align}
\end{enumerate}
\end{proposition}

\begin{figure}
\begin{minipage}{0.47\textwidth}
\centering
\includegraphics[height=4.5cm]{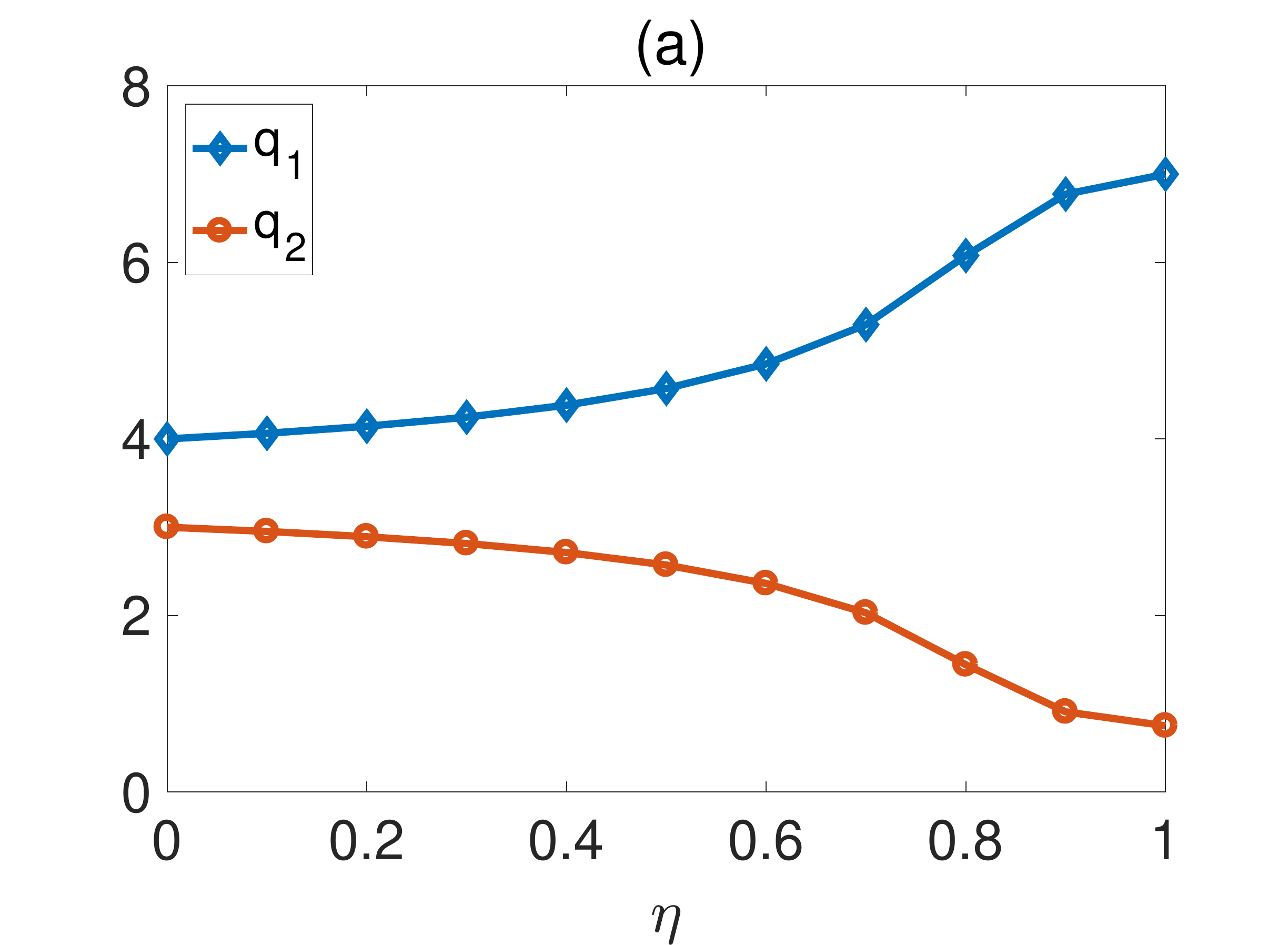}
\end{minipage}\hfill
\begin{minipage}{0.47\textwidth}
\centering
\includegraphics[height=4.5cm]{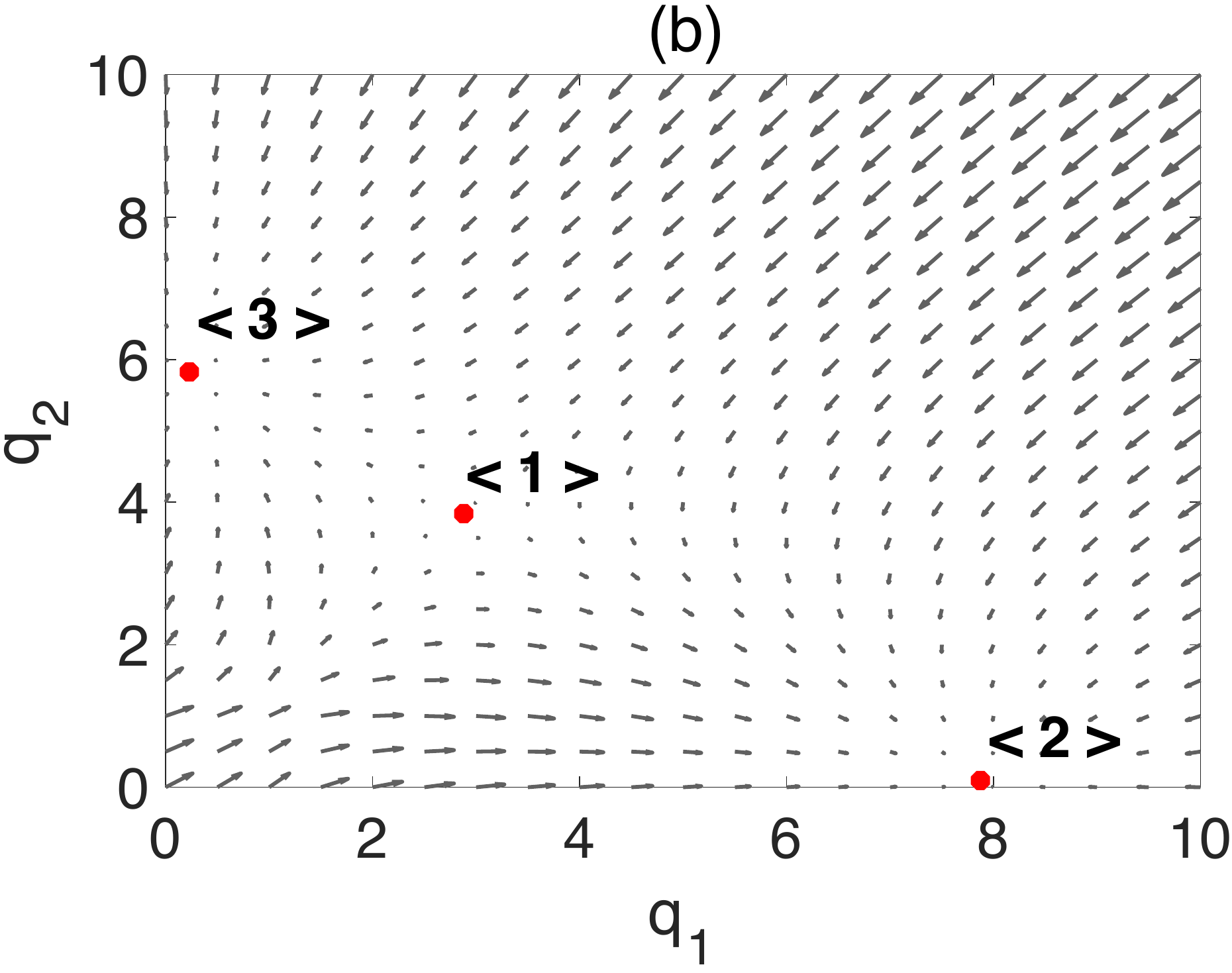}
\end{minipage}
\caption{An example for the polynomial reward-matching model with $K=2$, $(\lambda_1,\lambda_2) =(8, 6) $, and $\alpha_0=1$. Plot (a): the unique invariant state of the fluid solutions as a function of $\eta$, with $0\le \eta \le 1$. The agent tends to choose equally likely the two actions as $\eta \to 0$. (b) Vector field for the drift of the fluid solution, when $\eta =2$. The red dots mark the three invariant states. }
\label{fig:eta}
\vspace{-15pt}
\end{figure}

We can draw a few interesting observations from Proposition \ref{prop:etaless1}. The result suggests that the equilibrium behavior of the agent is largely consistent with the original model and our intuition: as $\eta$ decreases towards $0$, the agent's choice distribution tends to become uniform over the actions, since when $\eta$ is small, the agent's choice probabilities depend weakly on the rewards of the actions. This is illustrated numerically in Figure \ref{fig:eta}-(a) via an example with two actions. The uniqueness also suggests that the stochastic approximation portion of our results (Theorem \ref{thm:main}) should carry over to the regime of $\eta<1$ using the same arguments, without significant difficulty. Overall, these observations show that the insights from our original model should be fairly robust under small perturbation of $\eta$ in the direction towards $0$. 

Proposition \ref{prop:etaless1} also reveals interesting new insights that are less obvious. Firstly, Item 1 of the proposition shows that the phenomenon of ``complete oblivion'' observed in the original model (see Theorem \ref{thm:main} and Section \ref{sec:impThm}), where the agent chooses uniformly across actions as a result of a large exploration parameter $\alpha_0$, is not an artifact of $\eta$ being exactly equal to $1$, and in fact continues to exist \emph{whenever} $\eta<1$. Moreover, the threshold on $\alpha_0$ that marks the beginning of such oblivion is the same as before, and does not depend on $\eta$. 

Secondly, Item 2 of the proposition shows a regime not previously observed, and with a surprisingly clean expression: when $\alpha_0$ is very small, the probability of choosing action $k$ in the invariant state is {exactly proportional} to $\lambda_k^{\frac{\eta}{1-\eta}}$, and this expressions, remarkably, does not depend on $\alpha_0$ (to be contrasted with Eq.~\eqref{eq:ckFast1} which does depend on $\alpha_0$). This result would appear to lead to a contradiction with Theorem \ref{thm:main}: by increasing $\eta$ towards $1$, Eq.~\eqref{eq:qik-smeta} suggests that the agent will eventually concentrate $100$\% of her efforts on the best action and the degree of such concentration is independent of $\alpha_0$. In other words, how is it possible that the agent may choose the optimal action more frequently when $\eta<1$ than when $\eta=1$? To resolve this puzzle, we note that the regime in Item 2 requires  $\alpha_0\leq \lambda_K \frac{\lambda_K^\frac{\eta}{1-\eta}}{\sum_{i \in \calK}
\lambda_i^\frac{\eta}{1-\eta}}$, and the right-hand side, which is also equal to the invariant state of the worst action, $q^I_K$,  tends to $0$ as $\eta\to 1$. Therefore, as $\eta$ approach $1$ from below, the regime in Item 2 remains valid only if $\alpha_0$ vanishes accordingly, thus resolving the paradox. 

There is another important consequence of Item 2. Suppose that $\alpha_0$ is fixed at a sufficiently small value, say, $\frac{1}{10}\cdot \frac{\lambda_K^2}{\sum_{i}\lambda_i}$. Then, Item 2 becomes the only valid regime for \emph{all} sufficiently small $\eta$. As such, Eq.~\eqref{eq:qik-smeta} provides a precise, quantitative description as to \emph{how} the agent's choice probabilities become uniform, as $\eta\to 0$, namely, that the probability of choosing action $k$ will be equal to $\frac{\lambda_k^{\frac{\eta}{1-\eta}}}{\sum_{i \in \calK}\lambda_i^{\frac{\eta}{1-\eta}}}$ for all sufficiently small $\eta$.

Finally, Item 3 of the proposition addresses the remaining scenarios, where $\alpha_0$ is neither too small nor too big. Unfortunately, we no longer have a closed-form expression for this intermediate regime. Nevertheless, the invariant state can still be evaluated without much difficulty, and all coordinates of $q^I$ can be derived in closed-form as a function of $i^*$ and $q^I_{i^*}$.

{\bf Case 2: $\eta>1$.} We now look at the second case, where $\eta>1$. The dynamics of the fluid solutions turn out to be more complex than when $\eta<1$, and we do not yet have a good understanding. To illustrate the complexities that one would encounter in this regime, we next show a simple example demonstrating that the fluid solutions could admit \emph{multiple} invariant states. 

Consider the a setting with two actions ($K=2$), where $\eta >1$ and  $\alpha_0< \lambda_1^{\frac{1}{1-\eta}}/(\lambda_1^{\frac{\eta}{1-\eta}}+ \lambda_2^{\frac{\eta}{1-\eta}})$. The drifts of the fluid solutions for this example are visualized in Figure \ref{fig:eta}-(b).  The fluid solutions admit at least three distinct invariant states, marked by the red dots in Figure \ref{fig:eta}-(b). First, one can verify that the vector $q^{(1)}$, where $q^{(1)}_1 = \frac{\lambda_1^{\frac{1}{1-\eta}}}{\lambda_1^{\frac{\eta}{1-\eta}}+ \lambda_2^{\frac{\eta}{1-\eta}}}$ and  $q^{(1)}_2 = \frac{\lambda_2^{\frac{1}{1-\eta}}}{\lambda_1^{\frac{\eta}{1-\eta}}+ \lambda_2^{\frac{\eta}{1-\eta}}},$
constitutes an invariant state (labeled $\lt<1\rt>$ in Figure \ref{fig:eta}-(b)). However, $q^{(1)}$ is not the unique invariant state, and nor is it stable: a perturbation along the direction of $(\epsilon, -\epsilon)$ or $(-\epsilon, \epsilon)$, where $\epsilon>0$ is a small constant, can induce the fluid solution to move towards one of the other two (stable) invariant states, $q^{(2)}$ (labeled $ \lt< 2 \rt>$ in Figure \ref{fig:eta}-(b)) and $q^{(3)}$ (labeled $ \lt<3 \rt>$ in Figure \ref{fig:eta}-(b)), respectively. Here, $q^{(2)}$ satisfies $q^{(2)}_1 > \alpha_0> q^{(2)}_2$, $\lambda_1 \frac{(q^{(2)}_1)^\eta}{(q^{(2)}_1)^\eta + \alpha_0^\eta} =q^{(2)}_1$, and $q^{(2)}_2 = \lambda_2 \frac{\alpha_0^\eta}{(q^{(2)}_1)^\eta + \alpha_0^\eta}$. Analogous to $q^{(2)}$, the state $q^{(3)}$ satisfies $q^{(3)}_1 < \alpha_0 < q^{(3)}_2$, $\lambda_2 \frac{(q^{(3)}_2)^\eta}{(q^{(3)}_2)^\eta + \alpha_0^\eta} =q^{(3)}_2$ and $q^{(3)}_1 = \lambda_1 \frac{\alpha_0^\eta}{(q^{(3)}_2)^\eta + \alpha_0^\eta}$.

We offer some speculation as to why the multiplicity of invariant states emerges when $\eta>1$. In this case, the agent leans disproportionally more towards actions with higher rewards, which means that sufficiently substantial advantages in the initial recallable rewards $q(0)$ could lead an action $i$ to permanently dominate another action $j$, even if $\lambda_i < \lambda_j$. This logic suggests that, in principle, if the reward rates do not vary significantly across actions, then any measure that concentrates on a single action could be a stable invariant state for the choice process, so long as the fluid solution starts in a state where the action has a substantial amount of initial rewards compared to other actions; this seems to explain the two stable invariant states in Figure \ref{fig:eta}-(b), $\lt<2\rt>$ and $\lt<3\rt>$. The unstable invariant state $\lt< 1\rt>$, on the other hand, seems to stem from a different type of dynamics: it is a point where the action with the superior reward rate (in this case, $\lambda_1$) has fewer recallable rewards ($q_1(t)$), and the disadvantage in recallable rewards is exactly naturalized by the advantage in reward rate. The point is unstable because small perturbations could easily break this balance, inducing the advantage in the reward rate to dominate that of recallable rewards or \emph{vice versa}. 

What does this mean for the pre-limit, stochastic system? Note that while the fluid solutions can admit multiple invariant states when $\eta>1$, the steady-state distribution for any pre-limit stochastic system ($m<\infty$) is always unique: the exponential decay of recallable rewards in the pre-limit system ensures that $Q(\cdot)$ is positive recurrent, as it will return infinitely often to the state $0$. Nevertheless, the proceeding analysis of the fluid solutions suggests that the stochastic system will likely experience a qualitative change as well once $\eta$ is greater than $1$: when $m$ is large, instead of converging over time to the neighborhood of a single configuration, as is the case when $\eta\leq 1$, $Q(\cdot)$ is likely to alternate between multiple different configurations, each corresponding to one of the multiple invariant states in the fluid solutions.  

\subsubsection{General Choice Models} 
\label{sec:genChoi}
At a higher level, the reward matching rule is a  special case of a broader family of choice heuristics, in which action $k$ is chosen with a probability proportional to $w(Q_k(t) \vee \alpha)$, where $w: \rp \to \rp$ is a non-decreasing \emph{weight function}. The polynomial reward-matching rule studied earlier corresponds to $w(x) = x^\eta$, for some $\eta>0$. One may also consider $w(x) = \exp(c x)$ for some constant $c>0$, the sampling procedure becomes reminiscent of the celebrated logit model in choice theory (cf.~Chapter 3, \cite{train2009discrete}) as well as the exponential-weight algorithm in the multi-arm bandit literature (cf.~\cite{auer2002nonstochastic}). 

Can we extend our results beyond the polynomial reward-matching rule to include this even broader family of choice models? On the positive side, it is not difficult to check that the same analysis for the \emph{memory-deficient regime} can be readily extended to these models: the limiting probability of choosing action $k$ would become proportional to $w(\lambda_k)\vee \alpha_0 + (K-1)\alpha_0$, i.e., with $\lambda_k$ in our original result being replaced by $w(\lambda_k)$. Note that if we had the freedom of choosing $w(\cdot)$, then applying a highly skewed weight function, e.g., $w(x) = 10^{20^x}$, would induce the agent to almost always choose the best action in steady state, even in the memory-deficient regime. However, the downside is that with a skewed weight function, the agent is highly unlikely to give up the current choice at an update point, regardless of its reward rate. It will thus take a long time for the agent's time-average behavior to approach that of the steady state, resulting in bad transient performance. A similar trade-off seems to be present in the context of the exploration parameters, $\alpha_0$, where a smaller $\alpha_0$ tends to improve steady-state performance at the expense of a slower convergence to the steady-state dynamics. It would be interesting to obtain a more precious and rigorous understanding of the trade-off between the agent's steady-state  versus transient performance. 

On the negative side, extensions to general choice models appear to be significantly more difficult in the memory-abundant regime, as is evidenced by the proceeding analysis of the polynomial reward-matching model. Making progress on this front is likely to require more sophisticated analysis to cope with the complex dynamics that will arise in fluid solutions, such as when $\eta>1$ in the polynomial reward-matching rule.

\begin{figure}
\begin{minipage}{0.47\textwidth}
\centering
\includegraphics[scale=.24]{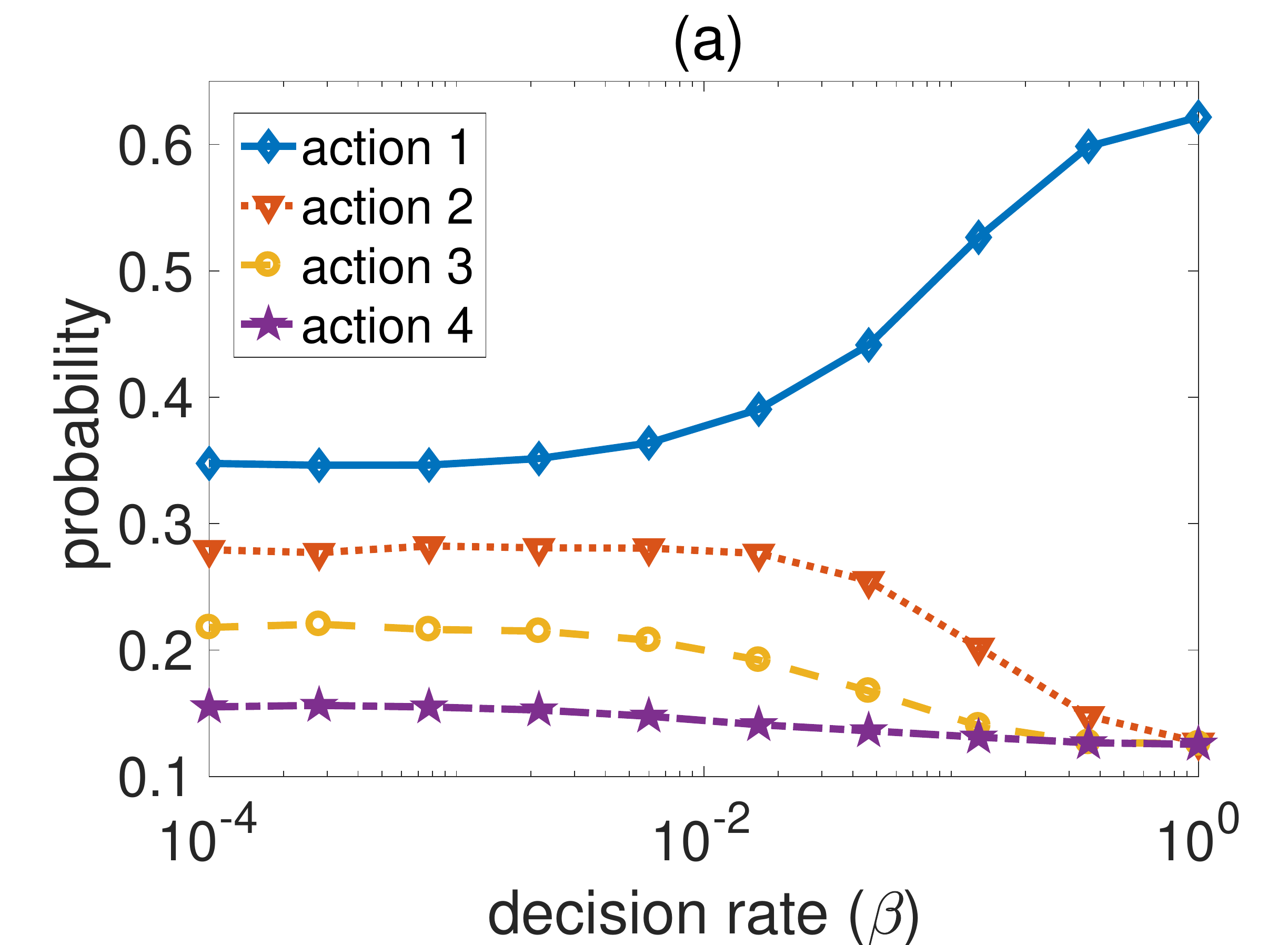}
\end{minipage}\hfill
\begin{minipage}{0.47\textwidth}
\centering
\includegraphics[scale=.24]{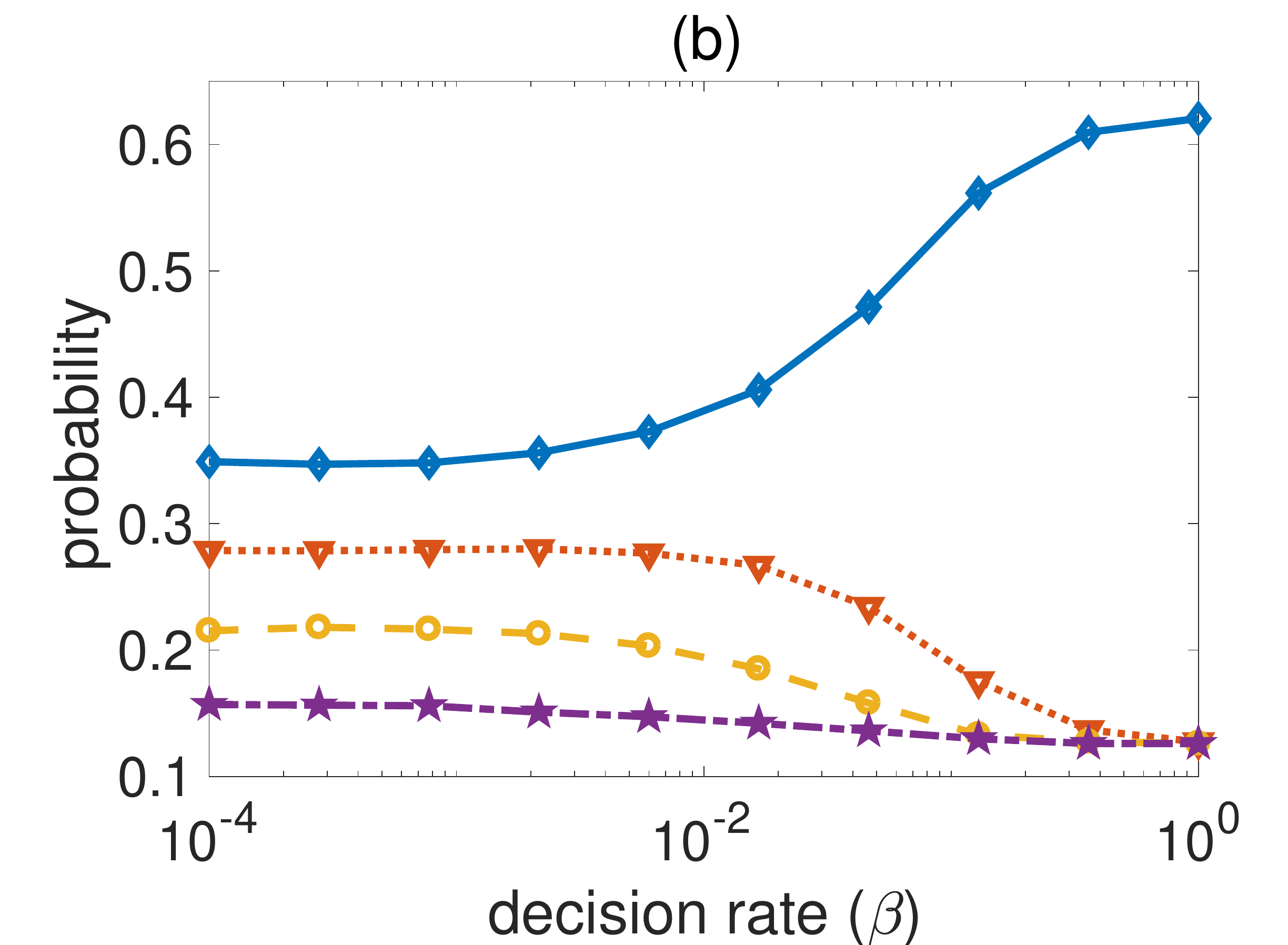}
\end{minipage}
\caption{Distribution of choices as a function of the update rate, $\beta$, with $K=4$, $(\lambda_1,\lambda_2, \lambda_3, \lambda_4) =(8, 6, 4, 2) $, $\alpha_0=1$, and $m=200$. In $(a)$, the lifespans of rewards are set to the constant $m$. In plot $(b)$ the lifespans follow  a Pareto distribution with scale ${m}/{2}$ and shape $2$. Each point is averaged over $5 \times 10^7$ to $5\times 10^5$ time units for $\beta$ ranging from $10^{-4}$ to $1$, respectively (longer duration for a slower update rate). }
\label{fig:const-levy}
\vspace{-15pt}
\end{figure}

\subsection{Non-Exponential Lifespan Distributions}
We have thus far assumed that the lifespans of rewards follow an exponential distribution, and it would be interesting to see whether the qualitative insights from our model would continue to hold under other, non-exponential lifespan distributions. In this subsection, we experiment with two other lifespan distributions: in the first case, we set the lifespan to be a constant, $m$, and in the second case, the lifespan follows a (heavy-tailed) Pareto distribution with scale $m/2$, shape $2$, and mean $m$. Figure~\ref{fig:const-levy} shows the steady-state distribution of choices under these two scenarios, with $m=200$. Interestingly, comparing the results in Figure \ref{fig:const-levy} to those from original model (Figure \ref{fig:SS-Sim}-(b)), we see that the distribution of choices is largely insensitive to the choice of lifespan distribution, and the theoretical predictions made in Theorem \ref{thm:main} appear to hold even under these non-exponential distributions. 

\subsection{Heterogeneous Memory Decay Rates}
\label{sec:heterorate}

As was alluded to in Example 2 in the Introduction, we may also consider the case where the memory decay rates are not constant across actions, so that the rewards associated with action $k$ depart at rate $\mu_k$. This generalization is useful for modeling service system where the rate of customer departure or attrition depends on the service type a customer is associated with. 

Analogous to the case of uniform decay rate, the scaling of Section \ref{sec:metrics} corresponds to setting $\mu = \mu^0_k m^{-1}$ in the $m$th system,  where $\{\mu^0_k\}_{k \in \calK}$ are positive constants.
Under this scaling, it is not difficult to verify that both Theorems \ref{thm:Ck} and \ref{thm:main} extend to this more general case, if we replace every occurrence of $\lambda_k$ with a decay-rate-weighted version, $\lambda_k / {\mu^0_k}$. In other words, as far as the limiting system is concerned, the version of the problem with heterogeneous memory decay rates is equivalent to one with uniform decay rates but appropriately weighted reward rates. We should note that this equivalence only holds in the limit as $m\to 0$, and is not exact in a  pre-limit system with a finite $m$.

\section{Discussion}
\label{sec:extensions}
We have proposed in this paper a stochastic action-reward model for studying the impact of imperfect memory in dynamic decision making. In the limiting regimes where the agent's memory span is large, our main results provide exact expressions for the steady-state choice probabilities of an agent following the so-called reward-matching rule, and demonstrate that these probabilities are highly sensitive to the relative scaling between two parameters: the rate of memory decay, $\mu$, and the rate at which the agent makes new choices, $\beta$. 

There is a number of potentially interesting extensions  of our formulation. For instance, the current paper focuses on the regime where  the rate of memory decay $\mu \to 0$, and alternatively, one may also look at an opposite limiting regime where $\mu$ tends to infinity, i.e., the agent becomes extremely forgetful and her memory span approaches zero. Note that since the total amount of recallable rewards in system is on the order of $\mu^{-1}$ (regardless of the value of $\beta$), in this limit we should see very few recallable rewards at any point in time and the process of recallable rewards frequently hitting the zero state. It is therefore reasonable to assume that the reward-matching rule will perform poorly when $\mu$ is very large, with the agent choosing all actions essentially uniformly at random. Unfortunately, we do not yet have the tools needed to characterize this regime, since the fluid approximations developed in this paper heavily rely on mean-field-type evolution of a large amount of recallable rewards and are unlikely to be applicable when $\mu$ is large.

The present reward-matching model is {symmetric}, in the sense that the agent does not have any inherent preferences of one action over another. One way to relax this assumption is to allow the exploration parameter $\alpha$ to depend on the action, so that  action $i$ will be sampled with probability $\propto (Q_k (t)\vee \alpha^{(k)})$. In the memory-deficient regime, it would not be difficult to show that this results in a steady-state choice distribution where action $k$ is chosen with probability proportional to $\lambda_k\vee \alpha_0^{(k)} + \sum_{i\neq k}\alpha_0^{(k)}$ . The memory-abundant regime will become more challenging to analyze, because the heterogeneity of $\alpha^{(k)}$ would in general break the monotonicity of the invariant-state recallable rewards, since now an action with a higher reward rate may not necessarily have a higher recallable reward in the invariant state. It would be an interesting future direction to better understand how to incorporate inherent preference asymmetries into the formulation. 

While the present paper focuses largely on the theory, there appears to be several interesting potential applications of our model. First, the reward matching rule can be viewed as a simplified cognitive model to capture humans' unconscious reinforcement behavior, and it will be interesting to see if the theory could be used to make predictions on a consumer's long-run choice behavior when switching among similar products or services (eg., Example 1 of Section \ref{sec:intro}), or to investigate the effect of memory loss on a player's learning and performance in repeated games similar to those studied by \cite{erev1998predicting}.  Second,  memory loss can be interpreted in a more metaphorical sense, and be used to model the departure of customers (eg., Example 2 of Section \ref{sec:intro}). Cast in this light, the reward matching policy may act as a simple and intuitive heuristic in a manager's algorithmic toolbox for dynamically choosing product offerings or demographic targets for advertisements. Finally, the dynamics of the recallable rewards in our system resemble that of a queueing network. For instance, the departures of rewards for each individual action function similarly to service completions in an infinite-server queue. There has been a growing literature in recent years on queuing systems where customers make their own choices regarding abandonments or the types of services they would like to receive (cf.~\cite{hassin2003queue, pender2016managing, ding2016asymptotic, dong2018impact}), and it would be interesting to see whether our model can be applied to the study of (strategic) choice behavior in these service systems. 

While memory loss is regarded as a given constraint in the current work, one could also ask whether it would be beneficial to artificially induce memory loss even when perfect recall is feasible. 
For instance, it has been observed in the  multi-arm bandit literature that policies with built-in regularization that penalizes distant past experience can achieve better regret with time-varying, or adversarially chosen, parameters (cf.~\cite{garivier2008upper,van2014follow,besbes2015non,keskin2016chasing}). This is because these policies tend to be more adaptive to the environment and not ``weighed down'' by past experience. While our theorems do not directly apply to time-varying reward rates, Theorem \ref{thm:fluidToInvar} shows that the fluid solution converges exponentially quickly to the invariant state from any bounded initial condition, suggesting that the reward matching rule could be an attractive algorithm for dynamic learning in non-stationary environments.

Finally, at a higher level, there appear to be other interesting  ramifications of ``memory loss'' in  organizations to be explored. For instance, departures of employees could lead to the loss of skills and expertise (cf.~\cite{benkard2000learning}), and such organizational forgetting could potentially affect a firm's efficiency in a significant way.

\section{Acknowledgment} The authors would like to thank the anonymous referees at Mathematics of Operations Research for their comments and input, and Professors Steven Callander, Dana Foarta, J.~Michael Harrison, David M.~Kreps, Mihalis Markakis, Michael Ostrovsky, Daniel Russo, Daniela Saban and Takuo Sugaya for the insightful discussions and feedback on the manuscript.  Yun would like to acknowledge the support from the National Research Foundation of Korea (NRF) grant No.~2019028324, funded by the government of Korea (MSIT).

\bibliographystyle{apalike}
\bibliography{rfm.bib} 

\appendix
\normalsize

\section{Comparison to Perfect Memory}
\label{app:perfect_memory}

We discuss in this appendix what could happen if there were no memory decay in our model, i.e., if $\mu=0$, and why it is different from the limiting regime considered in this paper, with $\mu \to 0$. 

When $\mu=0$, all recallable rewards  remain in the system indefinitely. If we view the recallable rewards sampled at the update points as a discrete-time process, and in addition set the exploration parameter $\alpha$ to 0, then our model becomes essentially the same as the choice process analyzed by \cite{beggs2005convergence}, where it is shown that the probability of choosing the best action converges to one as $t\to \infty$ under Luce's rule, as long as each action is associated with some strictly positive initial rewards. If $\alpha$ is a positive constant, because there is no memory decay, the effect of $\alpha$ disappears  as soon as the rewards for all actions exceed $\alpha$, and the same conclusion should hold.  

Therefore, one would expect that when $\mu=0$ and $\alpha \geq 0$, the choice probability under the reward-matching rule will concentrate on the best action as $t\to \infty$, regardless of the update rate, $\beta$. This is however different from the conclusion of Theorem \ref{thm:Ck}, which shows the existence of two distinct limiting steady-state probabilities, one of which does not exhibit concentration on the best action.  These observations thus suggest that our scaling regime do capture unique effects of {imperfect} memory. This is perhaps not too surprising in hindsight: if we had set $\mu$ to zero, any positive update rate $\beta$ would become, by definition, significantly greater than $\mu$, and hence the second (memory-deficient) regime in Theorems \ref{thm:Ck} and \ref{thm:main} could not have appeared when $\mu=0$.

\section{Technical Preliminaries}
\label{app:technicaLemma}

\begin{proposition}[Doob's inequality]  (cf.~Section 12.6, \cite{grimmett2001probability})
\label{prop:doob}
Let $\{X(t)\}_{t \geq 0}$ be a discrete- or continuous-time non-negative submartingale. Fix $T\geq 0$. We have that
$$\mathbb{P}\left(\sup_{0\le t\le T} X(t) \ge c \right) \le \frac{\E\left(X(T)\right)}{c}, \quad \forall c > 0. $$
\end{proposition}

\begin{proposition}[Gronwall's lemma] \label{prop:gronwall}
(cf.~Section 1.3, \cite{ames1997inequalities})
Let $f,b: \rp \to \rp$ be continuous and non-negative functions, and let $a: \rp \to \rp$ be a continuous, positive and non-decreasing function. If $f(t) \le a(t) + \int_0^t b(s)f(s) ds$, for all $t \in [0, T],$ then 
$$f(T) \le a(T) \exp\left(  \int_0^T b(t) \, dt \right).$$
\end{proposition}

\section{Proofs of Propositions}

\subsection{Proof of Proposition \ref{prop:bd1}}
\label{sec:proof-lemma-refl}

The proof of Proposition \ref{prop:bd1} is largely based on maximal inequalities for continuous-time martingales, and the main steps follow the approach developed in \cite{kurtz1978strong}. Define $X^m(t) \in \rp^\calK$, where
\begin{equation}
X_k^m(t) = \qol_k(t)- \lt( \qol_k(0) + \int_0^t G_k(\qol(s),\col(s)) ds \rt), \quad t\in \rp, k \in \calK.
\end{equation}
The proof proceeds in two stages: we first show, via Gronwall's lemma (Proposition \ref{prop:gronwall} in Appendix \ref{app:technicaLemma}), that bounding the deviation of $\qol_k(\cdot)$ from $V^m_k(\cdot)$ can be reduced to bounding the magnitude  of $X^m_k(\cdot)$. We then show that this can be accomplished by writing $X^m_k(\cdot)$ as a sum of two continuous-time martingales.   From the definition of $\{V^m(t) \}_{t\ge 0}$, we have that
\begin{align}
& \left\|\qol (t) - V^m(t) \right\|  \nln
=& \sum_{k\in \calK}\left|\qol_k (t) - q^0_k - \int_0^t G_k(V^m(s),\col(s)) ds \right| \cr
 \le &  \sum_{k\in \calK} \left|\qol_k(0) - q^0_k \right| + \sum_{k\in \calK}\left|X_k^m(t) + \int_0^t \left( G_k(\qol(s),\col(s))- G_k(V^m(s),\col(s)) \right)ds\right| \cr
\le & \left(\left\|\qol(0) - q^0 \right\| + \left\|X^m(t)\right\|\right) + \int_0^t \sum_{k\in \calK}\left| G_k(\qol(s),\col(s))- G_k(V^m(s),\col(s)) \right|ds \cr
\sk{a}{\le} & \left(\left\|\qol(0) - q^0 \right\| + \left\|X^m(t)\right\|\right) + \int_0^t \sum_{k\in \calK} \left| \qol_k(s)- V^m_k(s) \right|ds \nln
= & \left(\left\|\qol(0) - q^0 \right\| + \left\|X^m(t)\right\|\right)  + \int_0^t \left\| \qol(s)- V^m(s)\rt\|ds. \label{eq:Qv}
\end{align}
For step $(a)$, we observe that for all $c, k\in \calK$, $G_k(\cdot, c)$ is a $1$-Lipschitz continuous function. We now apply Gronwall's lemma (Proposition \ref{prop:gronwall} in Appendix \ref{app:technicaLemma}) to Eq.~\eqref{eq:Qv}, with $a(t) = \left\|\qol(0) - q^0 \right\| + \left\|X^m(t)\right\|$, $b(t)=1$, and  $f(t)=\left\|\qol (t) - V^m(t) \right\|$, and obtain that
\begin{align}
& \pb\lt( \sup_{0 \le t \le T} \left\|\qol (t) - V^m(t) \right\| > \epsilon \rt) \nln
 \leq & \pb\lt(   \sup _{0 \leq t\leq T} \left(\left\|\qol(0) - q^0 \right\| + \left\|X^m(t)\right\|\right)  e^t > \epsilon \rt) \nln 
\leq & \pb\lt(  \left\|\qol(0) - q^0 \right\| > \epsilon e^{-T}/2 \rt) + \pb\lt( \sup_{0\le t \le T} \left\|X^m(t) \right\| > \epsilon e^{-T}/2 \rt) .
\end{align} 
Because we have assumed that $\lim_{m\to \infty} \pb \left(\left\|\qol(0) - q^0 \right\| >\epsilon \right) = 0$ for all $\epsilon>0$, to prove Lemma~\ref{prop:bd1}, it therefore suffices to show  that
\begin{equation}
\lim_{m\to \infty} \pb\left( \sup_{0\le t \le T}  \left\|X^m (t)\right\| > \epsilon \right) =0 , \quad \forall \epsilon>0. 
\label{eq:doops-1}
\end{equation}

We now prove Eq.~\eqref{eq:doops-1} by expressing $X^m_k(\cdot)$ as the sum of two continuous-time martingales corresponding to the arrivals and departures of rewards, respectively. Fix $t\in \rp$ and $k\in \calK$, and denote by $A^m_k (t)$ and $D^m_k (t)$ the amount of rewards associated with action $k$ that have arrived and departed, respectively, during the interval $[0,t]$. In particular, we can write
\begin{equation}
Q^m_k(t) = Q^m_k(0)+ A^m_k(t) - D^m_k(t), \quad t\in \rp. 
\end{equation}
 We have that
\begin{align}
& \left| X_k^m(t) \right|  \nln
= & \lt|\frac{1}{m}(A^m_k(mt) - D^m_k(mt)) -  \int_0^t G_k(\qol(s),\col(s)) ds \rt| \nln
=& \left| \frac{1}{m}\left(A^m_k (mt) - \int_{0}^{mt} \lambda_k \mathbb{I}(C^m (s)=k)  ds \right) - \frac{1}{m}\left(D^m_k (mt) - \int_{0}^{mt}Q^m_k(s) m^{-1}  ds \right) \right|\cr
\le& \left| \frac{1}{m}\left(A^m_k (mt) - \int_{0}^{mt} \lambda_k \mathbb{I}(C^m (s)=k) ds \right) \right| + \left| \frac{1}{m}\left(D^m_k (mt) - \int_{0}^{mt} Q^m_k(s)  m^{-1} ds \right) \right|.\label{eq:xk}
\end{align}

For the remainder of the proof, we will focus on showing that, for all $\epsilon>0$, 
\begin{align}
\label{eq:X-1}
\lim_{m\to \infty} \pb\left( \sup_{0\le t \le T}  \left|\frac{1}{m} \left(A^m_k(mt)-\int_{0}^{mt}\lambda_k \mathbb{I}(C^m (s)=k) ds \right) \right| > \epsilon\right) =& 0,\quad \mbox{and}  \\ 
\lim_{m\to \infty} \pb\left( \sup_{0 \le t \le T} \left|\frac{1}{m}\left(D^m_k (mt) - \int_{0}^{mt}Q^m_k(s)m^{-1}  ds \right)\right| > \epsilon \right)  = & 0. \label{eq:X-2}
\end{align}
In light of Eq.~\eqref{eq:xk}, the above two equations  imply the validity of Eq.~\eqref{eq:doops-1}, which proves Proposition \ref{prop:bd1}.  The proofs for both Eqs.~\eqref{eq:X-1} and \eqref{eq:X-2} hinge upon the following technical result, which states that the sample path of a certain time-inhomogenous Poisson process stays close to its mean, with high probability. Note that a similar result for uniform-rate Poisson processes can be found in Lemma 7.6 of \cite{massey1998uniform}. The proof of the lemma is given in Appendix \ref{app:lem:poisson}, and uses a  similar line of argument as that of Theorem 2.2 in \cite{kurtz1978strong}. 

\begin{lemma}
\label{lem:poisson}
Fix $l \in \Z ^{K+1}$ and $T \in \rp$.  Let $\{ N(t) \}_{t\ge 0}$ be the counting process where $N(t) $ is the number of times in $[0,t]$ that the process $W^m(\cdot)$ jumps from state $w$ to $w+l$, for some $w \in \Z^{K+1}$. Denote by $\psi: \Z^{K+1} \to \rp$ the corresponding rate function of $N(\cdot)$, so that the instantaneous transition rate of $N(\cdot)$ at time $t$ is equal to $\psi(W^m(t))$. For all $\epsilon, \phi>0$, we have that
\begin{equation}
\pb\left( \sup_{0\le t \le T} \left|N(t)-\int_{0}^t {\psi} (W^m(s)) ds \right| \ge \epsilon \right) \le 2e^{-\phi T\cdot h(\epsilon / \phi T) }+\pb \left( \sup_{0\le t \le T}  \psi(W^m(t)) \ge \phi \right), \nnb
\end{equation}
where $h(x) = (1+x)\log (1+x)-x$. 
\end{lemma}

We now prove Eqs.~\eqref{eq:X-1} and \eqref{eq:X-2}. In the context of Lemma \ref{lem:poisson}, $A_k(\cdot)$ is the counting process with $l = e_k$, where $e_k$ is the vector $(K+1)\times 1$ vector whose $k$th coordinate is $1$ and all other coordinates zero, corresponding to an arrival to $Q^m_k(\cdot)$. The rate of $A_k(\cdot)$ at time $t$ is equal to $\lambda_k \mathbb{I}(C^m (t)=k)$, which is bounded from above by $\lambda_k$ for all $t\in \rp$. By applying Lemma \ref{lem:poisson}, with $\psi(W^m(t)) = \lambda_k \mathbb{I}(C^m (t)=k)$ and $\phi = \lambda_k$, we have that $\pb \left( \sup_{0\le t \le T} \psi(W^m(t)) > \phi \right) = 0$, and for all $\epsilon>0$, 
\begin{align}
&\lim_{m \to \infty} \mathbb{P}\left( \sup_{0\le t \le T}  \left|\frac{1}{m} \lt(A^m_k(mt)-\int_{0}^{mt}\lambda_k \mathbb{I}(C^m (s)=k) ds \rt) \right| \ge \epsilon \right) \nln
= & \lim_{m \to \infty} \mathbb{P}\left( \sup_{0\le t \le mT}  \left| A^m_k(t)-\int_{0}^{t}\lambda_k \mathbb{I}(C^m (s)=k) ds  \right| \ge  m \epsilon \right) \nln
\leq& \lim_{m \to \infty}  2\exp\left( -\lambda_k m T \cdot h\left( \frac{ \epsilon}{\lambda_k T}\right)\right)=0. 
\label{eq:r1-arriving}
\end{align}

The proof of Eq.~\eqref{eq:X-2} uses essentially the same idea, but the argument needs to be more delicate due to the fact that the rate of the counting process $D^m_k(t)$, which is equal to $Q^m_k (t) \mu = Q^m_k(t) m^{-1}$, is not bounded over the state space. Therefore, we first derive an upper bound on the tail probabilities of $Q^m_k (t)$, as follows. Let $\psi(\cdot)$ be the rate function for $D_k(\cdot)$ as in Lemma \ref{lem:poisson}, and $\phi = 2q^0_k + (1+\epsilon)\lambda_k T $. Fixing $\epsilon>0$, we have that, for all $r \in \zp$, 
\begin{align}
& \pb\left( \sup_{0\le t \le mT}  \psi(W^m(t)) \ge \phi \right)  =  \pb \left(\sup_{0\le t \le mT} Q^m_k (t)  m^{-1} \ge \phi \right) \nln
\leq &\pb(Q^m_k(0)\geq 2mq^0_k) + \max_{r= 1, \ldots, 2mq^0_k} \pb \left(\sup_{0\le t \le mT} Q^m_k (t) \ge   m  \phi \bbar Q^m(0) = r \right) \nln
 = &\pb(Q^m_k(0)\geq 2mq^0_k) + \max_{r= 1, \ldots, 2mq^0_k} \pb \left(\sup_{0\le t \le mT} Q^m_k (t) - Q^m_k \ge  m \phi  - r\bbar Q^m(0) = r \right) \nln
  \leq &\pb(Q^m_k(0)\geq 2mq^0_k) + \max_{r= 1, \ldots, 2mq^0_k} \pb \left(\sup_{0\le t \le mT} A^m_k (t)  \ge  m \phi - r\bbar Q^m(0) = r \right) \nln
\sk{a}{=}& \pb(Q^m_k(0)\geq 2mq^0_k) + \pb \left( \sup_{0\leq t\leq mT} A^m_k (t) \ge  (1+\epsilon)m\lambda_k  T  \right) \cr
\sk{b}{\leq}  & \pb(Q^m_k(0)\geq 2mq^0_k) \nln
&  + \pb \left(  \sup_{0 \leq t \leq mT} \lt(  A^m_k (t)  - \int_{0}^{t}\lambda_k \mathbb{I}(C^m(s)=k)ds \rt)  \,  \ge \,  m[(1+\epsilon)\lambda_k T  - \lambda_kT ] \right) \cr
\sk{c}{\leq} & \pb(m^{-1}Q^m_k(0)\geq 2q^0_k) +  2\exp\left(- \lambda_k m T\cdot h(\epsilon ) \right), \label{eq:qmax}
\end{align}
where step $(a)$ follows from the fact that $A^m_k(\cdot)$ is independent of $Q^m_k(0)$, and hence the maximum in the second term is attained by setting $r=2mq^0_k$. Step $(b)$ follows from $\lambda_k \mathbb{I}(C^m(s)=k) \leq \lambda_k$ for all $t\in \rp$, and $(c)$ from the inequality in Eq.~\eqref{eq:r1-arriving}, by replacing $\epsilon$ with $ \epsilon \lambda_kT$. 

We are now ready to establish Eq.~\eqref{eq:X-2}: 
\begin{align}
& \pb\left( \sup_{0 \le t \le T} \left|\frac{1}{m}\left(D^m_k (mt) - \int_{0}^{mt}Q^m_k (s) m^{-1} ds \right)\right| > \epsilon \right)  \nln
=&   \pb\left( \sup_{0 \le t \le mT} \left|D^m_k (t) - \int_{0}^{t}Q^m_k (s)  m^{-1}  ds \right| > m \epsilon \right) \nln
\sk{a}{\leq} &  2\exp\lt(-\phi mT\cdot h\lt(\frac{\epsilon}{\phi T} \rt) \rt)+\pb \left( \sup_{0\le t \le mT}  \psi(W^m(t)) \ge \phi \right) \nln
\sk{b}{\leq} &  2\exp\lt(-\phi mT\cdot h\lt(\frac{\epsilon}{\phi T} \rt) \rt)  \nln
& +  2\exp\left(- \lambda_k m T\cdot h(\epsilon ) \right) + \pb(m^{-1} Q^m_k(0)\geq 2q^0_k), 
\end{align}
where step $(a)$ follows from Lemma \ref{lem:poisson},   and $(b)$ from Eq.~\eqref{eq:qmax}. Note that $\phi, \epsilon$ and $T$ are positive constants, and by our assumption, $\lim_{m \to \infty}\pb(|m^{-1}Q^m_k(0)  - q^0_k|> \delta) = 0$ for all $\delta>0$. Therefore, Eq.~\eqref{eq:X-2} follows by taking the limit in the above inequality as $m\to \infty$. This completes the proof of Proposition \ref{prop:bd1}.

\subsection{Proof of Proposition~\ref{prop:bd2}}
\label{sec:proof-lemma-refl-1}

For the purpose of this proof, we will use an alternative representation of the fluid solutions using integral equations. Define the drift function, $F: \rp^\calK \to \rp^\calK$: 
\begin{equation}
 F_k(q) = \lambda_k p_k(q) - q_k, \quad q\in \calK, 
 \end{equation} 
where $p(\cdot)$ is defined in Eq.~\eqref{eq:pdef}. It can be verified from the definition of $p(\cdot)$ that there exists $l_F \in \rp$ such that $F_k(\cdot)$ is $l_F$-Lipschitz continuous for all $k\in \calK$. Fix $q^0 \in \rp^\calK$, let $\{q(t)\}_{t \in \rp}$, $q(t)\in \rp^\calK$, be a solution to the following integral equation: 
\begin{equation}
q_k(t) = q^0_k + \int_0^t F_k(q(s))ds, \quad k \in \calK. 
\label{eq:dif-2}
\end{equation}
Similar to Lemma \ref{lem:fluidUniqu}, it is not difficult to show that the function $q(\cdot)$ defined in Eq.~\eqref{eq:dif-2} exists, is unique, and coincides with the fluid solution with initial condition $q^0$.  We have that
\begin{align}
&  \left\| V^m(t)-q(t) \right\|  \nln
 =& \sum_{k\in \calK} \left|\int_0^t G_k(V^m(s),\col(s)) - F_k(q(s)) ds \right| \cr
\le&\sum_{k\in \calK}\left|\int_0^t G_k(V^m(s),\col(s)) -  F_k(V^m(s )) ds \right| +\int_0^t \left\| F(V^m(s )) - F (q(s)) \right\| ds \cr
\le&\sum_{k\in \calK}\left|\int_0^t G_k(V^m(s),\col(s)) -  F_k(V^m(s )) ds \right| + l_F \int_0^t \left\| V^m(s)-q(s) \right\| ds,
\label{eq:pf5}
\end{align}
where the last inequalities come from the fact that $F_k(\cdot)$ is $l_F$-Lipschitz continuous for all $k\in \calK$. In a manner analogous to Eq.~\eqref{eq:Qv} from the proof of Proposition \ref{prop:bd1}, by applying Gronwall's lemma (Proposition \ref{prop:gronwall} in Appendix \ref{app:technicaLemma}),  we have that, for all $\epsilon>0$, 
\begin{align}
& \pb\lt( \sup_{0 \le t \le T} \left\| V^m(t)-q(t) \right\| > \epsilon \rt)  \nln
\leq & \pb\lt(   \sup _{0 \leq t\leq T} \sum_{k\in \calK}\left|\int_0^t G_k(V^m(s),\col(s)) -  F_k(V^m(s )) ds \right| > \epsilon e^{-l_F T} \rt)  .
\label{eq:prop2Gronwall}
\end{align} 
Therefore, in order to establish Proposition \ref{prop:bd2}, it suffices to show that
\begin{equation}
\lim_{m\to\infty}\pb\left( \sup_{0\le t \le T} \sum_{k\in \calK} \left|\int_0^t G_k(V^m(s),\col(s)) -  F_k(V^m(s )) ds \right| > \epsilon \right) = 0 ,\quad \forall \epsilon>0. 
\label{eq:mean-lem2}
\end{equation}

In what follows, we will show \eqref{eq:mean-lem2} by using the discrete-time embedded process of $C^m(\cdot)$ and analyzing the system dynamics at the times when new choices are chosen.  We begin by introducing some notation. Fixing $i\in \zp$, we denote by $S_i$ the $i$th update point, i.e., time of the $i$th update in $C^m(\cdot)$, with $S_0 \bydef 0$, and
\begin{equation}
 \tau^m_i = m^{-1}(S_{i}-S_{i-1}), \quad i\in \N. 
 \end{equation} 
 Note that $\{\tau_i\}_{i\in \N}$ are independent exponential random variables with mean $(m\beta_m)^{-1}$. Let $\{\qol [i]\}_{i\in \N}$ be the discrete-time process, where $\qol[i]$ corresponds to the value of $\qol(\cdot)$ immediately following the $i$th update point: 
\begin{equation}
\qol_k[i] = m^{-1} Q^m_k(S_i) = \qol_k(S_i/m), \quad k \in \calK, i \in \N. 
\end{equation}
and $\{Z^m_k[i]\}_{i \in \N}$ be the process of indicator variables: 
\begin{equation}
Z^m_k[i] = \mathbb{I}(C^m(S_i)  = k), \quad k\in \calK, i \in \N. 
\end{equation}
That is, $Z^m_k[i] =1$ if action $k$ is selected on the $i$th update point. Finally, let $\{ \qwh(t) \}_{t\ge 0}$ be a right-continuous piece-wise constant process which coincides with $\qol(\cdot)$ at the points $\{S_i/m\}_{i \in \N}$: 
\begin{equation}
\qwh(t) = \qol[i], \quad \forall t \in [S_i/m , \, S_{i+1}/m). 
\end{equation}

Fix $k\in \calK$, and denote by $I_t$ the number of updates in $C^m(\cdot)$ by time $t$, i.e.,
\begin{equation}
I_t = \max\{i: S_i \leq t \}. 
\end{equation}
 By the triangle inequality, we have that
\begin{align}
& \lambda_k^{-1} \sup_{0\le t \le T} \left|\int_0^t  G_k(V^m(s),\col(s)) -  F_k(V^m(s )) ds \right|  \nln
=& \lambda_k^{-1}  \sup_{0\le t \le T}\left|\int_0^t \lambda_k \left(\mathbb{I}(\col(s)=k) - p_k(V^m(s )) \right) ds \right|   \cr
\le &  \sup_{0\le t \le T} \lt[ \left|\int_0^t  \left(\mathbb{I}(\col(s)=k) - p_k(\qwh (s)) \right) ds \right|  +  \left|\int_0^t  \left(p_k(\qwh (s))-p_k(V^m(s )) \right) ds \right| \rt] \nln
= &  \sup_{0\le t \le T} \lt[ \left|\sum_{i=1}^{I_{mt}+1}\tau^m_i \left(Z^m_k[i] - p_k(\qol [i]) \right) \right|  +  \left|\int_0^t  \left(p_k(\qwh (s))-p_k(V^m(s )) \right) ds \right| \rt] \nln
\leq&  \sup_{0\le t \le T} \left|\sum_{i=1}^{I_{mt}+1}\tau^m_i\left(Z^m_k[i] - p_k(\qol [i]) \right) \right|  \nln
& \quad +  \sup_{0\le t \le T} \left|\int_0^t  \left(p_k(\qwh (s))-p_k(V^m(s )) \right) ds \right|.
\label{eq:pf-GF}
\end{align}

We now derive upper bounds on tail probabilities for each term on the right-hand side of Eq.~\eqref{eq:pf-GF}. For the first term, define 
\begin{equation}
\xi^m_{n} = \sum_{i=1}^n\tau^m_i \left(Z^m_k[i] - p_k(\qol [i]) \right), \quad n \in \N. 
\end{equation}
Fix $m\in \N$ and $\epsilon>0$. Let
\begin{equation}
 J^m_\epsilon =  (1+\epsilon)m\beta_m T.
 \end{equation} 
(To avoid the excessive use of floors and ceilings, we assume that $J^m_\epsilon$ is a positive integer. The results extend easily to the general case.) Define the events
\begin{align}
\calA^m_\epsilon =& \lt\{ I_{mT}  < J^m_\epsilon\rt\}, \quad  \mbox{and} \quad  
\calB^m_\epsilon = \lt\{\max_{0\le n \le J^m_\epsilon } |\xi^m_n| \leq \epsilon  \rt\}. 
\end{align}
We have that
\begin{align}
& \pb\lt( \sup_{0\le t \le T} \left|\sum_{i=1}^{I_{mt}+1}\tau^m_i\left(Z^m_k[i] - p_k(\qol [i]) \right) \right| \geq  \epsilon \rt)  \nln
= & \pb\lt(\max_{0\le n \le I_{mT}} |\xi^m_{n}| \geq \epsilon\rt) \nln
\leq& 1-\pb(\calA^m_\epsilon \cap \calB^m_\epsilon) \nln
\leq& (1-\pb(\calA^m_\epsilon))+(1-\pb(\calB^m_\epsilon)). 
\label{eq:boundtoPAPB}
\end{align}
With the above equation in mind, we now proceed to demonstrate that both $)1-\pb(\calA^m_\epsilon))$ and ($1-\pb(\calB^m_\epsilon))$ converge to zero in the limit as $m\to \infty$. 
Note that $I_{mT}$ is a Poisson random variable with mean $m\beta_m T$.  Using elementary tail bounds on the Poisson distribution, we have that
\begin{equation}
1- \pb(\calA^m_\epsilon) = \pb(I_{mT} \geq (1+\epsilon)m\beta_m T)  \leq \frac{1+\epsilon}{\epsilon^2}(m\beta_m T)^{-1}, 
\label{eq:pbcalAbound}
\end{equation}
which converges to zero as $m\to \infty$. 

We next turn to the value of $\pb(\calB^m_\epsilon)$.  Recall from the definition of $Z^m_k[i]$ that
\begin{equation}
\E(Z^m_k[i] \bbar  \qol[i]) = p_k(\qol[i]). 
\end{equation}
It is therefore not difficult to verify that $\{\xi^m_n\}_{n\in \N}$ is a martingale, and our objective would be to derive an upper bound on its maximum upward excursion over $\{0,\ldots, J^m_\epsilon\}$. Unfortunately, we cannot apply the Azuma-Hoeffding  inequality, because the $i$th increment of $\{\xi^m_n\}_{n\in \N}$ involves the term $\tau^m_i$, which does not admit a bounded support. Instead, we will use the following upper bound on the moment generating function of $\xi^m_{n}$, whose proof is based on Doob's inequality and is given in Appendix \ref{app:lem:xiinduct}.  
\begin{lemma} 
\label{lem:xiinduct}
Fix $n \in \N$ and $\theta \in (0, m\beta_m)$. We have that
\begin{equation}
\E\left(\exp\left( \theta \xi^m_n\right) \right)  \le \exp \left(   \frac{ \theta^2 n }{ 4m\beta_m  \left( m\beta_m -\theta \right)  } \right). 
\end{equation}
\end{lemma}

We are now ready to establish an upper bound on the quantity $1-\pb(\calB^m_\epsilon)$. Recall that $J^m_\epsilon = (1+\epsilon)m\beta_mT$. Fix $\eta \in \lt(0, \min\{T,1\}/2 \rt)$, and let $\theta_0 = \eta\frac{\epsilon}{1+\epsilon}T^{-1} m\beta_m.$ In particular, $\theta_0\in (0, m\beta_m/2)$. We have that
\begin{align}
1-\pb(\calB^m_\epsilon) = & \mathbb{P}\left(\max_{0\le n \le J^m_\epsilon } \xi^m_n \ge \epsilon \right)\cr
\stackrel{(a)}{\le}  & \inf_{\theta >0 } \exp \left( -\theta \epsilon \right) \E\left(\exp\left( \theta \xi^m_{J^m_\epsilon}\right) \right) \cr
\stackrel{(b)}{\leq} &\exp \left( -\theta_0 \epsilon +   \frac{ J^m_\epsilon \theta_0^2/4 }{ m\beta_m  \left( m\beta_m -\theta_0 \right)  } \right) \cr
\sk{c}{\leq} & \exp \left( -\theta_0 \epsilon +   \frac{ J^m_\epsilon \theta_0^2 }{2 (m\beta_m)^2   } \right) \nln
= & \exp \left( -\theta_0 \lt( \epsilon -   \frac{ (1+\epsilon)T \theta_0 }{2 m\beta_m  }  \rt) \right) \nln
\sk{d}{\leq}& \exp \left( -\theta_0 \lt( \epsilon - \epsilon/2  \rt) \right) \nln
= & \exp \left( - \eta \frac{\epsilon^2}{2+2\epsilon }T^{-1}m\beta_m  \right). 
\label{eq:doob1}
\end{align}
Step $(a)$ follows from Doob's inequality and the fact that,  for any $\theta>0$, the sequence $\left\{ \exp\left(\theta \xi^m_n \right)\right\}_{n\in \N}$ is a submartingale, and $(b)$  from Lemma~\ref{lem:xiinduct} with $n = J^m_\epsilon$. Step $(c)$ is due to $\theta_0<m\beta_m/2$, and hence $m\beta_m -\theta_0>m\beta_m/2$. Finally, step $(d)$ follows from the definition of $\theta_0$ and that $\eta < 1$. 

With a derivation analogous to that of Eq.~\eqref{eq:doob1}, we have that 
\begin{equation}
\mathbb{P}\left(\min_{0\le n \le J^m_\epsilon } \xi^m_n   \leq - \epsilon \right) \le \exp \left( - \eta \frac{\epsilon^2}{2+2\epsilon }T^{-1} m\beta_m  \right) .\label{eq:doob2}
\end{equation}
Therefore, combining Eq.~\eqref{eq:boundtoPAPB} with Eqs.~\eqref{eq:pbcalAbound}, \eqref{eq:doob1} and \eqref{eq:doob2}, we have that
\begin{align}
& \pb\lt( \sup_{0\le t \le T} \left|\sum_{i=1}^{I_{mt}+1}\tau^m_i\left(Z^m_k[i] - p_k(\qol [i]) \right) \right| \geq  \epsilon \rt) \nln
\leq& (1-\pb(\calA^m_\epsilon))+(1-\pb(\calB^m_\epsilon)) \nln
\leq &   \frac{1+\epsilon}{T\epsilon^2}(m\beta_m )^{-1} + 2\exp \left( - \eta \frac{\epsilon^2}{2+2\epsilon } T^{-1} m\beta_m  \right). 
\end{align}
Since $m\beta_m \to \infty$ as $m\to \infty$, the above equation further implies that 
\begin{align}
 &  \lim_{m \to \infty} \pb\lt( \sup_{0\le t \le T} \left|\sum_{i=1}^{I_{mt}+1}\tau^m_i\left(Z^m_k[i] - p_k(\qol [i]) \right) \right| \geq  \epsilon \rt) \nln
  \leq & \lim_{m \to \infty}  \lt( \frac{1+\epsilon}{T\epsilon^2}(m\beta_m )^{-1} + 2\exp \left( - \eta \frac{\epsilon^2}{2+2\epsilon } T^{-1} m\beta_m  \right) \rt)   \nln
=& 0
\label{eq:pdfGFterm1}
\end{align}

We now bound the second term in Eq.~\eqref{eq:pf-GF}.  It is not difficult to verify, from Eq.~\eqref{eq:pdef}, that  there exists $l_p\in \rp$ such that $p_k(\cdot)$ is $l_p$-Lipschitz continuous for all $k\in \calK$.  We have that
\begin{align}
& \pb\left(\sup_{0\le t \le T}  \left|\int_0^t  \left(p_k(\qwh (s))-p_k( V^m(s )) \right) ds   \right| \geq \epsilon \right) \nln
 \leq &   \pb\left(  \sup_{0\le t \le T} \int_0^t l_p \left|\widehat{Q}_k^m (s) - V^m_k(s) \right| ds \geq \epsilon \right) \nln
\leq & \pb\left(  \int_0^T l_p \left|\widehat{Q}_k^m (s) - V^m_k(s) \right| ds \geq \epsilon \right). 
 \label{eq:doob4}
\end{align}
It therefore suffices to show that, if $m\beta_m \to \infty$ as $m \to \infty$, then
\begin{equation}
 \lim_{m\to \infty}\pb\left( \int_0^T \left|\widehat{Q}_k^m (s) - V_k^m(s) \right| ds  \geq \epsilon \right) = 0, \quad \forall \epsilon > 0. 
 \label{eq:EqtoVto01}
 \end{equation} 
To this end, we have that
\begin{align}
& \int_0^T \left|\widehat{Q}_k^m (s) - V_k^m(s) \right| ds \nln
= & \sum_{i=0}^{I_{mT}} \int_{S_{i}/m}^{S_{i+1}/m}\lt | \qwh_k(t) - V^m_k(s) \rt| ds \nln
\leq &  \sum_{i=0}^{I_{mT}} \int_{S_{i}/m}^{S_{i+1}/m}\lt | \qwh_k(t) - V^m_k(S_i/m) \rt| + \big| V^m_k(s)-V^m_k(S_{i}/m) \big | \, ds \nln
\sk{a}{=}& \sum_{i=0}^{I_{mT}} \tau^m_{i+1} \lt| \qol_k[i] - V^m_k(S_i/m)\rt| + \sum_{i=0}^{I_{mT}} \int_{S_{i}/m}^{S_{i+1}/m}|V^m_k(s)-V^m_k(S_{i}/m)|ds \nln
\leq&  T \sup_{0\leq s\leq T} \lt| \qol_k(s) - V^m_k(s) \rt| + \sum_{i=0}^{I_{mT}}  \tau^m_{i+1}\lt(\tau^m_{i+1} \sup_{0\leq s\leq T} \lt|G_k(V^m_k(s), \col(s)) \rt| \rt)  \nln
= & T \sup_{0\leq s\leq T} \lt| \qol_k(s) - V^m_k(s) \rt| + \lt(\sup_{0\leq s\leq T} \lt|G_k(V^m_k(s), \col(s)) \rt|\rt) \sum_{i=1}^{I_{mT}+1}  (\tau^m_i)^2\nln
\sk{b}{\leq}&  T \sup_{0\leq s\leq T} \lt| \qol_k(s) - V^m_k(s) \rt| + \lt(\lambda_k+\sup_{0\leq s\leq T} |V^m_k(s)| \rt) \sum_{i=1}^{I_{mT}+1}  (\tau^m_i)^2 \nln
\sk{c}{\leq} & T \sup_{0\leq s\leq T} \lt| \qol_k(s) - V^m_k(s) \rt| + [ q_0+\lambda_k(T+1) ] \sum_{i=1}^{I_{mT}+1} (\tau^m_i)^2. 
\label{eq:qwhtoVk}
\end{align}
In step $(a)$ we have invoked the property that $\qwh_k(\cdot)$ is piece-wise constant. Step $(b)$ follows from the definition of $G_k(\cdot, \cdot)$ in Eq.~\eqref{eq:gkdef}, and $(c)$ from the fact that $|V^m_k(t)| \leq q^0+\lambda_kt$ for all $t\in \rp$ (Eq.~\eqref{eq:vmkbound}). It remains to derive an upper bound on the tail probabilities of the term $\sum_{i=1}^{I_{mT}+1} (\tau^m_i)^2$, which is isolated in the form of the following lemma. The proof involves an elementary application of Markov's inequality, and is given in Appendix \ref{app:lem:xisqrTailBound}. 

\begin{lemma}
\label{lem:xisqrTailBound}
Suppose that $m \beta_m \to \infty $ as $m\to \infty$. We have that
\begin{equation}
\lim_{m \to \infty} \pb\lt( \sum_{i=1}^{I_{mT}+1} (\tau^m_i)^2 \geq \epsilon \rt) =0.  \quad \forall \epsilon>0. 
\end{equation}
\end{lemma}

We are now ready to prove Eq.~\eqref{eq:EqtoVto01}. By Eq.~\eqref{eq:qwhtoVk}, we have that
\begin{align}
& \lim_{m \to \infty}\pb\lt(\int_0^T \left|\widehat{Q}_k^m (s) - v_k^m(s) \right| ds \geq  \epsilon  \rt)  \nln
\leq & \lim_{m\to \infty} \pb\lt( \sup_{0\leq s\leq T} \lt| \qol_k(s) - V^m_k(s) \rt|  \geq  \frac{\epsilon}{2T} \rt)  \nln
& + \lim_{m\to \infty}\pb\lt(\sum_{i=1}^{I_{mT}+1} (\tau^m_i)^2  \geq  \frac{\epsilon}{2[ q_0+\lambda_k(T+1) ]} \rt) \nln
{=} & 0. 
\label{eq:limQkminusvk}
\end{align}
for all $\epsilon>0$, where the first inequality follows from a union bound, and the last step from Proposition \ref{prop:bd1} and Lemma \ref{lem:xisqrTailBound}. Substituting Eqs.~\eqref{eq:pdfGFterm1} and \eqref{eq:limQkminusvk} into Eq.~\eqref{eq:pf-GF}, we have that
\begin{align}
& \lim_{m\to\infty}\pb\left( \sup_{0\le t \le T} \sum_{k\in \calK} \left|\int_0^t G_k(V^m(s),\col(s)) -  F_k(V^m(s )) ds \right| > \epsilon\right) \nln
= &   \lim_{m \to \infty} \pb\lt( \sup_{0\le t \le T} \left|\sum_{i=1}^{I_{mt}+1}\tau^m_i\left(Z^m_k[i] - p_k(\qol [i]) \right) \right| \geq \epsilon/\lambda_k \rt)   \nln
&  \qquad  +  \lim_{m \to \infty} \pb \lt( \sup_{0\le t \le T} \left|\int_0^t  \left(p_k(\qwh (s))-p_k(V^m(s )) \right) ds \right| \geq \epsilon/\lambda_k  \rt)   = 0. 
\end{align}
This establishes Eq.~\eqref{eq:mean-lem2}, which, in light of Eq.~\eqref{eq:prop2Gronwall} completes the proof of Proposition \ref{prop:bd2}.

\subsection{Proof of Proposition \ref{prop:EQconditional}} 
\label{app:prop:EQconditional}

\bpf
We begin by showing the following, strengthened version of Lemma \ref{lem:stochDom1}, which states that a similar stochastic dominance property holds for $Q^m(\infty)$ even when conditioning on $C^m(\cdot)$ being of a specific value. 

 \begin{lemma}
 \label{lem:conditionalStochDom}
Let $W^m = (Q^m, C^m)$ be a random vector drawn from the steady-state distribution, $W^m(\infty)$. There exist constants ,  $m_0$ and $\gamma >0$, such that for all $m\geq m_0$, 
\begin{equation}
Q^m_i \bbar\{C^m = k \} \preceq \gamma  U^m_{\lambda_1}, \quad  \forall i, k\in \calK. 
\end{equation}
 \end{lemma}

\bpf  Let $\calQ^m = \{{x}/{m}: \,  x\in \zp^\calK\}$. We have that, for all $k\in \calK$, $n\in \N$, 
\begin{align}
\pb(C^m[n] = k ) {=} & \sum_{u\in \calQ^m} \frac{u_k \vee \alpha_0 }{\sum_{i \in \calK} ( u_i \vee \alpha_0 ) }  \pb \lt(\qol[n] = u \rt) \nln
\geq  &  \sum_{u\in \calQ^m} \frac{ \alpha_0 }{K\alpha_0+ \sum_{i \in \calK} u_i }  \pb \lt(\qol[n] = u \rt) 
\end{align}
By Eq.~\eqref{eq:Cm-discrete} of Lemma \ref{lem:pasta}, the above equation implies that
\begin{equation}
\pb(C^m = k ) = \lim_{n \to \infty} \pb(C^m[n] = k )\geq  \sum_{u\in \calQ^m} \frac{ \alpha_0 }{K\alpha_0+ \sum_{i \in \calK} u_i }  \pb \lt(\qol[\infty] = u \rt). 
\label{eq:cmklim}
\end{equation}
where the last step follows from the fact that $\frac{\alpha_0}{K\alpha_0+\sum_{i \in \calK} u_i}$ is always bounded from above by $1/K$.  

By Eq.~\eqref{eq:Umlamconverg}, we have that $U^m_{\lambda_1}/m \to \lambda_1$ almost surely as $m\to \infty$. Therefore, there exist $m_0$ and $y>0$, such that
\begin{equation}
\pb\lt( \frac{1}{m}U^m_{\lambda_1} \geq \lambda_1+y \rt) \leq \frac{1}{2K^2}, \quad \forall m\geq m_0. 
\label{eq:ytoK}
\end{equation}

Combining Eq.~\eqref{eq:Qmdom-discrete} in Lemma \ref{lem:pasta} with Eqs.~\eqref{eq:cmklim} and \eqref{eq:ytoK}, we have that, for all $m\geq m^0$, 
\begin{align}
& \pb(C^m = k ) \nln
\sk{a}{\geq} &   \sum_{u\in \calQ^m} \frac{ \alpha_0 }{K\alpha_0+ \sum_{i \in \calK} u_i }  \pb \lt(\qol[\infty] = u \rt) \nln
{\geq} &  \pb\lt( \max_{i\in \calK}\qol_i[\infty] \leq \lambda_1+y \rt)\frac{\alpha_0}{K\alpha_0+K(\lambda_1+y)} + \pb\lt( \max_{i\in \calK}\qol_i[\infty]  > \lambda_1+y \rt)\cdot 0 \nln
\sk{b}{\geq} &   \lt( 1-K \pb\lt(\frac{1}{m}U^m_{\lambda_1}  \leq \lambda_1+y \rt) \rt)\frac{\alpha_0}{K\alpha_0+K(\lambda_1+y)} \nln
\sk{c}{\geq} & \frac{\alpha_0\lt(1-{1}/{2K}\rt)}{K\alpha_0+K(\lambda_1+y)} \nln
= & \gamma^{-1}, 
\label{eq:cmklim2}
\end{align}
where $\gamma \bydef\frac{K\alpha_0+K(\lambda_1+y)}{\alpha_0\lt(1-{1}/{2K}\rt)}$. Step $(a)$ follows from Eq.~\eqref{eq:cmklim}, $(b)$ from Lemma \ref{lem:stochDom1} and a union bound, and $(c)$ from Eq.~\eqref{eq:ytoK}. 

Fix $x \in \zp$. We have that, for all $m\geq m_0$, 
\begin{align}
\pb\lt(Q^m_i  \geq x \bbar C^m = k \rt) = & \frac{\pb\lt(Q^m_i  \geq x ,  \, C^m = k \rt)}{\pb( C^m = k )}  
\leq  \frac{\pb\lt(Q^m_i  \geq x  \rt)}{\pb( C^m = k )}\nln
\sk{a}{\leq} &  \frac{\pb\lt(U^m_{\lambda_1}  \geq x  \rt)}{\pb( C^m = k )} \nln
\sk{b}{\leq} & \gamma \pb\lt(U^m_{\lambda_1}  \geq x \rt), 
\end{align}
for all $i,k \in \calK$, where step $(a)$ follows from Lemma \ref{lem:stochDom1}, and $(b)$ from Eq.~\eqref{eq:cmklim2}. Since the above inquality holds for all $x\in \zp$, this completes the proof of Lemma \ref{lem:conditionalStochDom}.  \qed

We now prove the convergence in  Eq.~\eqref{eq:convgProbQmh}.  Recall that the first update point, $S^m_1$, is exponentially distributed with mean $1/\beta_m$, and independent from $W^m(0)$. Define the event $\calE^m = \{S^m_1 \geq h(m)\}$.  We have that
\begin{equation}
\pb(\calE^m) = \pb(S^m_1 \geq h(m)) = \exp(-\beta_mh(m)) \to 1, \quad \mbox{ as $m \to \infty$.} 
\label{eq:Eepslimit}
\end{equation}
 Fix $i, k \in \calK$. Recall from Eq.~\eqref{eq:Umlamconverg} that $U^m_{\lambda_1}/m$ converges to $\lambda_1$ almost surely as $m\to \infty$. This implies that there exists $\upsilon>0$, independent of $m$, such that 
 \begin{equation}
 \limsup_{m \to \infty}\pb\lt( \qol_i[0] > \upsilon \bbar C^m[0] = k\rt) \sk{a}{\leq} \limsup_{m\to \infty}\pb\lt( U^m_{\lambda_1}/m > \gamma^{-1}\upsilon \rt)  =0,
 \label{eq:Qmbigv}
 \end{equation}
where step $(a)$ follows from Lemma \ref{lem:conditionalStochDom}.

Fix $x\in \rp$. We have that
\begin{align}
&  \limsup_{m \to \infty}\pb \lt( \qol_i [1] \geq x \bbar C^m[0]=k\rt) \nln
 \leq  & \limsup_{m \to \infty} \lt[\pb\lt(\qol_i[1] \geq x \bbar C^m[0]=k,  \, \calE^m \rt) + (1-\pb(\calE^m\bbar C^m[0]=k) )  \rt] \nln
\sk{a}{=} & \limsup_{m \to \infty} \lt[\pb\lt(\qol_i[1] \geq x \bbar C^m[0]=k,  \, \calE^m \rt)+ (1-\pb(\calE^m ) )  \rt] \nln
\sk{b}{=} & \limsup_{m \to \infty} \pb\lt(\qol_i[1] \geq x \bbar C^m[0]=k,  \, \calE^m \rt)
\label{eq:Qhmupper1}
\end{align}
where step $(a)$ follows from the independence between $\calE^m$ and $C^m[0]$, and $(b)$ from Eq.~\eqref{eq:Eepslimit}. 

We now bound the term on the right-hand side of Eq.~\eqref{eq:Qhmupper1}. Denote by $B(n,p)$ a binomial random variable with $n$ trials and a success probability of $p$ per trial, and by $U^m_{\lambda}(t)$ the number of jobs in system at time $t$ in an initially empty $M/M/\infty$ queue with arrival rate $\lambda$ and departure rate $1/m$. We have that
\begin{align}
& \pb\lt(\qol_i[1] \geq x \bbar C^m[0]=k,  \, \calE^m \rt) \nln
 \sk{a}{=} & \pb\lt(\frac{1}{m} \lt( U^m_{\lambda_k}(S^m_1)\mathbb{I}(i=k) + B\lt(Q^m_i(0), e^{-S^m_1/m} \rt) \rt) \geq x \bbar C^m[0]=k, \calE^m \rt) \nln
 \sk{b}{\leq} & \pb\lt(\frac{1}{m} \lt( U^m_{\lambda_k}(S^m_1)\mathbb{I}(i=k) + B\lt(Q^m_i(0), e^{-h(m)/m} \rt) \rt) \geq x \bbar C^m[0]=k \rt) \nln
 \sk{c}{\leq} & \pb\lt(\frac{1}{m} \lt( U^m_{\lambda_k}\mathbb{I}(i=k) + B\lt(Q^m_i(0), e^{-h(m)/m} \rt) \rt) \geq x \bbar C^m[0]=k \rt) \nln
{\leq} & \pb\lt( \frac{1}{m} \lt( U^m_{\lambda_k}\mathbb{I}(i=k) + B\lt( \upsilon m , e^{-h(m)/m} \rt) \rt) \geq x\bbar  C^m[0]=k, \qol_i(0) \leq \upsilon \rt) \nln
& +  \pb \lt( \qol_i(0) > \upsilon \bbar C^m[0]=k \rt) \nln
= & \pb\lt( \frac{1}{m} \lt( U^m_{\lambda_k}\mathbb{I}(i=k) + B\lt( \upsilon m , e^{-h(m)/m} \rt) \rt) \geq x \rt) \nln
&  + \pb \lt( \qol_i(0) > \upsilon\bbar C^m[0]=k \rt). 
 \label{eq:Qhmupper2}
\end{align}
For step $(a)$, note that each unit of reward initially present at time $t=0$ has probability of $\exp(-S^m_1/m)$ or remaining in the system by $t=S^m_1$. Therefore, the rewards in site $i$ at time $S^m_1$ satisfy the following decomposition: 
\begin{equation}
Q^m_i[1] \stackrel{d}{=} U^m_{\lambda_k}(S^m_1)\mathbb{I}(i=k)  + B\lt(Q^m_i(0), e^{-S^m_1/m} \rt).
\label{eq:Qdecomp} 
\end{equation}
The first term corresponds to the units of rewards at $t=S^m_1$ that had arrived during the interval $(0,S^m_1)$, and hence is non-zero only if $i=k$. The second term corresponds to those individuals initially present at $t=0$ who remained in the system by  $t=S^m_1$. Step $(b)$ follows from the definition of $\calE^m$, and  $(c)$   from the well-known fact that the number of jobs in system in an initially empty $M/M/\infty$ queue at any time is always stochastically dominated by its steady-state distribution.  

Because $h(m)/m\to \infty$ as $m\to \infty$, we have that $\lim_{m\to \infty} \E\lt(\frac{1}{m}B\lt( \upsilon m , e^{-h(m)/m} \rt)\rt) =0$. Applying Markov's inequality, we obtain that
\begin{equation}
 \frac{1}{m}B\lt( \upsilon m , e^{-h(m)/m} \rt) \stackrel{P}{\rightarrow} 0, \quad \mbox{as $m\to \infty$},
 \label{eq:binomialBound} 
\end{equation}
where $\stackrel{P}{\rightarrow}$ denotes convergence in probability.  Recall from Eq.~\eqref{eq:Umlamconverg} that, almost surely, 
\begin{equation}
\frac{1}{m}U^m_{\lambda_k}\mathbb{I}(i=k) \to \lambda_k\mathbb{I}(i=k) = q^*_{k,i}. 
\label{eq:UIndicate}
\end{equation}
Fix $\epsilon>0$, and substitute Eq.~\eqref{eq:Qhmupper2} into Eq.~\eqref{eq:Qhmupper1}. We have that
\begin{align}
& \limsup_{m \to \infty}\pb \lt( \qol_i [1] \geq q^*_{k,i}+\epsilon \bbar C^m[0]=k\rt)  \nln
{=}&  \limsup_{m \to \infty} \pb\lt(\qol_i[1] \geq q^*_{k,i}+\epsilon \bbar C^m[0]=k,  \, \calE^m \rt) \nln
\leq & \limsup_{m \to \infty} \pb\lt( \frac{1}{m} \lt( U^m_{\lambda_k}\mathbb{I}(i=k) + B\lt( \upsilon m , e^{-h(m)/m} \rt) \rt) \geq q^*_{k,i}+\epsilon \rt)  \nln 
& \quad + \limsup_{m \to \infty} \pb \lt( \qol_i(0) > \upsilon\bbar C^m[0]=k \rt) \nln
\sk{a}{=} & \limsup_{m \to \infty} \pb\lt( \frac{1}{m} \lt( U^m_{\lambda_k}\mathbb{I}(i=k) + B\lt( \upsilon m , e^{-h(m)/m} \rt) \rt) \geq q^*_{k,i}+\epsilon \rt) \nln
\sk{b}{=} & 0, 
\label{eq:Qhmupper3}
\end{align}
where step $(a)$ follows from Eq.~\eqref{eq:Qmbigv}, and $(b)$ from Eqs.~\eqref{eq:binomialBound} and \eqref{eq:UIndicate}. Using the same line of arguments as that in Eqs.~\eqref{eq:Qhmupper1} through \eqref{eq:Qhmupper3}, we can show that
\begin{equation}
 \limsup_{m \to \infty}\pb \lt( \qol_i [1] \leq q^*_{k,i}- \epsilon \bbar C^m[0]=k\rt)=0,
\end{equation}
which, along with Eq.~\eqref{eq:Qhmupper3}, yields that
\begin{equation}
 \limsup_{m \to \infty}\pb \lt( \lt| \qol_i [1] -  q^*_{k,i}\rt| > \epsilon \bbar C^m[0]=k\rt) = 0.
 \label{eq:probcon}
\end{equation}
Since  the above equation holds for all $i,k\in \calK$, this proves Eq.~\eqref{eq:convgProbQmh} in Proposition \ref{prop:EQconditional}. 

We now turn to Eq.~\eqref{eq:convgExpQmh}. Fix $i,k \in \calK$. Using essentially identical arguments as those for Eq.~\eqref{eq:probcon}, we can show that
\begin{equation}
 \limsup_{m \to \infty}\pb \lt( \lt| m^{-1} Q^m_i (h(m)) -  q^*_{k,i}\rt| > \epsilon \bbar C^m[0]=k\rt) = 0, \quad \forall \epsilon>0. 
\end{equation}
We have that, 
\begin{align}
Q^m_i(h(m)) | \{C^m[0]=k\} \sk{a}{\preceq}  & U^m_{\lambda_1}(h(m))+ Q^m_i(0) | \{C^m[0]=k\} \nln
\sk{b}{\preceq} &  U^m_{\lambda_1}(h(m)) + \gamma U^m_{\lambda_1} \nln
 \preceq & (\gamma+1) U^m_{\lambda_1}, 
\end{align}
Step $(a)$ follows from a decomposition similar to Eq.~\eqref{eq:Qdecomp}, by the writing the recallable rewards at time $h(m)$ as those who arrived after $t=0$, which is dominated by $U^m_{\lambda_1}(h(m))$, and those who were in the system at $t=0$, which is dominated by $Q^m_i(0) | \{C^m[0]=k\}$. Step $(b)$ follows from Lemma \ref{lem:stochDom1}. Since $U^m_{\lambda_1}$ is a Poisson distribution with mean $m\lambda_1$, it is not difficult show that there exists a random variable $Y \in \rp$, such that
\begin{equation}
m^{-1}Q^m_i(h(m)) | \{C^m[0]=k\} \sk{a}{\preceq} \frac{1}{m}U^m_{\lambda_1} (\gamma+1) \preceq Y, \quad \forall m \in \N. 
\label{eq:Qm1hmdom}
\end{equation}
Combining Eqs.~\eqref{eq:convgProbQmh} and \eqref{eq:Qm1hmdom}, the dominated convergence theorem implies that, for all $i,k \in \calK$, 
\begin{equation}
\lim_{m \to \infty} \E\lt( \lt| m^{-1}Q^m_i(h(m))- q^*_{k,i} \rt| \bbar C^m[0] = k\rt) = 0. 
\end{equation}
This shows Eq.~\eqref{eq:convgExpQmh}, and thus completes the proof of Proposition \ref{prop:EQconditional}. \qed

\subsection{Proof of Proposition \ref{prop:etaless1}}
\label{app:prop:etaless1}
\bpf 
Fix $\eta < 1$. Under the polynomial reward-matching model, the fluid solution satisfies
\begin{equation}
\dot{q}_k(t) = \lambda_k  \frac{(q_k \vee \alpha_0)^\eta}{\sum_{i\in \calK}(q_i\vee \alpha_0)^\eta} - q_k(t), \quad \forall k \in \calK,
\end{equation}
setting the left-hand side to $0$, we have that a state $q$ is an invariant state of the fluid solutions if
\begin{equation}
\frac{q_k}{(q_k \vee \alpha_0)^\eta} =   \lambda_k\frac{1}{\sum_{i\in \calK}(q_i\vee \alpha_0)^\eta} , \quad \forall k \in \calK. 
\label{eq:alphaqratio}
\end{equation}

We first show that the above equations admit a unique solution, $q^I$. That is, the fluid solutions admit a unique invariant state. Suppose, for the sake of contradiction, that there exist two distinct invariant states $q^I$ and $\tilde{q}^I$. Let 
\begin{equation}
Z= \sum_{i\in \calK}(q^I_i\vee \alpha_0)^\eta
\end{equation} 
and $\tilde{Z}= \sum_{i\in \calK}(\tilde{q}^I_i\vee \alpha_0)^\eta$ denote the denominators on the right-hand side of Eq.~\eqref{eq:alphaqratio} under $q^I$ and $\tilde{q}^I$, respectively. From Eq.~\eqref{eq:alphaqratio}, by considering  separately two cases depending on whether $q^I_k$ is smaller than $\alpha_0$, we have that the invariant state satisfies
\begin{equation}
q_k^I = \begin{cases} 
  \left(\frac{\lambda_k}{Z} \right)^{\frac{1}{1-\eta}} \quad  & (\geq \alpha_0), \quad \mbox{if} \quad \lambda_k \ge Z\alpha^{1-\eta}_0,\\
\lambda_k \frac{ \alpha_0^{\eta}}{Z} \quad & (< \alpha_0), \quad \mbox{if} \quad \lambda_k < Z\alpha^{1-\eta}_0,
\end{cases} \label{eq:alphaqratio2}
\end{equation}
which indicates that if $Z = \tilde{Z}$, then $q^I = \tilde{q}^I$.  Therefore, in order for  $q^I$ and $\tilde{q}^I$ to be distinct, we must have that $\tilde{Z} \neq Z$. Without loss of generality, let us assume that $\tilde{Z} > Z$. Because $\eta<1$, by Eq.~\eqref{eq:alphaqratio2}, we have that $q_k^I$ is a monotonically decreasing function of $Z$, for all $k$. We thus have that $\tilde{q}^I_k < q^I_k$ for all $k \in \calK$. This leads to a contradiction, since when $\tilde{q}^I_k < q^I_k$ for all $k \in \calK$, we will necessarily have that $\tilde{Z} $ is strictly less than $Z$. This proves that the solution to Eq.~\eqref{eq:alphaqratio} must be unique. 

We now find the unique invariant state $q^I$, and for now we assume that such $q^I$ exists. Note that when $\eta<1$, $q_k^I$ is a monotonically increasing function of $\lambda_k$ for $\lambda_k \ge 0$. Eq.~\eqref{eq:alphaqratio2} implies that $q_i \geq q_j$ if and only if $ \lambda_i \geq \lambda_j$, which further implies that $ q^I_1 \ge \dots \ge q^I_K$. Since $q^I_k$ is non-increasing in $k$, we may define $i^*$ as the unique index such that 
\begin{equation}
q^I_{i^*} \ge \alpha_0, \quad \mbox{and} \quad q^I_{i^*+1} < \alpha_0,
\end{equation}
where we define $i^*=0$ if $q^I_1 < \alpha_0$, and $i^* = K$ if $q^I_K \geq \alpha_0$. 

We now consider different values of $\alpha_0$. Suppose that $\alpha_0 \geq \lambda_1/K$. It is not difficult to verify that in this case $i^*=0$, $q^I_k < \alpha_0$ for all $k\in \calK$, and $Z = \sum_{i\in \calK}(q^I_i\vee \alpha_0)^\eta = K\alpha_0^\eta$. It follows from Eq.~\eqref{eq:alphaqratio2}  that we must have $\lambda_1 < Z\alpha^{1-\eta} = K\alpha_0$, or equivalently, $\alpha_0 > \lambda_1/K$, and that 
\begin{equation}
p^I_k= {\lambda_k}/{K}, \quad \forall k \in \calK.
\end{equation}
This proves Item 1 in the proposition. 

Consider next the other extreme where $\alpha_0$ is so small that $i^*=K$ and $q^I_k \geq \alpha_0$ for all $k$. By Eq.~\eqref{eq:alphaqratio2}, this is to say that
\begin{equation}
\lambda_K \geq Z\alpha_0^{1-\eta}.
\label{eq:alpha0cond1}
\end{equation}
In this case, we have that $Z = \sum_{i \in \calK} \lt(\frac{\lambda_{i}}{Z}\rt)^\frac{\eta}{1-\eta},$ which leads to, after rearrangement, 
\begin{equation}
Z= \lt(\sum_{i \in \calK} \lambda_i^{\frac{\eta}{1-\eta}}\rt)^{1-\eta}.
\label{eq:Zexp1}
\end{equation}
Substituting the value of $Z$ from Eq.~\eqref{eq:Zexp1} into \eqref{eq:alpha0cond1}, we obtain the condition on $\alpha_0$: 
\begin{equation}
\alpha_0\leq \lambda_K \frac{\lambda_K^\frac{\eta}{1-\eta}}{\sum_{i \in \calK}
\lambda_i^\frac{\eta}{1-\eta}}.
\end{equation}
The expression of $q^I_k$ in Eq.~\eqref{eq:qik-smeta} is obtained by substituting Eq.~\eqref{eq:Zexp1} into the top line of Eq.~\eqref{eq:alphaqratio2}. This proves Item 2 of the proposition. 

Finally, fix $\alpha_0$ such that $\lambda_K \frac{\lambda_K^\frac{\eta}{1-\eta}}{\sum_{i \in \calK}
\lambda_i^\frac{\eta}{1-\eta}}  <\alpha_0 \leq \lambda_1/K$. In this case, we have that $1\leq i^* \leq K-1$. 
By Eq.~\eqref{eq:alphaqratio2}, we have that for all $k\geq i^*+1$, 
\begin{align}
q^I_k = &\lambda_k \frac{\alpha_0^\eta}{Z} \nln
= & \frac{\lambda_{i^*}}{Z} \cdot \frac{\lambda_k}{\lambda_{i^*}}\alpha^\eta \nln
= & \lt(q^I_{i^*} \rt)^{1-\eta}\frac{\lambda_k}{\lambda_{i^*}}\alpha_0^\eta, 
\end{align}
where the last equality follows from the fact that $q^I_{i^*}\geq \alpha_0$, and hence $q^I_{i^*} = \lt(\frac{\lambda_{i^*}}{Z}\rt)^{\frac{1}{1-\eta}}$. This yields
 \begin{align}
q^I_k = \left\{ \begin{array}{ll}
          q^I_{i^*} \left( \frac{\lambda_k}{\lambda_{i^*}}\right)^{\frac{1}{1-\eta}}, & \quad k=1, \ldots, i^*-1,\\
          \lt(q^I_{i^*} \rt)^{1-\eta}\frac{\lambda_k}{\lambda_{i^*}}\alpha_0^\eta , & \quad  k = i^*+1, \ldots, K. \\
         \end{array},  \right.
\label{eq:qIketa2}
\end{align}
which proves Eq.~\eqref{eq:thm-smalleta} in Item 3. It remains to identify the values of $i^*$ and $q^I_{i^*}$. By Eq.~\eqref{eq:alphaqratio2}, in order to have $q^I_{i^*} \ge \alpha_0$, it is necessary and sufficient to have
\begin{align}
\lambda_{i^*} \ge& \left( \sum_{i\in \calK}(q^I_i\vee \alpha_0)^\eta \right) \alpha_0^{1-\eta} \cr
\sk{a}{=}& \left( (K-i^* ) \alpha_0^\eta + \sum_{i=1}^{i^*}(q^I_i)^\eta \right) \alpha_0^{1-\eta} \cr
\sk{b}{=}& \left( (K-i^* ) \alpha_0^\eta + \sum_{i=1}^{i^*}\left( \frac{\lambda_i}{\lambda_{i^*}}\right)^{\frac{\eta}{1-\eta}}\lt(q^I_{i^*} \rt)^\eta \right) \alpha_0^{1-\eta} \cr
\sk{c}{\ge} & \left( (K-i^* )  + \sum_{i=1}^{i^*}\left( \frac{\lambda_i}{\lambda_{i^*}}\right)^{\frac{\eta}{1-\eta}} \right) \alpha_0 \nln
= & g(i^*),
 \label{eq:istar1}
\end{align}
where the equalities $(a)$ and $(b)$ are based on Eq.~\eqref{eq:qIketa2}, and inequality $(c)$ uses the fact that $q^I_{i^*} \ge \alpha_0$. Analogously, in order to have $q^I_{i^*} < \alpha_0$,  it is necessary and sufficient to have
\begin{align}
\lambda_{i^* + 1} <& \left( \sum_{i\in \calK}(q^I_i\vee \alpha_0)^\eta \right) \alpha_0^{1-\eta} \cr
=& \left( (K-i^* ) \alpha_0^\eta + \sum_{i=1}^{i^*}(q^I_i)^\eta \right) \alpha_0^{1-\eta} \cr
=& \left( (K-i^* ) \alpha_0^\eta + \sum_{i=1}^{i^*}\left( \frac{\lambda_i}{\lambda_{i^*}}\right)^{\frac{\eta}{1-\eta}}\lt(q^I_{i^*} \rt)^\eta \right) \alpha_0^{1-\eta} \cr
\sk{a}{<} & \left( (K-i^* -1 ) \alpha_0^\eta + \sum_{i=1}^{i^*+1}\left( \frac{\lambda_i}{\lambda_{i^*+1}}\right)^{\frac{\eta}{1-\eta}}\alpha_0^\eta \right) \alpha_0^{1-\eta}
\nln
= & g(i^*+1),
\label{eq:istar2}
\end{align}
where inequality $(a)$ is derived from the bottom line of Eq.~\eqref{eq:qIketa2}: 
$$  \lt(q^I_{i^*} \rt)^\eta = \left( \frac{q^I_{i^* + 1}}{\alpha_0^\eta}\right)^{\frac{\eta}{1-\eta}}\left( \frac{\lambda_{i^*}}{\lambda_{i^*+1}}\right)^{\frac{\eta}{1-\eta}} > \alpha_0^\eta \left( \frac{\lambda_{i^*}}{\lambda_{i^*+1}}\right)^{\frac{\eta}{1-\eta}}.$$
Eqs.~\eqref{eq:istar1} and \eqref{eq:istar2} thus give a charaterization of $i^*$ in terms of the problem primitives, and establish Eq.~\eqref{eq:gitoistar}. Note that the proceeding analysis has already shown that $i^*$ is unique and lies in $\{1, \ldots, K-1\}$, although the uniqueness can also be derived easily by noticing that $g(\cdot)$ is a non-decreasing function (since $\lt(\frac{\lambda_j}{\lambda_i}\rt)^{\frac{\eta}{1-\eta}}\geq 1$ whenever $j\leq i$) and the $\lambda_i$'s are non-increasing in $i$. 

Finally, we show that $q^I_{i^*}$ exists and is given by the unique solution to Eq.~\eqref{eq:qstarval}. Substituting the expressions for the $q^I_k$'s from Eq.~\eqref{eq:qIketa2} into Eq.~\eqref{eq:alphaqratio} leads to Eq.~\eqref{eq:qstarval}. In particular, $q^I_{i^*}$ is the solution to the following equation: 
\begin{equation}
\lt( q^I_{i^*} \rt)^{1-\eta}= \frac{\lambda_{i^*}}{(K-i^*)\alpha_0^\eta + \lt(q^I_{i^*} \rt)^\eta \sum_{j=1}^{i^*} \left( \frac{\lambda_j}{\lambda_{i^*}}\right)^{\frac{\eta}{1-\eta}}}. 
\label{eq:qisolution}
\end{equation}
To see that such a solution exists and is unique for any fixed $i^* \in \{1, \ldots, K-1\}$, note that since $\eta<1$, the left-hand side of the above equation is a strictly increasing function in $q^I_{i^*}$, which grows from $0$ to $\infty$ as $q^I_{i^*}$ varies from $0$ to $ \infty$; the right-hand side, on the other hand, is a strictly decreasing function in $q^I_{i^*}$, which, as $q^I_{i^*}$ varies from $0$ to $\infty$, decreases from $\frac{\lambda_{i^*}}{(K-i^*)\alpha_0^\eta}$ to $0$. Together, they imply that Eq.~\eqref{eq:qisolution} must admit a unique solution in $\rp$.  
This completes the proof of Proposition \ref{prop:etaless1}.
\qed

\section{Proofs for Lemmas}

\normalfont

\subsection{Proof of Lemma \ref{lem:fluidUniqu}}
\label{app:lem:fluidUniqu}
\bpf The existence and uniqueness of the fluid solution follow from Picard's existence theorem (Section 2, Chapter 1, \cite{coddington1955theory}) by verifying that the right-hand side of Eq.~\eqref{eq:drift} is uniformly Lipschitz-continuous in $q(t)$ over $\rp^\calK$, and (trivially) continuous in $t$. To show the fluid solution's continuous dependence on initial condition, note that because of the Lipschitz continuity of $p_k(\cdot)$, there exists a constant $l>0$, such that for initial conditions $x,y\in \rp^\calK$ and $t\in \rp$, we have that
 \begin{align}
 \|q(x,t) - q(y,t)\| \leq & \|x-y\| + \int_{0}^t \|( p(q(x,s))- q(x,s)) - ( p(q(y,s))- q(y,s))\|ds \nln
 \leq &  \|x-y\| + l\int_{0}^t \|q(x,s)-q(y,s)\| ds \nln
 \leq&  \|x-y\|e^{lt}, 
 \end{align}
where the last inequality follows from  Gronwall's lemma (Proposition \ref{prop:gronwall}). Therefore, $\lim_{x\to y} q(x,t) = q(y,t)$ for all $y\in \rp^\calK$. This completes the proof of the lemma. 
\qed

\subsection{Proof of Lemma \ref{lem:stochDom1}}
\label{app:lem:stochDom1}
\bpf 
We will use a simple coupling argument as follows. Fixing $k\in \calK$, the evolution of the process $U^m_{\lambda_k}(\cdot)$ corresponds to that of $Q^m_k(\cdot)$ if action $k$ were selected for all $t\in \rp$, and it follows, given the same initial condition, that $Q^m_k(t)$ is stochastically dominated by $U^m_k(t)$ for all $t\in \rp$. Since $U^m_k(\cdot)$ is positive recurrent for all $k\in \calK$, we know that  $Q^m(\cdot)$ is also positive recurrent, which in turn implies the positive recurrence of $\wol(\cdot)$, because the evolution of $C^m(\cdot)$ is derived by sampling from $\calK$ based solely on the value of $Q^m(\cdot)$. Lemma~\ref{lem:stochDom1} follows form the above-mentioned stochastic dominance, the fact that under any finite initial condition, $U^m_k(t)$ converges in distribution to $U^m_{\lambda_k}$ as $t\to \infty$, and the observation that  $U^m_{\lambda'} \preceq U^m_{\lambda}$ whenever $\lambda\geq \lambda'\geq 0$. 
\qed

\subsection{Proof of Lemma \ref{lem:disTight}}
\label{app:lem:disTight}

\bpf
To show Eq.~\eqref{eq:tightQ1}, observe that
\begin{align}
&\lim_{x\to \infty} \sup_{m \in \N } \, (1- \pi_W^m(\cqol_x) ) \nln
=& \lim_{x\to \infty} \sup_{m \in \N }  \mathbb{P}\left(\max_{k \in \calK} \qol_k(\infty) > x \right)\sk{a}{\leq} \lim_{x\to \infty} \sup_{m \in \N }  K\pb(U^m_{\lambda_1}/m>x )  \nln 
=  &  \lim_{x\to \infty} \sup_{m \in \N }  K e^{-m\lambda_1}\sum_{i=mx}^\infty \frac{(m \lambda_1)^i }{i!}  \sk{b}{\leq} \lim_{x\to \infty} \sup_{m \in \N }  K \sum_{i=mx}^\infty \lt( \frac{em \lambda_1 }{i}\rt) ^i \nln
{\leq} &  \lim_{x\to \infty} \sup_{m \in \N }  K \sum_{i=mx}^\infty \lt( \frac{em \lambda_1 }{mx}\rt) ^i \leq   \lim_{x\to \infty}K \sum_{i=x}^\infty \lt( \frac{e \lambda_1 }{x}\rt) ^i  \nln
\leq &  \lim_{x\to \infty} K2^{-(x-1)} \nln
= & 0, 
\label{eq:tightShow}
\end{align}
where step $(a)$ follows from Lemma \ref{lem:stochDom1} and the union bound, and $(b)$ from the elementary inequality $i! \geq (i/e)^i$. Eq.~\eqref{eq:tightShow} thus shows that for all $\epsilon$, $\inf_{m \in \N} \pi_W^m(\cqol_x) > 1-\epsilon$ for all sufficiently large $x$, which proves Eq.~\eqref{eq:tightQ1}. \qed

\subsection{Proof of Lemma \ref{lem:interchange1}}
\label{app:lem:interchange1}
\bpf  We first show the following uniform convergence property: 
\begin{equation}
\lim_{m \to \infty} \sup_{x\in S\cap \calQ^m} \pb\lt( \|q(x,t)- \qol(x,t)\| > \delta \rt) =0, \quad \forall \delta >0. 
\label{eq:uniformApprox}
\end{equation}
where $\qol(x,\cdot)$ denotes a process $\qol(\cdot)$  initialized with $\pb(\qol(0) = x)=1$. Suppose, for the sake of contradiction, that there exist $\delta>0$ and a sequence $\{x_i\}_{i \in \N}$, $x_i\in  S\cap \calQ^m$,  such that
\begin{equation}
\limsup_{i \to \infty} \pb\lt( \|q(x_i,t)- \qol(x_i,t)\| > \delta \rt) >0. 
\end{equation}
Because $S$ is compact, there exists a sub-sequence $\{x_{i_j}\}_{j \in \N} \subset \{x_i\}_{i \in \N}$ and $x^*\in S$ such that $x_{i_j} \to x^*$ as $j\to \infty$. We have that
\begin{align}
& \limsup_{j \to \infty} \pb\lt( \|q(x_{i_j},t)- \qol(x_{i_j},t)\| > \delta \rt)\nln
 \sk{a}{=}&   \limsup_{j \to \infty} \pb\lt( \|q(x^*,t)- \qol(x_{i_j},t)\| > \delta  \rt)  >0, 
\end{align}
where step $(a)$ follows from the fact that, for all $t\in \rp$, $q(x,t)$ is continuous with respect to $x$  (Lemma \ref{lem:fluidUniqu}). This leads to a contradiction with Theorem \ref{thm:MFE}, and hence proves Eq.~\eqref{eq:uniformApprox}.

\subsection{Proof of Lemma \ref{lem:pasta}}
\label{app:lem:pasta}
\bpf It is not difficult verify that $\{W^m[n]\}_{n\in \zp}$ is a time-homogeneous, aperiodic, and irreducible Markov chain. Because the continuous-time process,  $\{W^m(t)\}_{t\in \rp}$, is positive recurrent, so is its discrete-time counterpart, $\{W^m[n]\}_{n \in \N}$, and $W^m[n]$ converges to its steady-state distribution, $W^m[\infty]$, as $n\to \infty$.  Eq.~\eqref{eq:Qmdom-discrete} follows from the same argument as that of Lemma \ref{lem:stochDom1}. 
We now show Eq.~\eqref{eq:Cm-discrete}. Denote by $N^m(t)$ the index of the last update point by time $t$: 
\begin{equation}
N^m(t) = \sup\{n: S^m_n \leq t \}, 
\end{equation}
and by $T^m(t)$ its value
\begin{equation}
T^m(t) = S^m_{N^m(t)}. 
\end{equation}

Recall that the update points are generated according to a Poisson process, and it not difficult to show that, almost surely, $N^m(t) \to \infty$ and $T^m(t) \to \infty$, as $t\to \infty$. We thus have that, for all $k\in \calK$, 
\begin{align}
\pb(C^m[\infty] =k ) = & \lim_{n \to \infty } \pb(C^m[n] = k)  \nln
=&  \lim_{t\to \infty } \pb(C^m(T^m(t)) =k) = \pb(C^m(\infty) = k). 
\end{align}
The same argument applies for $Q^m[\infty]$ versus $Q^m(\cdot)$. This completes the proof of Lemma \ref{lem:pasta}. 
\qed

\subsection{Proof of Lemma \ref{lem:poisson}}
\label{app:lem:poisson}
\bpf 
Define
\begin{equation}
\gamma(t) = N(t) - \int_{0}^t \psi (W^m(s)) ds, \quad t\in \rp. 
\end{equation}
Let $\lt\{\calF_t\rt\}_{t \in \rp}$ be the natural filtration associated with $W^m(\cdot)$. It is not difficult to show, by the definition of $N(\cdot)$, that $\{\gamma(t)\}_{t \in \rp}$ is martingale with respect to $\{\calF_t\}_{t \in \rp}$.  Define the stopping time
\begin{equation}
T_s = \inf\{t: \psi(W^m(t))\ge \phi \}. 
\end{equation}
Let $\{\tilde N(t)\}_{t\in \N}$ be a counting process defined as: 
\begin{equation}
\tilde N (t) = N(t \wedge T_s) , \quad t\in \rp. 
\end{equation}
That is, $\tilde N(t)$ coincides with $N(t)$ up until $t=T_s$, and stays constant afterwards. Let
\begin{equation}
\tilde \psi(t) \bydef \psi(W^m(t))\mathbb{I}(t\leq T_s). 
\end{equation}
Then, it is not difficult to show that $\tilde N(\cdot)$ is a counting process whose instantaneous rate at time $t$ is $\tilde \psi(t)$, and the process 
$$\tilde \gamma (t) = \tilde N (t) - \int_{0}^t \tilde\psi (s) ds, \quad t\in \rp, $$
is a martingale with respect to $\{\calF_t\}_{t\in \rp}$.  From the definition of $\tilde \gamma(\cdot)$ and $\gamma(\cdot)$, we have that, for all $\epsilon>0$, 
\begin{align}
& \pb\left( \sup_{0\le t \le T} \left|\gamma(t) \right| \ge \epsilon \right) \nln
=&\pb\left( \sup_{0\le t \le T} \left|\gamma(t) \right| \ge \epsilon , T_s > T \right)+\pb\left( \sup_{0\le t \le T} \left|\gamma(t) \right| \ge \epsilon , T_s \le T\right) \cr
\leq & \pb\left( \sup_{0\le t \le T} \left|\tilde \gamma(t) \right| \ge \epsilon \right)+\pb\left( T_s \le T\right)\nln
=& \pb\left( \sup_{0\le t \le T} \left|\tilde \gamma (t) \right| \ge \epsilon \right) +\pb \left( \sup_{0\le t \le T}  \psi(W^m(t)) \ge \phi \right).
\end{align}
To complete the proof, therefore, it suffices to show that 
\begin{equation}
\pb\left( \sup_{0\le t \le T} \left|\tilde \gamma (t) \right| \ge \epsilon \right)  \le 2e^{-\phi T\cdot h(\epsilon / \phi T) }. \label{eq:gammaprima}
\end{equation}

We now show Eq.~\eqref{eq:gammaprima} using Doob's inequality (Proposition \ref{prop:doob} in Appendix \ref{app:technicaLemma}), by following a line of arguments similar to that used in the proof of Theorem 2.2 of \cite{kurtz1978strong}. First, we introduce a representation of the Markov process $W^m(\cdot)$ using Poisson processes. Let $\{\Xi_{y}(\cdot)\}_{y\in \mathbb{Z}^{|\calK|+1}}$ a family of mutually independent unit-rate Poisson counting processes, indexed by $\mathbb{Z}^{|\calK|+1}$. For every $y\in \Xi_{y}$, let $\psi_{y}(\cdot)$ be the rate function of $W^m(\cdot)$ for the jump with value $y$, i.e., $\psi_{y}(w)$ is the instantaneous rate at which  $W^m(\cdot)$ jumps to state $w+y$ when in state $w\in \zp^{K+1}$. Then, the process $W^m(\cdot)$ can be expressed as a solution to the following integral equation: 
\begin{equation}
W^m(t) = W^m(0) + \sum_{y\in \zp^{K+1}} y \, \Xi_y\lt(\int_{0}^t \psi_{y}(W^m(s)) ds \rt), \quad t \in \rp. 
\end{equation}
In this representation, $\Xi_y\lt(\int_{0}^t \psi_{y}(W^m(s)) ds \rt)$ counts the number of jumps of value $y$ over the interval $[0,t]$. We thus have that, by setting $y=l$, 
\begin{equation}
N(t) = \Xi_l \lt(\int_{0}^t \psi(W^m(s)) ds \rt), \quad t\in \rp, 
\end{equation}
and 
\begin{equation}
\tilde \gamma (t) =  \tilde N (t) - \int_{0}^t \tilde\psi (s) ds= \Xi_l \lt(\int_0^t \tilde\psi (s) ds \rt)-\int_0^t \tilde\psi (s) ds, \quad t \in \rp. 
\label{eq:tildeGamPoi}
\end{equation}

Fix $\theta >0$. Since  $\tilde{\gamma}(\cdot)$ is a martingale and $\exp(\cdot)$ a positive convex function, $\{\exp(\theta \tilde{\gamma}(t))\})_{t\in \rp}$ is a submartingale. From Eq.~\eqref{eq:tildeGamPoi}, we have that
\begin{equation}
 \E \lt( \exp \lt( \tilde \gamma (T) \rt) \rt)  = \E \lt(\exp\lt( \Xi_l(\tau_T) - \tau_T \rt) \rt),
\end{equation}
where $\tau_t = \int_{0}^t \tilde\psi(s)ds$. Note that $\tau_T$ is a stopping time with respect to $\Xi_l(\cdot)$. Since by definition $\tilde \psi(t) \leq \phi$ for all $t$, we have that $\tau_T \leq \phi T$. Applying the optional sampling theorem for submartingales indexed by partially ordered sets (cf.~\cite{washburn1981optional}) to $\lt\{\exp(\theta \tilde \gamma(t))\rt\}_{t\in \rp}$, we have that
 \begin{align}
  \E \lt( \exp \lt( \theta \tilde \gamma (T) \rt) \rt)    \leq \E \lt(\exp\lt(\theta \Xi_l(\phi T) - \theta \phi T \rt) \rt) \leq & \exp\lt((e^\theta-\theta-1)\phi T\rt), \label{eq:Ndoob1}
  \end{align}
where the last inequality follows from the fact that  $ \Xi_l(\phi T) $ is a Poisson random variable with mean $\phi T$ whose moment generating function is given by $\E(\exp( a \Xi_l(\phi T))) = \exp(\phi T(e^a-1))$, $a \in \R$. Analogously, we can show that
\begin{equation}
 \E \left( \exp\left( -\theta \tilde \gamma(T)  \right) \right) \leq \exp\left(\left( e^{-\theta} + \theta -1\right) \phi T \right)\leq \exp\left(\left(e^\theta -\theta -1 \right)\phi T \right), \label{eq:Ndoob2}
\end{equation}
where the last inequality follows from the  fact that $e^\theta - e^{-\theta} \ge 2\theta$ for all $\theta \ge 0$. 

Note that for any $\theta >0$, both $\left\{\exp\left( \theta \gamma(t) \right) \right\}_{t\ge 0}$ and $\left\{\exp\left( -\theta \gamma(t) \right) \right\}_{t\ge 0}$ are non-negative submartingales. Using Doob's inequality and Eqs.~\eqref{eq:Ndoob1} and Eq.~\eqref{eq:Ndoob2}, we have that, for all $\epsilon, \theta>0$, 
\begin{align}
& \pb \left( \sup_{0\le t \le T} \tilde \gamma(t)  \ge \epsilon \right) + \mathbb{P} \left( \sup_{0 \le t \le T} -\tilde \gamma(t)  \ge \epsilon \right) \nln
= & \pb\lt(\exp(\theta \tilde \gamma(T)) \geq \exp(\theta \epsilon) \rt) + \pb\lt(\exp(- \theta \tilde \gamma(T)) \geq \exp(\theta \epsilon) \rt) \nln
\le&\frac{\mathbb{E} \left( \exp\left( \theta \tilde \gamma(T) \right) \right) }{\exp(\theta \epsilon)} + \frac{\mathbb{E} \left( \exp\left( -\theta\tilde \gamma(T) \right) \right) }{\exp(\theta \epsilon)} \cr
\le& 2\exp\left(\phi T \left(e^\theta -\theta -1 \right) - \theta \epsilon \right). \label{eq:doob-bd1}
\end{align}
By setting $\theta = \log \left(1+{\epsilon}/{\phi T} \right)$ in Eq.~\eqref{eq:doob-bd1}, we conclude that
$$\mathbb{P}\left( \sup_{0\le t \le T} \left|\tilde \gamma (t) \right| \ge \epsilon \right) \le 2\exp(-\phi T\cdot h(\epsilon/\phi T) ),$$
where $h(x) \bydef (1+x)\log(1+x)-x$.  This completes the proof. \qed


\subsection{Proof of Lemma~\ref{lem:xiinduct}}
\label{app:lem:xiinduct}
\bpf We show the result by induction. For the base case, we extend the definition of $\{\xi^m_i\}_{i \in \N}$ by letting $\xi^m_0\bydef0$, and it is not difficult to see that the inequality holds when $n = 0$. 

Fix $i \in \{0, \ldots, n-1\}$, and suppose that
\begin{equation}
\E\left(\exp\left( \theta \xi^m_{i}\right) \right)  \le \exp \left(   \frac{ i \theta^2/4 }{ m\beta_m  \left( m\beta_m -\theta \right)  } \right).
\label{eq:induct0}
\end{equation}
In light of the base case, it then suffices to show that the above equation implies
\begin{equation}
\E\left(\exp\left( \theta \xi^m_{i+1}\right) \right)  \le \exp \left(   \frac{ (i+1) \theta^2/4 }{ m\beta_m  \left( m\beta_m -\theta \right)} \right).
\end{equation}

Let $\{\calC_l\}_{l \in \zp}$ be the natural filtration induced by $\{\tau_l^m,Z^m_k[l],\qol [l]\}_{l\in \N}$, with $\calC_0\bydef \emptyset$. We have that
\begin{align}
& \E\left(\exp\left( \theta \xi^m_{i+1}\right) \bbar \calC_i\right) \nln
=& \E \left( \exp\left( \theta \xi^m_{i}\right)  \left.\exp\left(  \theta\tau^m_{i+1}\left( Z^m_k[i+1] - p_k(\qol [i+1]) \right) \right) \right| \mathcal{C}_i \right) \cr
\sk{a}{=} &  \exp\left( \theta \xi^m_{i}\right)  \E \left( \left.\exp\left(  \theta\tau^m_{i+1}\left( Z^m_k[i+1] - p_k(\qol [i+1]) \right) \right) \right| \calC_i \right), 
\label{eq:induct1}
\end{align}
where step $(a)$ follows from $\xi^m_i$ being $\calC_i$-measurable. We now develop an upper bound for the second term on the right-hand side of Eq.~\eqref{eq:induct1}, as follows.
\begin{align}
&\E \left( \left.\exp\left(  \theta\tau^m_{i+1}\left( Z^m_k[i+1] - p_k(\qol [i+1]) \right) \right) \right| \calC_i \right) \cr
=& \E\left(\left. \E\left(\left. \exp\left(  \theta\tau^m_{i+1}\left( Z^m_k[i+1] - p_k(\qol [i+1]) \right) \right)\right|\qol [i+1] \right)\right| \calC_i\right) \cr
\stackrel{(a)}{=}&  \E\left(\left.\E\left(p_k(\qol[i+1]) \exp\lt(  \theta\left( 1 - p_k(\qol[i+1]) \right)\tau^m_{i+1} \rt) \rt. \rt. \rt.  \nln
& \qquad \lt. \lt. \lt.+\lt(1-p_k(\qol[i+1])\rt)\exp\lt( -\theta p_k(\qol[i+1])\tau^m_{i+1} \rt) \Bbar \qol[i+1] \right)\right|\calC_i\right)  \cr 
\stackrel{(b)}{=} & \E\left(\left. \frac{p_k(\qol[i+1]) m\beta_m}{m\beta_m -\theta \left( 1 - p_k(\qol[i+1]) \right) } +\frac{\lt(1-p_k(\qol[i+1])\rt)m\beta_m}{m\beta_m + \theta p_k(\qol[i+1]) }\right|\calC_i\right) \cr
= & \E \left( \left. \frac{m\beta_m}{m\beta_m + \theta p_k(\qol[i+1]) } \cdot \frac{m\beta_m + 2 \theta p_k(\qol[i+1])-\theta  }{m\beta_m +  \theta p_k(\qol[i+1])-\theta   }\right|\calC_i\right) \cr
= & \E \left(\left. 1+  \frac{ \theta^2 p_k(\qol[i+1])\lt(1-p_k(\qol[i+1])\rt)}{\left( m\beta_m + \theta p_k(\qol[i+1]) \right) \left( m\beta_m +  \theta p_k(\qol[i+1])-\theta \right)  }\right| \calC_i \right)  \cr
\stackrel{(c)}{\le} & \E\left(\left.\exp \left(  \frac{ \theta^2 p_k(\qol[i+1])\lt(1-p_k(\qol[i+1])\rt)}{\left( m\beta_m + \theta p_k(\qol[i+1]) \right) \left( m\beta_m +  \theta p_k(\qol[i+1])-\theta \right)  } \right)\right| \calC_i\right)  \cr
\sk{d}{\leq} & \exp \left(    \frac{ \theta^2/4 }{ m\beta_m  \left( m\beta_m -\theta \right)  } \right). 
\label{eq:induct2}
\end{align}
Step $(a)$ follows from the fact that, for a given value of $\qol[i+1]$, $Z^m_k[i+1]$ is a Bernoulli random variable with $\pb\lt(Z^m_k[i+1]=1\rt) = p_k(\qol[i+1])$, and is independent from $\tau^m_{i+1}$. For step $(b)$, note that $\tau_{i+1}$ is an exponential random variable with mean $(m\beta_m)^{-1}$, independent from $\calC_i$ and $\qol[i+1]$. Its moment generating function is given by $\mathbb{E}\left( e^{h\tau_{i+1}} \right) = \frac{\beta_m}{\beta_m- h}$ for all $h\le \beta_m$, where $h$, in our case, corresponds to $\theta(1 - p_k(\qol[i+1])) $ and $\theta p_k(\qol[i+1])$, for the two terms respectively. Step $(c)$ stems from the fact that $1+x \le e^x $ for all $x\in \rp$. Finally, for step $(d)$ we have used the fact that $ p_k(\qol[i+1]) \in [0,1]$ by definition, and that $\theta < m\beta_m$; the exponent is hence bounded from above by setting $p_k(\qol[i+1])$ to $1/2$ and $0$ in the numerator and denominator, respectively. 

Substituting Eq.~\eqref{eq:induct2} into Eq.~\eqref{eq:induct1}, and invoking the induction hypothesis of Eq.~\eqref{eq:induct0}, we have that
\begin{align}
\E\left(\exp\left( \theta \xi^m_{i+1}\right) \right) 
= & \E\lt( \E\left(\exp\left( \theta \xi^m_{i+1}\right) \bbar \calC_i \right) \rt) \nln
=  & \E\lt(  \exp\left( \theta \xi^m_{i}\right) \E \left( \left.\exp\left(  \theta\tau^m_{i+1}\left( Z^m_k[i+1] - p_k(\qol [i+1]) \right) \right) \right| \calC_i \right)\right)  \nln
\leq & \exp \left(   \frac{ i \theta^2/4 }{ m\beta_m  \left( m\beta_m -\theta \right)  } \right) \exp \left(    \frac{ \theta^2/4 }{ m\beta_m  \left( m\beta_m -\theta \right)  } \right) \nln
\leq & \exp \left(   \frac{ (i+1) \theta^2/4 }{ m\beta_m  \left( m\beta_m -\theta \right)  } \right). 
\end{align}
This proves our claim.   \qed

\subsection{Proof of Lemma \ref{lem:xisqrTailBound}}
\label{app:lem:xisqrTailBound}

\bpf Fix $\epsilon>0$ and $m\in \N$. We have that
\begin{align}
&  \pb\lt( \sum_{i=1}^{I_{mT}+1} (\tau^m_i)^2 \geq \epsilon \rt)  \nln
 \leq& \pb\lt( \sum_{i=1}^{I_{mT}+1} (\tau^m_i)^2 \geq \epsilon \, , \, I_{mT} < (1+\epsilon)m\beta_mT  \rt) +\pb\lt( I_{mT} \geq (1+\epsilon) m\beta_mT \rt)  \nln
\sk{a}{\leq} & \pb\lt( \sum_{i=1}^{(1+\epsilon)m\beta_mT} (\tau^m_i)^2 \geq \epsilon \rt) + \frac{1+\epsilon}{T\epsilon^2}(m\beta_m )^{-1} \nln
\sk{b}{\leq} & \frac{ (1+\epsilon )m\beta_mT \cdot \E((\tau^m_1)^2)}{\epsilon} + \frac{1+\epsilon}{T\epsilon^2}(m\beta_m )^{-1} \nln
\sk{c}{=}  & { 2(1+\epsilon ) T} (\epsilon   m\beta_m)^{-1}+ \frac{1+\epsilon}{T\epsilon^2}(m\beta_m )^{-1}, 
\end{align} 
where step $(a)$ follows from Eq.~\eqref{eq:pbcalAbound}, $(b)$ from the Markov's inequality, and $(c)$ from $\tau^m_i$ being an exponentially distributed random variable with mean $(m \beta_m)^{-1}$ and hence $\E((\tau^m_1)^2) = 2(m\beta_m)^{-2}$. Because $m \beta_m \to \infty$ as $m\to \infty$, the claim follows. \qed

\end{document}